\documentclass[11pt]{amsart}

\usepackage[a4paper,margin=1in]{geometry} 
\usepackage{setspace}
\usepackage[T1]{fontenc}
\usepackage[utf8]{inputenc}

\usepackage{amsmath,amssymb,amsfonts,amsthm,bm}

\usepackage{graphicx}
\usepackage{booktabs}
\usepackage{multirow}

\usepackage{indentfirst}

\usepackage{subcaption}

\usepackage[authoryear]{natbib}

\usepackage[hidelinks]{hyperref}

\usepackage{tikz}
\usetikzlibrary{shapes.geometric, arrows, positioning, fit, calc}
\usepackage{boldline}
\usepackage{multirow}

\usepackage{ctable}
\usepackage{float}

\usepackage{varioref}
\usepackage{xr-hyper}
\usepackage{hyperref}
\hypersetup{colorlinks,linkcolor=blue,citecolor=blue,filecolor=red,urlcolor=blue}

\usepackage{algorithm}
\usepackage{algpseudocode} 

\usepackage{pgfplots}
\pgfplotsset{compat=1.18}

\def\P{\mathbb{P}}
\def\F{\mathbb{F}}
\def\E{\mathbb{E}}
\def\H{\mathrm{H}}
\def\T{\mathrm{T}}
\def\calT{\mathcal{T}}
\def\calN{\mathcal{N}}
\def\calF{\mathcal{F}}
\def\calB{\mathcal{B}}
\def\N{\overline{N}}
\def\bbN{\mathbb{N}}
\def\K{\overline{K}}
\def\L{\mathbb{L}}
\def\R{\mathbb{R}}
\def\M{\mathbb{M}}
\def\Var{\mathrm{Var}}

\definecolor{myblue}{rgb}{0.2,0,0.9}
\definecolor{blue_violet}{rgb}{0.54, 0.17, 0.89}
\definecolor{darkgreen}{rgb}{0,0.35,0}

\allowdisplaybreaks

\newtheorem{theorem}{Theorem}[section]
\newtheorem{proposition}[theorem]{Proposition}
\newtheorem{problem}[theorem]{Problem}
\newtheorem{corollary}[theorem]{Corollary}
\newtheorem{lemma}[theorem]{Lemma}
\numberwithin{equation}{section}

\theoremstyle{definition}
\newtheorem{remark}[theorem]{Remark}
\newtheorem{definition}[theorem]{Definition}
\newtheorem{assumption}[theorem]{Assumption}

\DeclareMathOperator*{\Lip}{Lip}
\DeclareMathOperator*{\LipH}{LipH}
\DeclareMathOperator*{\size}{Size}
\DeclareMathOperator*{\Growth}{Gr}
\DeclareMathOperator*{\Holder}{Hol}
\DeclareMathOperator*{\Tr}{Tr}
\DeclareMathOperator{\Transpose}{T}
\DeclareMathOperator*{\arginf}{arg\,inf}
\DeclareMathOperator*{\argmin}{arg\,min}

\DeclareMathOperator*{\esssup}{ess\,sup}

\title[DeepMartingale]{DeepMartingale: Duality of the Optimal Stopping Problem with Expressivity and high-dimensional hedging}

\author{Junyan Ye and Hoi Ying Wong}
\thanks{Department of Statistics and Data Science, The Chinese University of Hong Kong, Shatin, N.T., Hong Kong. Emails: \texttt{junyanye@link.cuhk.edu.hk}, \texttt{hywong@cuhk.edu.hk}.}

\date{}

\begin{document}

\maketitle

\begin{abstract}

We propose \textit{DeepMartingale}, a deep-learning framework for the dual formulation of discrete-monitoring optimal stopping problems under continuous-time models. Leveraging a martingale representation, our method implements a \emph{pure-dual} procedure that directly optimizes over a parameterized class of martingales, producing computable and tight \emph{dual upper bounds} for the value function in high-dimensional settings without requiring any primal information or Snell-envelope approximation. We prove convergence of the resulting upper bounds under mild assumptions for both first- and second-moment losses. A key contribution is an expressivity theorem showing that \textit{DeepMartingale} can approximate the true value function to any prescribed accuracy $\varepsilon$ using neural networks of size at most $\tilde{c} d^{\tilde{q}}\varepsilon^{-\tilde{r}}$, with constants independent of the dimension $d$ and accuracy $\varepsilon$, thereby avoiding the curse of dimensionality. Since expressivity in this setting translates into scalability, our theory also motivates estimating the dimension scaling law to guide architecture design and the training setup in deep learning-based numerical computation and the choice of rebalancing frequency for the related hedging strategy. The learned martingale representation further yields a practical and dimension-scalable \emph{deep delta hedging strategy}. Numerical experiments on high-dimensional Bermudan option benchmarks confirm convergence, expressivity, scalable training, and the stability of the resulting upper bounds and hedging performance.

\vspace{0.5em}


\end{abstract}


\section{Introduction}\label{sec:Intro}

Optimal stopping problems are commonly studied from two complementary viewpoints: \emph{primal} and \emph{dual}. For maximization problems, primal methods search over stopping rules and typically deliver lower bound of the value, whereas dual methods optimize over martingales to produce upper bounds and an associated hedging strategy.

On the primal side, classical simulation-based algorithms include least-squares Monte Carlo methods \citep{carriere96,longstaff01,vanroy01,Clement2002} and variants combined with policy iteration \citep{bender08}. A well-known limitation is the dependence on carefully chosen basis functions, whose complexity can grow rapidly with the dimension, leading to instability in high-dimensional settings \citep{CSW2019}. On the dual side, upper bounds for the Snell envelope can be obtained by minimizing over suitable classes of martingales \citep{Haugh04,Roger02}. Early implementations often relied on nested Monte Carlo \citep{andersen04,kolodko04}, while later work developed faster non-nested alternatives \citep{belome09,brown10}. The pure-dual paradigm of \citep{roger10} (see also \citep{puredual-mf}) is particularly attractive since it does not require an accurate approximation of the Snell envelope and can yield hedging information without bias from primal estimates; however, existing dual simulation methods still face substantial challenges in high dimensions.

The strong empirical performance of deep neural networks (DNNs) has stimulated their application to a broad range of problems in finance, including optimal stopping. DNN-based numerical PDE solvers have been developed in \citep{han18,RAISSI2019}, and expressivity-based analyses showing how certain PDE classes can avoid the curse of dimensionality have been established in \citep{Jentzen20,grohs21,Jentzen23}. 
For primal optimal stopping, \cite{KohlerNN2010} pioneered the use of neural networks within least-squares regression for continuation values. \citep{Becker19} then used neural networks to learn approximate stopping policies in a semimartingale setting, and \citep{Becker20} proposed direct approximation of the primal value function (equivalently, the continuation value), with discrete-time expressivity guarantees provided in \citep{gonon23}. However, many models in finance are specified in continuous time, while exercise opportunities are monitored at discrete time points. Related approaches include \cite{lapeyre2020neuralnetworkregressionbermudan}, which follows the LS-style methodology of \cite{longstaff01} and dispenses with the dynamic programming principle, and \citep{reppen23}, which studies direct neural approximations of the free boundary in a continuous-time framework but requires a prescribed boundary.


Despite this progress, two gaps remain central for high-dimensional applications: (i). existing deep-learning approaches that combine primal and dual perspectives \citep{guo2024simultaneousupperlowerbounds,yang2024deepprimaldualbsdemethod} do not, in general, provide a dual-side expressivity theory or numerical guarantees that truly overcome the curse of dimensionality; in particular, they do not deliver a theoretically grounded \emph{dimension scaling law} to guide architecture design, the training setup and the choice of the rebalancing frequency for the related hedging strategy as $d$ grows. In our implementation, these methods fail to produce reliable results in high-dimensional settings under continuous-time models with discrete-time monitoring; 
(ii). reliable, dimension-scalable, and accurate \emph{high-dimensional hedging} strategies for Bermudan-type options remain largely unaddressed. The deterioration observed for \citep{guo2024simultaneousupperlowerbounds,yang2024deepprimaldualbsdemethod} in high dimensions further highlights the need for a dimension-scalable and accurate hedging approach.

Motivated by these gaps, we propose \textit{DeepMartingale}, a pure-dual deep-learning framework for optimal stopping under continuous-time It\^{o} diffusions with discrete monitoring. The method directly optimizes over a parameterized class of martingales to compute a dual upper bound, and the resulting learned martingale representation naturally yields a practical, dimension-scalable deep delta hedging strategy. Leveraging martingale theory and our network design, we establish a complete expressivity theory that guarantees the dimension-scalability of \textit{DeepMartingale} (i.e., avoidance of the curse of dimensionality) and supports a corresponding numerical estimation of the \emph{dimension scaling law}, thereby guiding architecture design and training setup for \textit{DeepMartingale}, as well as the choice of rebalancing frequency for the related \textit{deep delta hedging}.

\subsection{Key contributions of this paper}
\begin{enumerate}
\item \textbf{A pure-dual \textit{DeepMartingale} framework for continuous-time models with discrete monitoring.}
Section~\ref{sec:deep_mtg} introduces \textit{DeepMartingale}. We prove convergence when minimizing either a first-moment (upper-bound) loss or a second-moment loss (Corollary~\ref{coro:tight_upper} and Proposition~\ref{pro:L2-to-L1-loss}). This formulation is well suited for Bermudan-style contracts, hedging, and operational stopping problems, where the state evolves in continuous time while decisions occur at discrete monitoring dates.
\item \textbf{Dual expressivity theory with dimension-dependent complexity bounds (including stochastic-integration discretization), and dimension-scaling law estimation for architecture and training design.}
We develop a comprehensive \emph{expressivity} theory for the dual (martingale) problem, addressing two components:
(i) for It\^{o} diffusions whose coefficients satisfy suitable expressivity conditions, Section~\ref{subsec:numerical_express} establishes expressivity results for the stochastic-integral approximation in the martingale representation (Theorem~\ref{theorem:express_N_0_new});
(ii) building on \cite{Jentzen23,gonon23} and an infinite-width random-feature / reproducing kernel Banach space (RKBS) perspective \citep{BARTOLUCCI2023194,RKBS2024}, Section~\ref{subsec:expressivity} derives dual expressivity results under structural assumptions (Theorem~\ref{thm:express_deep_mtg}), including affine It\^{o} diffusions as an important special case (Theorem~\ref{thm:AID-log_express}).
Together, these results establish the dimension-scalability (avoidance of the curse of dimensionality) of \textit{DeepMartingale} and motivate the dimension-scaling law estimation in Section~\ref{subsubsec:dimension-scaling-law}, which infers the implied expressivity order in the dimension $d$ from low-dimensional experiments and guides architecture design, training setup, and the choice of rebalancing frequency for the hedging strategy in high dimensions.
\item \textbf{High-dimensional numerical reliability: scalable training, upper bounds, and dimension-scalable deep delta hedging.}
Section~\ref{subsec:experiments} presents numerical results on high-dimensional Bermudan option benchmarks, demonstrating:
(i) accurate, reliable, and computable \emph{dual upper bounds} in moderate to high dimensions under scaled computational budgets, compared with recent methods \cite{guo2024simultaneousupperlowerbounds,yang2024deepprimaldualbsdemethod,puredual-mf} (Section~\ref{subsubsec:experiment_compare_after_scaling});
(ii) effective and dimension-scalable \emph{deep delta hedging} induced by the learned martingale (Sections~\ref{subsubsec:numerical_delta_hedge_evaluation} and \ref{subsubsec:experiment_hedge}), evaluated via a \emph{worst-case} hedging error metric consistent with adversarial exercise decisions.
Under comparable computational budgets, our \textit{deep delta hedging} approach delivers stable, reliable, and computable hedging performance in moderate to high dimensions relative to \cite{guo2024simultaneousupperlowerbounds,yang2024deepprimaldualbsdemethod,puredual-mf}.
\end{enumerate}


\subsection{Notations}\label{sec:notations}

Fix a finite time horizon \(T>0\). Let \((\Omega,\mathcal{F},\mathbb{F},\mathbb{P})\) be a filtered probability space with continuous-time filtration \(\mathbb{F}:=(\mathcal{F}_t)_{t\in[0,T]}\). For each \(t\in[0,T]\), we denote by
\(
\mathbb{F}_t:=(\mathcal{F}_s)_{s\in[t,T]}
\)
the filtration restricted to the time interval \([t,T]\).

For the discrete monitoring setup, fix the number of monitoring rights \(N\in\mathbb{N}_+\). Define $ \overline{N} := \{ 0,\ldots,N \} $, $ \overline{N}^{-1} := \overline{N} \setminus \{ N \} $ , and for \(n \in \overline{N}\), set \( \overline{N}_n:=\{ n,\ldots,N \} \), $ \overline{N}_n^{-1} := \overline{N}_n \setminus \{ N \} $. The monitoring points are specified by the uniform grid \( \T^N := \{t_n: t_n=nT/N,\; n\in \overline{N} \} \), and we set \( \T^N_n := \{t_m: t_m=mT/N,\; m\in \overline{N}_n \} \). Denote the corresponding discrete-monitoring filtrations by
\(
\mathbb{F}^N:=(\mathcal{F}_{t_m})_{m=0}^{N},\;
\mathbb{F}^N_n:=(\mathcal{F}_{t_m})_{m=n}^{N},\; n\in \overline{N}.
\)

We write \(\mathbb{E}[\cdot]\) for expectation under $ \mathbb{P} $ and $ \Var(\cdot) $ for the variance. Moreover,
\(\mathbb{E}_n[\cdot]:=\mathbb{E}[\cdot\mid\mathcal{F}_{t_n}]\)
denotes conditional expectation and \( \Var_n(\cdot) \) denotes conditional variance on $ \calF_{t_n} $. For any random variable $X$, we write $ \E^X[\cdot ] $ for conditional expectation given $X$.
For vectors in Euclidean space, \(\|\cdot\|\) denotes the Euclidean norm; for real matrices, \(\|\cdot\|_{\H}\) denotes the Hilbert--Schmidt norm and $ \| \cdot \|_p $ denote the $ p $-induced norm (operator norm). $ \calB(\R^d) $ denote the Borel $\sigma$-field on $ \R^d $.
We adopt the conventions \(\inf\varnothing:=+\infty\) and \(\sum_{k\in\varnothing}:=0\).

For any \(p\ge 1\), \(k\in\mathbb{N}\), and \(t\in[0,T]\), we introduce the following sets.
\begin{itemize}
\item \(L^p(\mathcal{F}_t;\mathbb{R}^k)\) denotes the space of \(\mathbb{R}^k\)-valued, \(\mathcal{F}_t\)-measurable random variables \(\xi\) such that
\(
\|\xi\|_{L^p}^p:=\mathbb{E}\|\xi\|^p < \infty.
\)
If it is $ \R^{k\times k} $-valued, then \(
\|\xi\|_{L^p;\H}^p:=\mathbb{E}\|\xi\|^p_{\H}
\).
\item \(\mathbb{L}^p(\mathbb{R}^k)\) denotes the space of \(\mathbb{R}^k\)-valued, \(\mathbb{F}\)-adapted processes \(Z=(Z_t)_{t\in[0,T]}\) such that
\(
\|Z\|_{\mathbb{L}^p}^p:=\mathbb{E}\int_0^T \|Z_t\|^p\,dt<\infty.
\)
For \(0\le s\le t\le T\), \(\mathbb{L}^p_{s,t}(\mathbb{R}^k)\) is defined by
\(
\|Z\|_{[s,t],p}^p:=\mathbb{E}\int_s^t \|Z_u\|^p\,du<\infty.
\)
\item \(\mathbb{L}_N^p(\mathbb{R}^k)\) denotes the space of \(\mathbb{R}^k\)-valued, \(\mathbb{F}^N\)-adapted discrete-time processes \(X=(X_{t_n})_{n=0}^N\) such that \(X_{t_n}\in L^p(\mathcal{F}_{t_n};\mathbb{R}^k)\) for all \(n\in \overline{N} \).
Similarly, \(\mathbb{L}_{n,N}^p(\mathbb{R}^k)\) denotes the corresponding class of processes restricted to \(\T^N_n\).
\item Given a Borel measure $ \rho $ on $ \mathbb{R}^{k_1} $, \(L^p_{k_1,k_2}(\rho)\) denotes the space of functions $F: \mathbb{R}^{k_1} \to \mathbb{R}^{k_2} $ such that $\| F\|^p_{p;\rho} := \int_{\mathbb{R}^{k_1} } \|F(x) \|^p \rho(dx) < \infty $. Also, $ \mathbb{M}^p_p(\rho) := \int_{\mathbb{R}^k } \|x \|^p \rho(dx) $.
\item \(\mathcal{T}^N\) denotes the set of all discrete \(\mathbb{F}^N\)-stopping times \(\tau\) taking values in \( \T^N \). Similarly, \(\mathcal{T}^N_n\) denotes the set of all discrete \(\mathbb{F}^N_n\)-stopping times \(\tau\) taking values in \(\T^N_n\).
\item \(\mathcal{M}^N\) denotes the set of all discrete-time martingales \(M \in \mathbb{L}_N^1(\mathbb{R})\) (with respect to \(\mathbb{F}^N\)).
Similarly, \(\mathcal{M}^N_{n}\) denotes the class of such martingales on \(\T^N_n\).
\end{itemize}

Let $d \in \bbN_{+} $ denote the dimension. For most arguments, we omit the dependence on $d$ in the notation for simplicity. For expressivity-related statements, we explicitly indicate the dependence on $d$ using subscripts/superscripts, e.g.\ from $X$ to $X^d$. To avoid technical distractions, we assume $d\ge 3$ in the expressivity analysis; this restriction does not induce any essential difference.

\section{Problem formulation and preliminary analysis of the dual form}
\label{sec:problem_formulate}
In this section, we formulate the discrete-monitoring, continuous-observation optimal stopping problem and provide preliminary results that motivate our dual approach. In particular, we recall weak/strong duality, the surely optimal Doob martingale, and the backward recursion representation of the dual upper bound.

Let $(X_t)_{0\le t\le T}$ be a continuous-time Markovian process defined on a filtered probability space $(\Omega,\mathcal{F},\mathbb{F},\mathbb{P})$.
Given a discounted payoff function $g(t,x)$ that is Lipschitz continuous with respect to $x$, and satisfying \( \mathbb{E}|g(t_n,X_{t_n})|^2 < \infty \) for all $ n\in \overline{N} $, the optimal stopping problem is to maximize
\begin{equation}\label{eq:primal_problem}
    Y^{*}_0= \sup_{ \tau \in \mathcal{T}^N}\mathbb{E}\bigl[g(\tau,X_{\tau})\bigr]. \tag{\textbf{P}}
\end{equation}
The corresponding Snell envelope is defined by, for $ n\in \overline{N} $,
\[
    Y^{*}_n := \esssup_{\tau_n \in \calT_N^n }\mathbb{E}_n[g(\tau_n,X_{\tau_n})]  .
\]
Although our primary interest is the discrete-monitoring setting, continuous-monitoring can be approximated by increasing the monitoring frequency, without changing the underlying continuous-time dynamics. Nevertheless, much of the deep optimal stopping literature is formulated in discrete time, and thus it is natural to align discrete observation times with monitoring times.

\subsection{Duality and the Doob martingale}
Following the dual formulation in \citep{Roger02}, \citep{Haugh04}, and \citep{belome09}, we recall the following duality results.

\begin{lemma}[Duality]\label{lemma:duality}
For any $ n\in \overline{N} $, $M \in \mathcal{M}^N_n $, we have the following:
\begin{enumerate}
    \item[(i)] (Weak Duality)
    \begin{equation}
        Y^{*}_n \le \mathbb{E}_n \Bigl[\max_{n \le m \le N}\bigl( g(t_m,X_{t_m} ) - M_{t_m} + M_{t_n}\bigr)\Bigr]
        \label{eq_weak_dual}
    \end{equation}
    \item[(ii)] (Strong Duality \& Dual Problem)
    \begin{equation}
        Y_n^{*}=\inf_{M\in \mathcal{M}^N_n} \mathbb{ E}_n \Bigl[\max_{n\le m \le N} \bigl(g(t_m,X_{t_m} ) - M_{t_m} + M_{t_n}\bigr) \Bigr] .  \tag{\textbf{D}}
        \label{eq_strong_dual}
    \end{equation}
\end{enumerate}
\end{lemma}

By the Doob decomposition for $ Y^{*}_n $, there exists a Doob martingale $ M^{*} $ and an increasing predictable process $ A^{*} $ such that
\begin{equation*}
    Y_n^{*}=Y_0^{*}+M_{t_n}^{*}-A_{t_n}^{*} , \quad n\in \overline{N} .
\end{equation*}
By simple verification, we have $ Y^{*},M^{*} \in \L^2_N(\mathbb{R}) $.
The following lemma shows that the Doob martingale $M^{*} $ is a surely optimal candidate for \eqref{eq_strong_dual}.

\begin{lemma}[Surely Optimal \cite{schoen13}] \label{lemma:sure_optimal} For all $ n\in \N $
    \begin{equation}
        Y_n^{*}=\max_{n \le m \le N}\bigl(g(t_m,X_{t_m}) - M_{t_m}^{*} + M_{t_n}^{*} \bigr) \quad \mathbb{P}\text{-a.s.}
        \label{eq_sure_optimal}
    \end{equation}
\end{lemma}

\subsection{Backward recursion}
We use the backward recursion formulation in \citep{schoen13}, which, similarly to \citep{Becker19}, enables a step-by-step characterization from the terminal time.

We define a sequence of functions $\widetilde{U}_n : \mathcal{M}^N_n \rightarrow \mathbb{R}$, $n \in \N $ such that $\widetilde{U}_N(\cdot) \equiv g(t_N,X_{t_N}) $ and, for $  n \in \N^{-1} $,
\begin{equation}\label{eq:Upper_bound_M}
    \widetilde{U}_n(M) = \max_{n \le m \le N}\bigl(g(t_m,X_{t_m}) - M_{t_m} + M_{t_n}\bigr) ,
\end{equation}
for any $M \in \mathcal{M}^N_n $. The quantity $\widetilde{U}_n(M)$ is an upper bound for $Y^{*}_n$ under $ \mathbb{E}_n[\cdot] $ \citep{schoen13}, and in particular $ Y^{*}_n = \widetilde{U}_n(M^{*}) = \mathbb{E}_n[\widetilde{U}_n(M^{*})] $.

Moreover, for all $ M \in \mathcal{M}^N_n $, \( \mathrm{Var}_n ( \widetilde{U}_n(M^{*}) ) = 0 \le \mathrm{Var}_n ( \widetilde{U}_n(M) )  \),
and when $ Y^{*}_n \ge 0 $, $ \mathbb{P} $-a.s.,
\[
    \mathbb{E}\big| Y^{*}_n \big|^2 = \mathbb{E}\big|\mathbb{E}_n[Y^{*}_n] \big|^2 \le \mathbb{E}\big|\mathbb{E}_n[\widetilde{U}_n(M)] \big|^2 \le \mathbb{E}\big|\mathbb{E}_n[\widetilde{U}_n(M)] \big|^2 + \mathbb{E}\big[ \mathrm{Var}_n(\widetilde{U}_n(M)) \big] = \mathbb{E}\big|\widetilde{U}_n(M)\big|^2
\]
which shows that $ \mathbb{E}\big|\widetilde{U}_n(M)\big|^2 $ is an upper bound for $ \mathbb{E}\big| Y^{*}_n \big|^2 $, and that $ \mathbb{E}\big|\widetilde{U}_n(M^{*})\big|^2 $ attains the equality.

\begin{remark}
    The condition $ Y^{*}_n \ge 0 $ for all $ n \in \N $, $ \mathbb{P} $-a.s. is common in applications. In particular, it is guaranteed whenever the payoff $ g $ is non-negative.
\end{remark}

Let $\xi_n (M) := M_{t_{n+1}}-M_{t_n} \in \mathcal{F}_{t_{n+1}} $ be the martingale increments. Then the following backward recursion holds true \citep{schoen13}. For $n \in \N^{-1} $,
\begin{equation}\label{eq:backward_recursion}
    \widetilde{U}_n(M) = g(t_n,X_{t_n}) + \bigl(\widetilde{U}_{n+1}(M) - \xi_{n}(M)-g(t_n,X_{t_n}) \bigr)^{+}.
\end{equation}
By direct manipulation, we obtain the following basic error propagation bounds associated with \eqref{eq:Upper_bound_M}.

\begin{lemma}[Error Propagation] \label{lem-error-propagate}
    For any $ n \in \N^{-1} $, $M_1,M_2 \in \mathcal{M}^N_n $,
    \begin{equation}\label{eq_error1}
        \bigl|\widetilde{U}_n(M_1) -\widetilde{U}_n(M_2)\bigr| \le \bigl|\widetilde{U}_{n+1}(M_1) - \widetilde{U}_{n+1}(M_2) \bigr| + \bigl|\xi_{n}(M_1)- \xi_{n}(M_2) \bigr| ,
    \end{equation}
    \begin{equation}\label{eq_error2}
        \begin{split}
         \big\| \widetilde{U}_n(M_1) -\widetilde{U}_n(M_2) \big\|_{L^2}
         &\le \big\|\widetilde{U}_{n+1}(M_1) -\widetilde{U}_{n+1}(M_2)\big\|_{L^2} + \big\|\xi_n(M_1) -\xi_n(M_2)\big\|_{L^2} .
        \end{split}
    \end{equation}
\end{lemma}

\section{Numerical approximation for the Doob martingale}
\label{sec:Numerical_Approximation}
This section develops the theory that underpins our approximation of the Doob martingale under a Brownian filtration. We first recall the martingale representation and then introduce a numerical integration scheme on each interval $[t_n,t_{n+1}]$ for $n \in \N^{-1}$. We subsequently establish convergence and derive expressivity guarantees for the proposed numerical approximation.

We focus on an It\^o diffusion. The filtered probability space $(\Omega, \mathcal{F}, \mathbb{F}, \mathbb{P} )$ supports a $d$-dimensional Brownian motion $W=(W^1,\ldots,W^{d})^{\top}$, where the augmented filtration $\mathbb{F}$ is now generated by $W$. We consider the $\mathbb{R}^d$-valued solution $X$ to the SDE
\begin{equation}
    dX_t = a(t,X_t)\,dt + b(t,X_t)\,dW_t, \quad X_0 =x_0 \in \mathbb{R}^d ,  \label{eq:SDE}
\end{equation}
where $a(t,x)$ and $b(t,x)$ are Lipschitz continuous in $x$ and $\frac{1}{2}$-H\"{o}lder continuous in $t$. In this setting, by the Markov property, $Y^{*}_{t_n} = V_n(X_{t_n})$ for some $\mathcal{B}(\mathbb{R}^d)$-measurable function $V_n$.

\subsection{Martingale representation and numerical integration scheme}
Since $M^{*} \in \L^2_N(\mathbb{R})$, the martingale representation theorem implies that there exists a process $Z^{*} \in \L^2(\mathbb{R})$ such that
\begin{equation}\label{eq:Martingale_represent}
    M_{t_n}^{*} =\int_0^{t_n} Z^{*}_{s} \cdot d W_s\, , \quad n \in \N .
\end{equation}

\begin{remark}[Connection with \textit{delta hedge} strategy]\label{rem:connection_delta_hedging}
As is well known (see, e.g., \citep{belome09}), the Doob integrand \(Z^{*}\) is naturally linked to the delta hedge for Bermudan options. In a Markovian market (so that \(Z^{*}_s=Z^{*}(s,X_s)\)), if \(b(t,x)\) is invertible, then for \(t\in[t_n,t_{n+1})\) the strategy
\[
\vartheta(t,x):=b^{-1}(t,x)Z^{*}(t,x)
\]
is the delta hedge for the European claim starting at \(t_n\) with terminal payoff \(Y^{*}_{n+1}\), and it replicates this claim on \([t_n,t_{n+1})\). If the market is further complete, then $ \vartheta $ is a perfect hedge for Bermudan options. Pasting these hedges over all monitoring intervals and closing the position when the option is exercised at some \(t_n\) (paying \(g(t_n,X_{t_n})\) to the buyer) yields a strategy that attains \(Y^{*}_0\) under optimal exercise, hence perfectly replicates the Bermudan option. In Section~\ref{sec:deep_mtg}, we construct a “deep” version of \(Z^{*}\), which serves as a “deep” delta hedge.
\end{remark}

Motivated by \citep{belome09}, we approximate \eqref{eq:Martingale_represent} via a local numerical integration scheme. For each interval $[t_n,t_{n+1}]$, $n \in \N^{-1}$, we introduce $K \in \mathbb{N}_{+}$ uniform subintervals of length $\Delta t := t^n_{k+1}-t^n_{k} = T/(N K)$. The mesh points are given by
$ \T^{n,K} := \{ k\Delta t , \; k\in \K \cup \{ K \} \} $
and $ \K := \{ 0,\ldots,K-1 \} $.
The Brownian increments are denoted by $\Delta W_{t^n_k} = W_{t^{n}_{k+1}} - W_{t^{n}_{k}}$, $k\in \K$.
Let $ \E_{n,k}[\cdot] := \mathbb{E}[ \cdot \,|\,\mathcal{F}_{t^n_k} ] $, and define
\begin{equation}
    \widehat{Z}^{*,K}_{t^n_k}  := \frac{1}{\Delta t } \mathbb{E}_{n,k} \big[Y^{*}_{n+1} \Delta W_{t^n_k} \big] , \; k \in \K , \; n\in \N^{-1} , \; \; \widehat{M}^{*,K}_{t_n}  := \sum_{m=0}^{n-1} \sum_{k=0}^{K-1} \widehat{Z}^{*,K}_{t^m_k} \cdot \Delta W_{t^m_k} ,\; n\in \N ,
    \label{eq_Z}
\end{equation}
with $ \widehat{Z}^{*,K,n}_t := \sum_{k=0}^{K-1} \widehat{Z}^{*,K}_{t^n_k} 1_{[t^n_k , t^n_{k+1} )}(t) $ for $ t \in [t_n,t_{n+1} ] $.
Since $\widehat{M}^{*,K} \in \mathcal{M}^N$, we define the corresponding approximate Snell envelopes by
\begin{equation*}
    \widehat{Y}^{*,K}_n := \mathbb{E}_n \bigl[ \widetilde{U}_n(\widehat{M}^{*,K})  \bigr], \quad   n \in \N^{-1} .
\end{equation*}
Therefore, $\widehat{Y}^{*,K}$ is an upper bound for $Y^{*}$ by the weak duality \eqref{eq_weak_dual}.

The convergence of this numerical approximation is ensured by the following theorem.

\begin{theorem}[Numerical Approximation Convergence \citep{belome09}] \label{theorem:discrete_convergence}
As $K \rightarrow \infty$, for all $ n\in \N^{-1} $
    \begin{equation}
         \big\| \widehat{Z}^{*,K,n} - Z^{*} \big\|_{[t_n,t_{n+1}],2}
              , \; \mathbb{E}\Bigl[\max_{0\le n \le N}|M^{*}_{t_{n}}-\widehat{M}^{*,K}_{t_{n}} |^2 \Bigr] , \; \mathbb{E}\bigl|Y^{*}_n - \widehat{Y}^{*,K}_n \bigr|^2 \to 0 .
            \label{eq_Z_converge}
    \end{equation}
\end{theorem}

Suppose we can construct, backward in time for $n=N-1,N-2,\ldots,0$, a sequence $\tilde{\xi}_n \in \mathcal{P}_n := \{\xi \in \mathcal{F}_{t_{n+1}} : \mathbb{E}_n[\xi ]=0\}$ that approximates $\xi_n(\widehat{M}^{*,K})$ arbitrarily well in $L^2$. Defining \( \widetilde{M}_n = \sum_{m=0}^{n-1} \tilde{\xi}_m \) for all $ n \in \N^{-1} $,
we can then approximate $\widehat{Y}^{*,K}_n$ by $\widetilde{U}_n(\widetilde{M})$ via Lemma \ref{lem-error-propagate}, and consequently approximate $Y^{*}_n$ by Theorem \ref{theorem:discrete_convergence}. By weak duality \eqref{eq_weak_dual}, $ \mathbb{E}_n [\widetilde{U}_n(\widetilde{M}) ]$ remains an upper bound for $Y^{*}_n$.

To clarify the dependence in the optimization problem, we also write
$ \widetilde{U}_{n}(\tilde{\xi}_n ; \widetilde{M}^{n+1} ) := \widetilde{U}_{n}(\widetilde{M}^n) $
with $ \widetilde{M}^{n} := \sum_{m=n}^{N-1} \tilde{\xi}_m = \tilde{\xi}_n + \widetilde{M}^{n+1} $. We then consider the following backward minimization problem.

\begin{problem}[Pure Dual Backward minimization problem]\label{pb:backward_minimization}
    Given $n \in \N^{-1}$ and optimal martingale differences $\tilde{\xi}_{n+1},\ldots,\tilde{\xi}_{N-1}$, we solve
    \begin{equation}\label{eq_optim}
        \tilde{\xi}_n \in \arginf_{\xi_n \in \mathcal{P}_n}\mathbb{E} \big[\widetilde{U}_{n}({\xi}_n ; \widetilde{M}^{n+1}) \big] .
    \end{equation}
If $ Y^{*}_n \ge 0 ,\; n=0,\ldots,N $, $ \mathbb{P} $-a.s., we instead solve
  \begin{equation}\label{eq_optim_second_moment}
        \tilde{\xi}_n \in \arginf_{\xi_n \in \mathcal{P}_n}\mathbb{E}\big|\widetilde{U}_{n}({\xi}_n ; \widetilde{M}^{n+1}) \big|^2 .
    \end{equation}
\end{problem}

\begin{remark}
    Clearly, Doob martingale $ M^{*} $ solve the above two minimization problems simultaneously.
\end{remark}

\subsection{Expressivity}\label{subsec:numerical_express}
We now investigate the expressivity of the numerical integration scheme, with a particular focus on the choice of $K$. This requires additional structural conditions on the It\^o diffusion so that expression rates can be derived. To make the dependence on the dimension transparent, we attach subscripts/superscripts $d$ where needed.

For any $x$-variate function $ f $, denote
$\Lip f:= \sup_{x\neq y}\frac{\| f(x) - f(y) \|}{\|x-y\|} $
as the minimal Lipschitz constant of $f$ (when $f$ takes values in a matrix (tensor) space, $\|f(x)\|$ is $\| f(x) \|_{\H}$).
If $f$ is $\mathbb{R}^{d\times d}$-valued, we further define
$ \Lip^i f:= \sup_{x\neq y} \frac{\| f^{i}(x) - f^{i}(y) \|  }{\|x-y\| } $
(component-wise Lipschitz, where $ f^i $ is the $ i $-th column of $ f $), and
$ \LipH f := ( \sum_{i=1}^d |\Lip^i f |^2 )^{\frac{1}{2}} $
(the stronger Lipschitz norm). Similarly, for time-variate functions, we denote by $\Holder f$ the minimal $\frac{1}{2}$-H\"{o}lder constant of $f$ in $t$.

The numerical integration of \eqref{eq:Martingale_represent} yields
\begin{equation*}
    M^{*;d}_{t_{n+1}} - M^{*;d}_{t_n} = Y^{*;d}_{n+1} - \mathbb{E}_n\bigl[Y^{*;d}_{n+1}  \bigr] = \int_{t_n}^{t_{n+1}} Z^{*;d}_s \,dW^d_s   .
\end{equation*}
Following \citep{belome09}, we consider the associated non-driver decoupled FBSDE: for $ n\in \N^{-1} $ and $x \in \mathbb{R}^d $,
\begin{equation}\label{eq:problem_related_FBSDE}
    \begin{aligned}
        X^d_t &= X^d_{t_n} + \int_{t_n}^t a^d(u,X^d_u) \,du + \int_{t_n}^t b^d(u,X^d_u) \,dW^d_u ,\\
        Y^d_t &= V^d_{n+1}(X^d_{t_{n+1}}) - \int_{t}^{t_{n+1}} Z^{*;d}_u \,dW^d_u, \quad t_n \le t \le t_{n+1}.
    \end{aligned}
\end{equation}
In order to control \eqref{eq:problem_related_FBSDE}, we require expressivity properties for $V^d$. These are obtained by backward recursion starting from the expressivity of the payoff $g^d$; details are provided in the Appendix.

For generality, we also consider the following non-driver decoupled FBSDE with terminal function $\bar{g}^d$: for $x \in \mathbb{R}^d$,
\begin{equation}\label{eq:general_FBSDE}
    \begin{aligned}
        X^d_t &= x + \int_{0}^t a^d(s,X^d_s) \,ds + \int_{0}^t b^d(s,X^d_s) \,dW^d_s,  \\
        Y^d_t &= \bar{g}^d(X_T) - \int_{t}^T Z^d_s \,dW^d_s, \quad 0 \le t \le T.
    \end{aligned}
\end{equation}
As discussed in \citep{zhangjianfeng17}, solvability of \eqref{eq:general_FBSDE} requires suitable conditions on $(a^d,b^d)$ and $\bar{g}^d$. The same type of conditions is also essential for deriving dimension-explicit approximation rates. We therefore strengthen the assumptions in \cite{zhangjianfeng17} and adopt the following \emph{structural condition}.

\begin{assumption}\label{ass:N_0_structural_ass_model_deter}
    The functions $a^d(t,x)$, $b^d(t,x) $ in \eqref{eq:SDE} satisfy the condition that, for any $t\in [0,T]$, $x\in \mathbb{R}^d $,   
    \[
         \Lip a^d(t,\cdot) \le c (\log d)^{\frac{1}{2}} , \; \; \LipH  b^d(t,\cdot) \le C (\log d)^{\frac{1}{4}}
    \]
    \begin{equation*}
         \|a^d(t,0)\|,\|b^d(t,0)\|_{\H} , \Holder a^d(\cdot,x), \Holder b^d(\cdot,x)  \le C d^Q
    \end{equation*}
    for some positive constants $C,Q$, independent of $d$.
\end{assumption}

\begin{remark}\label{rem:stronger_Lip_remark}
    The volatility term requires the stronger Lipschitz condition $ \LipH  b^d(t,\cdot) \le C (\log d)^{\frac{1}{4}} $. Indeed, for a matrix-valued function $ f $,
    \[
        (\Lip  f )^2 = \sup_{x\neq y} \frac{ \sum_{i=1}^d \| f^i(x) - f^i(y) \|^2 }{\|x-y\| ^2 } \le \sum_{i=1}^d \sup_{x\neq y} \frac{\| f^i(x) - f^i(y) \|^2}{\|x-y\|^2} \le \sum_{i=1}^d (\mathrm{Lip}^i f)^2 = (\LipH f)^2 .
    \]
    Hence,
    \( \Lip  b^d  \le \LipH  b^d \le C (\log d)^{\frac{1}{4}} \le C (\log d)^{\frac{1}{2}} \).
\end{remark}

Assumption~\ref{ass:N_0_structural_ass_model_deter} imposes a sparsity structure on high-dimensional dynamics. Nevertheless, it accommodates a broad class of commonly used models. In particular, both geometric Brownian motion and constant volatility models satisfy this assumption. A more detailed discussion is deferred to Remark~\ref{rem:discussion_of_expression_rate_model}.

\begin{assumption}\label{ass:N_0_structural_ass_g_deter}
     The function $\bar{g}^d$ in \eqref{eq:general_FBSDE} satisfies \( \Lip \bar{g}^d , \;
        |\bar{g}^d(0)|  \le C d^Q \)
for some positive constants $C,Q$, independent of $d$.
\end{assumption}

\begin{remark}
    The constants $C$ and $Q$ in Assumptions \ref{ass:N_0_structural_ass_model_deter} and \ref{ass:N_0_structural_ass_g_deter} can be chosen to coincide by taking maxima if necessary.
\end{remark}

Let $[0,T]$ be uniformly partitioned with spacing $h:=T/N_h\le 1$, $t_i:=ih$, $i \in \N_h$, and define $ \overline{Z}^{h;d}_t := \sum_{i=0}^{N_h-1} {Z}^{h;d}_{t_i} 1_{[t_i , t_{i+1} )}(t) $. The following theorem shows that the discretization level is expressive, providing a dimension-explicit counterpart of classical SDE/BSDE estimates; proofs are given in the Appendix.

\begin{theorem}[Expressivity of Numerical Integration Approximation] \label{thm:Numerical_integration_est}
    Under Assumption~\ref{ass:N_0_structural_ass_model_deter} and Assumption \ref{ass:N_0_structural_ass_g_deter}, there exist positive constants $\tilde{B},\tilde{Q} $ such that for \eqref{eq:general_FBSDE}, we have
    \[
        \big\| Z^d - \overline{Z}^{h;d} \big\|^2_{\L^2} \le \tilde{B}d^{\tilde{Q}} (1+\|x\|^2) h  , \; \; \forall \; h \in [0,1] .
    \]
\end{theorem}

We can now translate this estimate into an explicit expressivity bound on the number of integration points $K$.

\begin{theorem}[Expressivity of integration points] \label{theorem:express_N_0_new}
    Under Assumption~\ref{ass:N_0_structural_ass_model_deter} on every $[t_n,t_{n+1}]$, $n\in \N^{-1} $,  Assumption~\ref{ass:N_0_structural_ass_g_deter} for the payoff function on $ \T^N $, i.e., $\{g^d(t_n,\cdot) \}_{n=0}^{N}$, and $ \| x_0 \| \le C d^Q $ (the same constants as in above assumptions), there exist positive constants ${B}^{*}, {Q}^{*} $ independent of $d$, such that for any $\varepsilon>0 $, there exists an $K_{d,\varepsilon} \in \mathbb{N}$ satisfying \(K_{d,\varepsilon} \le {B}^{*} d^{{Q}^{*}} \varepsilon^{-1}  \),
    so that, for all $n \in \N^{-1} $,
    \begin{equation*}
        \big\| \widehat{Z}^{*,K_{d,\varepsilon},n ; d} - Z^{*;d} \big\|^2_{[t_n,t_{n+1}],2} \le \varepsilon .
    \end{equation*}
\end{theorem}

\begin{remark}
    The expressivity result above for the numerical approximation can be extended to European option hedging problems with only minor modifications; we leave a full treatment for future work. Moreover, it suggests that the rebalancing frequency in the hedging strategy can be chosen to grow at most polynomially as the dimension $d$ increases.
\end{remark}

\section{DeepMartingale}\label{sec:deep_mtg}
This section presents our DNN architecture and the corresponding approximation of \eqref{eq_Z}. By establishing a universal approximation theorem (UAT), we obtain a tight dual upper bound with a theoretical convergence guarantee. Since our construction relies on DNN parameterizations, we refer to it as the \textit{DeepMartingale} approach. We further demonstrate its expressivity. Although the expressivity analysis is formulated via the value function, the method itself does not require any primal information; hence, it is a \textit{pure dual} approach in the spirit of \cite{roger10,puredual-mf}, and it provides an efficient and accurate high-dimensional hedging strategy.

Motivated by \eqref{eq_Z} and \eqref{eq_Z_converge}, our goal is to construct NNs that approximate $\widehat{Z}^{*,K,n}$ for $n\in \N^{-1}$. By the Markov property of $X$, the following lemma justifies a functional representation of $\widehat{Z}^{*,K,n}$ in terms of $X$. With a slight abuse of notation, we use $\widehat{Z}^{*,K,n}$ to denote both the process and its function representation.

\begin{lemma} \label{lemma:represent_z_x}
    Given any $n \in \N^{-1} $, there exists a Borel measurable $\widehat{Z}^{*,K,n} : \mathbb{R}^{1+d} \rightarrow \mathbb{R}^{d} $, such that for all $ k\in \K $, \( \widehat{Z}^{*,K}_{t^n_k} = \widehat{Z}^{*,K,n}(t^n_k, X_{t^n_k} )     \), and specifically,
    \begin{equation}\label{eq:Z_def}
            \widehat{Z}^{*,K,n}(t^n_k, X_{t^n_k} ) = \frac{1}{\Delta t}\mathbb{E}^{X_{t^n_k}}\bigl[ V_{n+1}  (X_{t_{n+1}} ) \Delta W_{t^n_k} \bigr] .
    \end{equation}
\end{lemma}

\begin{remark}\label{rem:continuous_representation_sde}
    Lemma~\ref{lemma:represent_z_x} can be proved by using the representation \( X_t(\omega) = F(X_r(\omega), r , t ,\omega ),\; t\ge r \) in \cite[Proof of Theorem~7.2]{SDE03}, and by observing that $\omega \mapsto (V_{n+1} \circ F ) (x, t^n_k, t_{n+1}, \omega ) \Delta W_{t^n_k} (\omega)$ satisfies Lemma~\ref{lemma:independent_measurable_relation} provided in the Appendix. We omit the detailed proof.
\end{remark}

\begin{remark}
    \label{remark:priori_approx}
    Our expressivity analysis for \textit{DeepMartingale} is guided by \eqref{eq:Z_def}. In particular, if an NN $\widehat{V}$ approximates $V_{n+1}$ and a random NN $\widehat{f}(\cdot,\omega)$ approximates the It\^o process $X$, then $Z^{*}$ can be approximated via expectation-expressivity arguments similar to those in \citep{Jentzen23} and \citep{gonon23}. A detailed discussion is provided in Subsection~\ref{subsec:expressivity}.
\end{remark}

\subsection{Neural network architecture}
Let \(\Theta:=\bigcup_{m\ge1}\mathbb{R}^m\). For each \(n \in \N^{-1} \), let \(I\ge1\) and \(q_1,\ldots,q_I\in\mathbb{N}\), $ q_0 := 1+d $, $ q_{I+1}:= d $. Define the FNN \(z_n^{\theta_n}:\mathbb{R}^{1+d}\to\mathbb{R}^d\), \(\theta_n\in\mathbb{R}^{Q_n}\), by
\begin{equation}
z^{\theta_n}_n(t,x)=a^{\theta_n}_{I+1}\circ \varphi_{q_I}\circ a^{\theta_n}_{I}\circ\cdots\circ \varphi_{q_1}\circ a^{\theta_n}_1(t,x),
\label{neural_net}
\end{equation}
where \(\varphi_j:\mathbb{R}^j\to\mathbb{R}^j\) is a component-wise non-constant activation and \(a_i^{\theta_n}(u)=A_i u+b_i\) are affine maps with \(A_i\in\mathbb{R}^{q_i\times q_{i-1}}\) for \(i=1,\ldots,I+1\),
and \(b_i\in\mathbb{R}^{q_i}\) for \(i=1,\ldots,I+1\). The parameter dimension is \( Q_n=\sum_{i=1}^{I+1}\big( q_i q_{i-1} + q_i  \big) \).

Motivated by \eqref{eq_Z}, we construct \textit{DeepMartingale} as follows. Let $ \theta:=(\theta_{n} )_{n=0}^{N-1} \in \Theta^N $,
\begin{equation}
    \xi^{\theta_n,K}_n := \sum_{k=0}^{K-1} z^{\theta_n}_n(t^n_k,X_{t^n_k}) \cdot \Delta W_{t^n_k} , \; n \in \N^{-1} ,\quad M^{\theta,K}_n := \sum_{m=0}^{n-1} \xi^{\theta_m}_m ,\; n \in \N ,
      \label{deep_mart}
\end{equation}
Clearly, $\xi^{\theta_n,K}_n \in \mathcal{P}_n$, and hence $M^{\theta,K}\in \mathcal{M}^N$.

\begin{remark}[Delta hedge by DeepMartingale]
\label{rem:delta_hedge_deepMart}
Comparing \eqref{eq:Martingale_represent} with Remark~\ref{rem:connection_delta_hedging}, the network $z^{\theta_n}_n$ also induces a ``deep'' delta hedge for Bermudan options. Given $ n\in \N^{-1} $,
\[
     \vartheta^{\theta_n}_n(t,x) := b^{-1}(t,x) z^{\theta_n}_n(t,x) ,\; t\in [t_n,t_{n+1}] .
\]
Moreover, the same construction applies to classical European option hedging. From a practitioner perspective, $z^{\theta_n}_n$ corresponds to a continuous-time hedging rule, but it can only be trained from ``less'' frequent observations indexed by $K$. The actual rebalancing frequency may be chosen independently, taking potential market frictions into account. In the literature, e.g.\ \cite{puredual-mf,guo2024simultaneousupperlowerbounds}, it is common to align the rebalancing frequency with the observation frequency; we adopt the same convention in our DeepMartingale implementation.
\end{remark}

\subsection{Convergence under bounded activation}
We now outline the logic used to prove convergence. We first introduce an auxiliary metric for measuring the approximation error of the dual upper bound. Since the induced finite Borel measure does not necessarily have compact support, we restrict attention to bounded activation functions $\varphi$, which enables us to apply the UAT in \citep{hornik91}. We then establish a UAT for the integrand under this measure, which yields the convergence result for \textit{DeepMartingale}.

\subsubsection{Auxiliary metric for applying UAT}
By the error propagation Lemma~\ref{lem-error-propagate}, it is necessary to work with a metric that is compatible with the required $L^2$ approximation. For each $n \in \N^{-1}$, we define the following finite Borel measures. Let \( \mu^n_k \) be the distribution of $ X_{t^n_k} $, and define
\(
    \lambda^n_k(A) :=  1_A(t^n_k) \Delta t^n_k ,\; A\in \mathcal{B}(\mathbb{R})
\),
\begin{equation}
    \mu^{K}_n(A) := \mathbb{E}\Bigl[\sum_{k=0}^{K-1} \mathbf{1}_{A}(t^n_k, X_{t^n_k}) \Delta t \Bigr] = \sum\limits_{k=0}^{K-1} (\lambda^n_k \otimes \mu^n_k)(A)  ,\quad A\in \mathcal{B}(\mathbb{R}^{1+d})  . \label{eq:measure_mu_def}
\end{equation}

Using standard arguments from measure theory, we obtain the following lemma characterizing integration with respect to $\mu^K_n$.

\begin{lemma} \label{lemma: inte_relation_mu}
    For each $n \in \N^{-1} $ and any Borel measurable function $f: \mathbb{R}^{1+d} \rightarrow \mathbb{R}^m $, if $\mathbb{E}\|f(t^n_k,X_{t^n_k})\|   < \infty$, $k \in \K $, then $f \in L^1_{1+d , m}(\mu^K_n) $. In addition, for $m=1$,
    \begin{equation}
        \mathbb{E}\Bigl[\sum_{k=0}^{K-1} f(t^n_k,X_{t^n_k}) \Delta t \Bigr] = \int_{\mathbb{R}^{1+d}} f(t,x) \, \mu^{K}_n (dt dx) .
        \label{eq_E_int_relation}
    \end{equation}
\end{lemma}

By Lemma~\ref{lemma: inte_relation_mu}, we can verify that $\widehat{Z}^{*,K,n} \in L^2_{1+d,d}(\mu^{K}_n)$. We next construct NN approximators for $\widehat{Z}^{*,K,n}$ in $L^2_{1+d,d}(\mu^{K}_n)$.

\subsubsection{UAT \& Convergence}
The following UAT establishes convergence of the NN approximation to the integrand $\widehat{Z}^{*,K,n}$ in $L^2_{1+d,d}(\mu^{K}_n)$, as well as convergence of the induced martingale increments to $\xi^{*,K}_n:=\xi_n(\widehat{M}^{*,K})$ in $L^2(\mathcal{F}_{t_{n+1}},\R )$.

\begin{theorem}[Universal Approximation Theorem] \label{theorem:inte_L2_approx}
    For any $\varepsilon >0 , K\in \bbN_{+} $, there exist parameters $ \theta^K_{\varepsilon} := (\theta^K_{n,\varepsilon})_{n=0}^{N-1} \in \Theta^N  $, such that for any $n \in \N^{-1} $,
    \begin{equation}
        \big\|  z^{\theta^K_{n,\varepsilon}}_n - \widehat{Z}^{*,K,n}   \big\|_{2;\mu^K_n} \; ,\; \big\|  \xi^{\theta^K_{n,\varepsilon},K}_n - \xi^{*,K}_n \big\|_{L^2}  < \varepsilon   .
        \label{eq_integrand_approx}
    \end{equation}
\end{theorem}

Combining \eqref{eq_E_int_relation} and \eqref{eq_error1} with \eqref{deep_mart}, we obtain the following approximation result for the deep dual upper bound.

\begin{theorem}[Approximation Capability of \textit{DeepMartingale}]\label{theorem:theta_approx}
    For any $\varepsilon>0 , K \in \mathbb{N}_{+} $, there exists a DeepMartingale $ M^{\theta^K_{\varepsilon} ,K} $ such that for each $n \in \N $,
    \begin{equation*}
        \big\| \widetilde{U}_n(M^{\theta^K_{\varepsilon} ,K}) - \widetilde{U}_n(\widehat{M}^{*,K}) \big\|_{L^2} \le (N-n)\varepsilon .
    \end{equation*}
\end{theorem}

Since $\mathbb{E} [\widehat{Y}^{*,K}_n] = \mathbb{E}[\widetilde{U}_n(\widehat{M}^{*,K})] $, the tightness of the \textit{DeepMartingale} upper bound for $Y^{*}$ follows immediately.

\begin{corollary}[DeepMartingale Solvability of Pure Dual Problem~\eqref{pb:backward_minimization}] \label{coro:tight_upper}
    For all $ n \in \N^{-1} $,
    \begin{itemize}
        \item[(i)] \( \mathbb{E} [Y^{*}_n] = \lim_{K \to \infty} \inf_{\theta \in \Theta^N } \mathbb{E}\bigl[ \widetilde{U}_n(M^{\theta,K})\bigr] \).
        \item[(ii)] If $ Y^{*}_n \ge 0 ,\; n=0,\ldots,N $, $ \mathbb{P} $-a.s., then \(  \mathbb{E} |Y^{*}_n|^2 =  \lim_{K \to \infty} \inf_{\theta \in \Theta^N }\mathbb{E}\big| \widetilde{U}_n(M^{\theta,K}) \big|^2 \).
    \end{itemize}
\end{corollary}

For any function \(F:\Theta^N \to\mathbb{R}\), define the $\varepsilon$-optimal set
\[
    \mathrm{argmin}_{\varepsilon} F:=\{\theta \in \Theta^N :\ F(\theta) \le \inf_{\vartheta \in \Theta^N }F(\vartheta)+\varepsilon\}.
\]
We next state a proposition that justifies the use of the second-moment loss in the pure dual problem~\eqref{pb:backward_minimization}.

\begin{proposition}\label{pro:L2-to-L1-loss}
     Let $ \varepsilon_{K} \downarrow 0 $ as $ K \uparrow \infty $. If $ Y^{*}_0 \ge 0 ,\; n\in \N^{-1} $, $\mathbb{P} $-a.s., then given any $ K \ge 1 $,
        let
        $ \theta^{K} \in \argmin_{\varepsilon_{K}}
        \mathbb{E}\big| \widetilde{U}_n(M^{\theta,K}) \big|^2  $, then
        \begin{itemize}
            \item[(i)] \( \lim_{K \to \infty} \mathbb{E}\big[\mathrm{Var}_n(\widetilde{U}_n\big(M^{\theta^{K},K}) \big) \big] = 0 \).
            \item[(ii)] If further $ (M^{\theta^{K},K})_{K\ge1} $ is uniformly integrable,
            then \( \mathbb{E}[Y^{*}_n] = \lim_{K \to \infty}\mathbb{E}\big[ \widetilde{U}_n(M^{\theta^{K},K}) ] \).
        \end{itemize}
\end{proposition}

\begin{remark}
    Proposition~\ref{pro:L2-to-L1-loss} shows that, under suitable conditions, minimizing the second-moment objective also drives the minimization of the (first-moment) dual upper bound. The uniform integrability condition have been discussed in \cite{schoen13}. This convergence behavior is consistent with our numerical experiments; see Section~\ref{sec:Implement}. We emphasize that Proposition~\ref{pro:L2-to-L1-loss} does not depend on neural networks or on the choice of activation function; it is a general convergence statement.
\end{remark}

In summary, Problem~\ref{pb:backward_minimization} can be solved by \textit{DeepMartingale} with a convergence guarantee once the activation is bounded.

\subsection{Expressivity \& Convergence under unbounded ReLU activation}  \label{subsec:expressivity}
Here, we establish the expressivity result of \textit{DeepMartingale}, namely that it yields a tight dual upper bound while the NN size grows at most polynomially in $d$ and $\varepsilon$. This provides a theoretical guarantee that \textit{DeepMartingale} can overcome the curse of dimensionality.
Therefore, expressivity in our setting directly translates into scalability. This observation also motivates our numerical estimation of a dimension scaling law, which we then use to scale the network architecture and the training setup in the deep learning-based numerical computations (see Section~\ref{sec:Implement}).
To develop the theory, we first introduce a random neural network (RanNN) framework in an infinite-width setting with an RKBS treatment, which ensures both generality and a mathematically well-posed parameter space. This framework can be viewed as a multilayer, infinite-width extension of the random NN architecture proposed in \citep{GononRandomNN23}. We then prove an expressivity result for value function approximation under structural conditions using RanNN. Based on such value-function approximation, we construct a ``deep integrand'', equivalently a ``deep delta hedge''. Under strengthened structural conditions, we further derive the expressivity of \textit{DeepMartingale}. These strengthened conditions are not restrictive, as they are satisfied by many practical models; see Subsection~\ref{subsubsec:affine_ito} and Remark~\ref{rem:discussion_of_expression_rate_model} for discussions.

\subsubsection{Infinite-width neural network with RKBS treatment and random neural network}
To rigorously formulate our framework, we define RanNNs as neural networks whose parameters are random variables. Such networks have also been investigated from a computational perspective in \citep{DeepPrimalRandom}. Due to the randomness in width, we adopt an infinite-dimensional RKBS approach \citep{BARTOLUCCI2023194,RKBS2024}, under which the metric structure and measurability of the parameter space arise naturally. All proofs in this subsection are deferred to the Appendix.

Let $\ell^2(\mathbb{N}_{+}) $ be the space of square-summable sequences. Set $\Theta_0=\{0,\ldots,d_1\}$, $\Theta_{I+1}=\{1,\ldots,d_3 \}$, $\Theta_i=\mathbb{N}_{+}$ for $1\le i\le I $, $\mathcal{X}_0 =\mathbb{R}^{d_1}$, $\mathcal{X}_{I+1}=\mathbb{R}^{d_3}$, and $\mathcal{X}_i=\ell^2(\mathbb{N})$ for $1\le i\le I $. Following \citep{RKBS2024}, we denote by $\mathcal{M}(\Theta_i, \mathcal{X}_{i+1}) $ the Banach space of $\mathcal{X}_{i+1}$-valued vector measures on $\Theta_i$, $0 \le i \le I $, endowed with the total variation norm $\|\mu\|_{\text{TV}} = |\mu|(\Theta_i)$ for $\mu\in \mathcal{M}(\Theta_i, \mathcal{X}_{i+1})$, where $|\mu|$ is a bounded positive measure on $\Theta_i$. We use $ \| \cdot \|_{\mathcal{X}_i} $ to denote the norm on $ \mathcal{X}_i $.

For a given depth $I\ge 1 $, we generalize the usual composition of finite-dimensional nonlinear vector functions to the infinite-dimensional setting, as illustrated by the following graph.

\begin{center}
\begin{tikzpicture}[
    arrow/.style={->,>=latex},
    node distance=1.2cm,
    scale=0.7,
    transform shape
]
    \node (R1) at (0,0) {$\mathbb{R}^{d_1}$};
    \node (P1) at (2,0) {$\ell^2(\mathbb{N})$};
    \node (dots) at (4,0) {$\cdots$};
    \node (P2) at (6,0) {$\ell^2(\mathbb{N})$};
    \node (R2) at (8,0) {$\mathbb{R}^{d_3}$};
    
    \draw[arrow] (R1) -- node[above] {$f_1$} (P1);
    \draw[arrow] (P1) -- (dots);
    \draw[arrow] (dots) -- (P2);
    \draw[arrow] (P2) -- node[above] {$f_{I+1}$} (R2);
    
    \draw[arrow] (R1) to[bend right=30] node[below] {$f^{\text{deep}}$} (R2);
\end{tikzpicture}
\end{center}

Specifically, for $1\le i\le I$, define
\[\rho_0(x,n):=\begin{cases}
    1 & n=0 \\
    x_n & n=1,\ldots,d_1
\end{cases},\; x \in \mathbb{R}^{d_1}, \quad  \rho_i(x,n):=\begin{cases}
    1 & n=0 \\
    \sigma(x_{n-1}) & n\ge 1
\end{cases} ,\; x\in \ell^2(\mathbb{N})\]
For $\mu_i\in \mathcal{M}(\Theta_i, \mathcal{X}_{i+1})$, $0 \le i\le I $, it is clear that \( \mu_i  =\sum_{m=0}^{M} w^{i+1}_m \delta_m \),
where $M=d_1 $ if $i=0 $, and $M=\infty $ otherwise; $w^{i+1}_m \in \mathcal{X}_{i+1}$, $0\le i\le I$, $m\ge 0$ and $\delta_m $ is the Dirac delta function.

\begin{definition}[Infinite-width Neural Network]\label{def:infinite_width}
    We call $f=f_{I+1}\circ f_{I} \circ \cdots \circ f_1 :\mathbb{R}^{d_1}\rightarrow \mathbb{R}^{d_3} $ the infinite-width neural network with depth $I$
    if for $1\le i\le I-1 $ and $y\in \mathcal{X}_{i-1} $,
    \begin{equation*}
        f_i (y) = \int_{\Theta_{i-1}} \rho_{i-1}(y,\theta_{i-1})\,d\mu_{i-1}(\theta_{i-1}).
    \end{equation*}
\end{definition}

As shown in \citep{RKBS2024}, Definition~\ref{def:infinite_width} is equivalent to the following familiar operator form. Let $W^{i+1}:\mathcal{X}_i \rightarrow \mathcal{X}_{i+1}$, $0\le i\le I $, be bounded linear operators such that \(  W^{i+1}x = \sum_{m=1}^{M} w^{i+1}_m x_{m-1} \),
where $M=d_1 $ if $i=0 $, and $M=\infty $ otherwise. Setting $b^{i+1} = w^{i+1}_0$, $0\le i\le I$, we have $f_1(x)=\sum_{m=0}^{d_1} w^1_m \rho_0(x,m) = W^1 x + b^1$ and, for $1\le i\le I-1 $,
\[
    f_{i+1}(y) = \sum_{m=0}^{\infty} w^{i+1}_m \rho_{i}(y,m) = W^{i+1}(\sigma(y)) + b^{i+1}.
\]
This form coincides with the standard feed-forward neural network (FNN) structure if we truncate $\ell^2(\mathbb{N})$ to a Euclidean subspace.

Under this formulation, we parametrize DNNs via a Banach space of vector-valued measures with finite total variation, which equips the parameter space with a natural metric and ensures measurability. To properly define RanNNs, we introduce random parameters in the network. Let $\mathcal{U}:=\prod_{i=0}^{I} \mathcal{M}(\Theta_i,X_{i+1}) $ be the product parameter space, equipped with the product metric induced by the total variation norms on $\mathcal{M}(\Theta_i,X_{i+1})$. We view the infinite-width NN as a function of both inputs and parameters, $f:\mathbb{R}^{d_1}\times \mathcal{U} \rightarrow \mathbb{R}^{d_3} $. Since bounded linear operators in finite-width NNs correspond to finite-dimensional matrices, the total variation norm is consistent with the Hilbert--Schmidt norm in that setting. By Borel measurability under the product metric of $\mathbb{R}^{d_1}\times \mathcal{U}$, continuity of $f$ with respect to $x\in \mathbb{R}^{d_1}$ and $\mu\in \mathcal{U}$ follows.

\begin{proposition}[Continuity]\label{prop:contin_RKBS}
    The infinite-width NN is continuous w.r.t.\ $x\in \mathbb{R}^{d_1}$ and $\mu\in \mathcal{U}$.
\end{proposition}

We now define the linear growth functional (minimal rate of linear growth): for any function $ f: \mathbb{R}^{d_1} \to \mathbb{R}^{d_2} $, define \( \Growth(f) := \sup_{x\in \mathbb{R}^{d_1}} \frac{\|\tilde{f}(x)\|}{1+\|x\|} \). The size functional $ \size(\cdot) $ of an infinite-width network $ f^{\mathrm{deep}} $ with depth $I$ (including standard ReLU networks) is defined as the total number of non-zero entries across layers: \(   \size\big(f_i(*,\mu^{'}_{i-1})\big) = \sum_{m=0}^{M}\sum_{n=1}^{M^{'}} \mathbf{1}_{(w^{i+1}_{mn} \neq 0)} \), where $ M^{'} = d_3 $ if $ i=I+1 $ and $ \infty $ otherwise.

We can now formally define RanNNs and the associated random functionals.

\begin{definition}[Random Feed-Forward Neural Network]\label{def:random_NN}
    In the probability space $(\Omega,\mathcal{F},\mathbb{P}) $, let $\mu(\cdot) := (\mu_0(\cdot),\ldots,\mu_I(\cdot) ): \Omega \rightarrow \mathcal{U} $ be a $\mathcal{F}/\mathcal{B}(\mathcal{U}) $-random variable. For an infinite-width NN $f $ with depth $I$ (Definition~\ref{def:infinite_width}), $\tilde{f}: \mathbb{R}^{d_1}\times \Omega \rightarrow \mathbb{R}^{d_3} $ is called a random feed-forward neural network (RanNN) of depth $I$ w.r.t.\ $\mu$ if $\tilde{f}(x,\omega)=f(x,\mu(\omega)) $. Here and below, we do not distinguish between $\tilde{f}(x,\mu) $ and $\tilde{f}(x,\omega)$.
    The size, growth rate, and Lipschitz random variable of $\tilde{f} $ are defined as follows:
    \[
        \size(\tilde{f}):\omega \mapsto \size(\tilde{f}(\cdot,\omega)) ,\; \Growth(\tilde{f}): \omega \mapsto \Growth(\tilde{f}(\cdot,\omega)) ,\;  \Lip(\tilde{f}): \omega \mapsto \Lip(\tilde{f}(\cdot,\omega)) .
    \]
\end{definition}

The measurability of $\size$, $ \Growth $ and $ \Lip $ follows directly.

\begin{proposition}[Measurability of $\size$, $ \Growth $, and $ \Lip $]\label{prop:size_measurable}
    For any RanNN (Definition~\ref{def:random_NN}) $ \tilde{f} $,
    $\size(\tilde{f})$, $\Growth(\tilde{f})$, and $\Lip(\tilde{f}) $ are random variables.
\end{proposition}

\subsubsection{Structural Framework}
We now collect the expressivity assumptions imposed on the dynamics and the payoff structure. This includes discrete-time and continuous-time (pathwise) NN representations for the dynamics, as well as NN approximation assumptions for the payoff function. These assumptions form the basis for the subsequent expressivity analysis of \textit{DeepMartingale}.

To obtain an extended expressivity result for the value function approximation, we introduce a non-specific dynamic assumption to maintain theoretical generality. From now on, we attach superscripts/subscripts $d$ whenever needed to emphasize dimension dependence. Fix $ p>0 $.

\begin{assumption}[Dynamic Assumption with Order $p$]\label{ass:dynamic_ass}
We make the following assumption.
\begin{enumerate}
    \item[(1)]
    (\textit{Independent innovation})
    The process \(X^d\) is \(\mathbb{F}^N\)-adapted (not restricted to diffusion \eqref{eq:SDE}). For each \(n \in \N^{-1} \), there exists a \(\mathcal{B}(\mathbb{R}^d)\otimes\mathcal{F}_{t_{n+1}}\)-measurable map \(f_n^d:\mathbb{R}^d\times\Omega\to\mathbb{R}^d\) such that (a).
    \(
X^d_{t_{n+1}}=f_n^d(X^d_{t_n},\cdot)\quad \mathbb{P}\text{-a.s.},
    \)
    and (b). for every \(x\in\mathbb{R}^d\), \(f_n^d(x,\cdot)\) is independent of \(\mathcal{F}_{t_n}\) and \(f_n^d(x,\cdot)\in L^1(\mathcal{F}_{t_{n+1}},\R^d) \).

    \item[(2)] (\textit{Representation}) There exist constants $c,q > 0$ independent of $d$, such that for $ n\in \N^{-1} $,
    \begin{enumerate}

    \item[(a)] \( \| \Growth(f^d_n) \|_{L^p} \leq c d^{q} \).

    \item[(b)] there exists a RanNN (Definition~\ref{def:random_NN}) $\widehat{f}^d_n:\mathbb{R}^d \times \Omega \rightarrow \mathbb{R}^d$ with depth $I^d_n \leq cd^q $ such that $f^d_n$ can be represented by $\widehat{f}^d_n $, i.e., $f^d_n(x,\cdot) = \widehat{f}^d_n(x,\cdot)$, $\forall x\in \mathbb{R}^d$, $\mathbb{P}$-a.s.

    \item[(c)] RanNN approximator $\widehat{f}^d_n$ in (b) satisfies \( \mathbb{E}[\size(\widehat{f}^d_n)] \leq c d^q \).
    \end{enumerate}

\end{enumerate}
\end{assumption}

\begin{remark}
    Clearly, the independent innovation assumption (1) in Assumption~\ref{ass:dynamic_ass} implies the $\F^N$-Markov property of $X^d$. In particular,
    $\mathbb{E}_n[ X^d_{t_{n+1}}  ] = \mathbb{E}^{X^d_{t_n}} [ X^d_{t_{n+1}}  ] = \mathbb{E}[ f^d_n(x,\cdot ) ]|_{x =X^d_{t_n} }$.
\end{remark}

\begin{remark}
    The above assumption is stronger than that in \citep{gonon23}. However, in our setting no expressivity condition on a pathwise Lipschitz constant is required, which allows certain continuous-time processes (in contrast to the discrete-time models considered in \citep{gonon23}), such as affine It\^o processes, to fall within our framework for deriving value function approximation results. In \citep{gonon23}, the assumption is formulated in terms of a pathwise Lipschitz expressivity property (i.e., an approximation rate holding for $\mathbb{P}$-a.s.\ $\omega$), which cannot be directly applied in a continuous-time setting.
\end{remark}

We also formulate a continuous-time (pathwise) dynamic assumption for our integrand approximation, which accommodates more frequent observations between monitoring points, including continuous-time observations. Fix $\tilde{p} > 0 $.

\begin{assumption}[Pathwise Dynamic Assumption with Order $\tilde{p}$]\label{ass:stronger_ass}
We make the following assumption.
    \begin{enumerate}
        \item[(1)] (\textit{It\^o Diffusion}) $ X^d $ is the unique strong solution of \eqref{eq:SDE}. In particular, $ X^d $ admits a representation $ X^{d}_t = f^d(X^d_s,s,t,\cdot)$ for all $ 0\le s\leq t\leq T $ as argued in Remark~\ref{rem:continuous_representation_sde}, where $ (x,\omega) \mapsto f^d(x,s,t,\omega) =: f_s^{t;d}(x,\omega) $ is $ \calB(\R^d) \otimes \mathcal{F}_t  $-measurable, and $  f^d(x,s,t,\cdot)$ is independent of $ \mathcal{F}_{s} $ for all $x\in \R^d $.

        \item[(2)] (\textit{Representation})  There exist $\bar{c},\bar{q} > 0 $ independent of $d$, such that for any $n \in \N^{-1} $, \( t_n\leq s< t\leq t_{n+1} \),
        \begin{enumerate}
            \item[(a)] \( \| \Growth(f_s^{t;d}) \|_{L^{\tilde{p}}} \le \bar{c}d^{\bar{q}} \).

            \item[(b)] There exists a RanNN (Definition~\ref{def:random_NN}) $\widehat{f}_s^{t;d}: \mathbb{R}^d \times \Omega \rightarrow \mathbb{R}^d$ with depth $I_s^{t;d} \leq \bar{c}d^{\bar{q}} $ such that $f_s^{t;d} $ can be represented by $\widehat{f}_s^{t;d}$, i.e., $f_s^{t;d}(x,\cdot) = \widehat{f}_s^{t;d}(x,\cdot)$, $\forall x\in\mathbb{R}^d$, $\mathbb{P}$-a.s.

            \item[(c)] RanNN approximator $\widehat{f}_s^{t;d}$ in (b) satisfies \( \mathbb{E}[\size(\widehat{f}_s^{t;d})] \leq \bar{c}d^{\bar{q}} \).
        \end{enumerate}
    \end{enumerate}
\end{assumption}

\begin{remark}
    The above assumption covers a large class of models used in financial markets. We further discuss concrete examples in Subsection~\ref{subsubsec:affine_ito} and Remark~\ref{rem:discussion_of_expression_rate_model}.
\end{remark}

Similar to \citep{gonon23}, we make the following assumption on the obstacle (payoff) function.

\begin{assumption}[Assumption on payoff]\label{ass:g_new}
    There exist constants $c, q, r >0 $ independent of $d$, such that for any $\varepsilon>0, n \in \N $, there exists a deep ReLU network $\widehat{g}^{d,\varepsilon}_n : \mathbb{R}^d \rightarrow \mathbb{R}$ that satisfies
    \begin{itemize}
        \item[(i)]  \( |\widehat{g}^{d,\varepsilon}_n(x) - g^d(t_n,x)|  \le \varepsilon c d^q(1+\|x\| ) \) for all \( x\in \R^d  \) ;
        \item[(ii)] \( \Lip(\widehat{g}^{d,\varepsilon}_n ) ,\;  |g^d(t_n,0)| \le c d^{q } \) and \( \size(\widehat{g}^{d,\varepsilon}_n) \le c d^{q}\varepsilon^{-r} \).
    \end{itemize}
\end{assumption}

\begin{remark}
    The constants $c, q$ in Assumptions~\ref{ass:dynamic_ass} and \ref{ass:g_new} can be chosen to coincide by taking maxima if necessary.
    By inspection, Assumption~\ref{ass:N_0_structural_ass_g_deter} for $\{ g(t_n,\cdot) \}_{n=0}^N $ is naturally satisfied under Assumption~\ref{ass:g_new}, which in turn ensures the expressivity of the numerical integration scheme.
\end{remark}

\subsubsection{Expressivity of the value function approximation}
Expressivity for the value function approximation of the optimal stopping problem has been studied in \citep{gonon23}, which extends the Black--Scholes PDE setting developed in \citep{Jentzen23}. Both references provide useful tools and intuition for our discrete monitoring points with continuous-time observation setting, but they do not directly cover the present framework. We therefore extend their methodology to obtain an expressivity result for the value function deep ReLU approximation, which forms the basis of our subsequent \textit{DeepMartingale} expressivity analysis.

To align with our \textit{DeepMartingale} technical setup, we introduce the finite Borel measures for all $n\in \N^{-1}$ and $K \in \bbN_{+}$, that is 
\[ \tilde{\rho}^{K;d}_{n+1}(A):= \mathbb{E}\Bigl[ \sum_{k=0}^{K-1}\|\Delta W^d_{t^n_k}\|^2 \mathbf{1}_{A}( X^d_{t_{n+1}} ) \Bigr] ,\; \; A \in \mathcal{B}(\mathbb{R}^d)  .  \]
The following lemma ensures that expressivity can be evaluated under $ \tilde{\rho}^{K;d}_{n+1} $.

\begin{lemma}\label{lem:prob_measure_bound}
Under Assumption~\ref{ass:dynamic_ass} with order $p > 0 $ and $ \|x^d_0\| \le cd^q $ (same constants in the above assumption), let
define $ \bar{\rho}^{K;d}_{n+1} = (\frac{N}{T}d^{-1}) \tilde{\rho}^{K;d}_{n+1} $ for all $ n\in \N^{-1}  $, $ K\in \bbN_{+} $. Then, $ \bar{\rho}^{K;d}_{n+1} $ are probability measures on $ \calB(\R^d) $, and for any $ \bar{p} \in (0,p) $, there exists positive constants $ \widehat{k}_{n+1} , \widehat{p}_{n+1} $ independent of $d, K$, such that \(  \M_{\bar{p}}(\bar{\rho}^{K;d}_{n+1}) \le  \widehat{k}_{n+1} d^{\widehat{p}_{n+1}} \) for all $ K \in \bbN_{+} , n\in \N^{-1} $.
\end{lemma}

We now state our main theorem on deep ReLU approximation of the value function under discrete monitoring points with continuous-time observation. Importantly, the result holds uniformly for an arbitrary $K \in \N_{+} $.

\begin{theorem}[Value Function Approximation with Expressivity]
\label{theorem:neural_approx_V_new}
    Under Assumption~\ref{ass:dynamic_ass} with order $p>2 $ and Assumption~\ref{ass:g_new}. Moreover, suppose $ \|x^d_0\| \le cd^q $ (same constants in above assumptions). Then, for any $n \in \N^{-1} $, there exist constants $c_{n+1},q_{n+1},\tau_{n+1} \ge 1 $ independent of $d , K $, such that for any $\varepsilon>0$, there exists a deep ReLU network $\widehat{V}^{d,\varepsilon}_{n+1} $ satisfies for all $K\in \mathbb{N}_{+} $,
    \[
      \big\| \widehat{V}^{d,\varepsilon}_{n+1} - V^d_{n+1}\|_{2;\tilde{\rho}^{K;d}_{n+1}} \le \varepsilon  ,\; \;  \text{and} \; \; \size(\widehat{V}^{d,\varepsilon}_{n+1}),\;  \Growth(\widehat{V}^{d,\varepsilon}_{n+1}) \le c_{n+1}d^{q_{n+1}} \varepsilon^{-\tau_{n+1}} .
    \]
\end{theorem}

The proof of Theorem~\ref{theorem:neural_approx_V_new} is given in the Appendix. It follows as a direct corollary of the general statement below (Theorem~\ref{theorem:recursive_express_new}) together with Lemma~\ref{lem:prob_measure_bound}. Theorem~\ref{theorem:recursive_express_new} can be viewed as an extension of \citep{Jentzen23} and \citep{gonon23} to the setting of discrete monitoring with continuous-time observation. Since it directly follows the elegant proof procedure of \citep[Theorem~3.6]{gonon23}, we omit it in this paper.

\begin{theorem}[Value Function Approximation with Expressivity, General Form]
\label{theorem:recursive_express_new}
Under Assumption~\ref{ass:dynamic_ass} with order $p>2 $ and Assumption~\ref{ass:g_new}, for any given $\bar{p} \in (2,p] $, $k_1,p_1 \ge 1 $ independent of $d,K$ with sequences $k_{n+1} := c(1+k_{n}) $, $p_{n+1}=p_n+q$, $n \in \N^{-1} $ ($ c,q $ are from above assumptions), there exist constants $c_{n+1}, q_{n+1}, \tau_{n+1} \ge 1$ independent of $ d,K $, $ \forall \; n \in \N^{-1}$, such that for any family of probability measures $(\rho_{n+1})_{n=0}^{N-1}$ on $ \mathcal{B}(\mathbb{R}^d) $ satisfying \( \M_{\bar{p}}(\rho_{n+1}) \le k_{n+1}d^{p_{n+1}}  \),
and for any $\varepsilon>0 $, there exists a deep ReLU network $\widehat{V}^{d,\varepsilon}_{n+1}:\mathbb{R}^d \rightarrow \mathbb{R} $ such that
\[
   \big\| \widehat{V}^{d,\varepsilon}_{n+1} -  V^d_{n+1} \big\|_{2;\rho_{n+1}} \le \varepsilon   ,\; \;  \text{and} \; \; \size(\widehat{V}^{d,\varepsilon}_{n+1}) ,\;  \Growth(\widehat{V}^{d,\varepsilon}_{n+1}) \le c_{n+1}d^{q_{n+1}}\varepsilon^{-\tau_{n+1}}  .
\]
\end{theorem}

\subsubsection{Expressivity of integrand approximation}
Given the above value function approximation, we next construct an approximator for the integrand $\widehat{Z}^{*,K,n;d} $, motivated by Remark~\ref{remark:priori_approx}.

We begin with the following technical lemma, which extends \cite[Lemma~4.10]{gonon23}.

\begin{lemma}
    \label{lemma:prob_new}
    Let $U$ be a non-negative random variable. Given $N_1\in \mathbb{N} $ and the non-negative integer sequence $\{J_n \}_{n=1}^{N_1} $, then for any $n=1,\ldots,N_1 $, $X_n^{i}, i=1,\ldots,J_n $ are i.i.d random variables. Suppose $\mathbb{E}[U ]\le M_0,\; \mathbb{E}|X^1_n| \le M_n ,\; n =1,\ldots,N_1 $, and $M_n>0 $ for all $n$. Then, let 
    \[
        B_{N_1} := \prod_{n=1}^{N_1} \Big\{ \max\limits_{i=1,\ldots,J_n}|X^i_n| \le  (N_1+2) J_nM_n \Big\},
    \]
    we have
    \[
        \mathbb{P}  (U \le  M_0 ) >0 ,\quad  \mathbb{P}\big(  U\le (N_1+2) M_0 ,\; B_{N_1} \big) > 0  .
    \] 
\end{lemma}

We now construct a joint approximation of $\widehat{Z}^{*,K,n;d} $ by a family of deep ReLU networks, one at each observation time $t^n_k$.

\begin{theorem}[Joint Network Approximation for Integrand] \label{theorem:z_express_approx_new}
    Under Assumption~\ref{ass:stronger_ass} with order $\tilde{p} > 2 $, Assumption~\ref{ass:dynamic_ass} (2).(a) with order $p\ge 2 $ and $\|x^d_0\| \le cd^q $, and Assumption~\ref{ass:g_new}, suppose that for any $n \in \N^{-1} $, there exist constants $c_{n+1},q_{n+1},\tau_{n+1} \ge 1$ independent of $d , K $, such that for any $\varepsilon>0$, there exists a deep ReLU network $\widehat{V}^{d,\varepsilon}_{n+1}: \mathbb{R}^d \rightarrow \mathbb{R} $ that satisfies for all $K\in \mathbb{N}_{+} $,
    \[
       \big\| \widehat{V}^{d,\varepsilon}_{n+1} - V^d_{n+1}\|_{2;\tilde{\rho}^{K;d}_{n+1}} \le \varepsilon ,\; \;  \text{and} \; \; \size(\widehat{V}^{d,\varepsilon}_{n+1}),\;  \Growth(\widehat{V}^{d,\varepsilon}_{n+1}) \le c_{n+1}d^{q_{n+1}} \varepsilon^{-\tau_{n+1}} .
    \]
    Then, for any $n \in \N^{-1} $, there exist constants $\widehat{c}_n,\widehat{q}_n,\widehat{\tau}_n,\widehat{m}_n \ge 1 $ independent of $d$, such that for any $\varepsilon>0$, $K\in \mathbb{N} $, there exist a family of deep ReLU networks $\gamma^{n;d,\varepsilon,K}_k: \mathbb{R}^d\rightarrow \mathbb{R}^d$, $k \in \K $ and their joint function (connecting by spline) $\widehat{z}^{d,\varepsilon,K}_n(t,x) $ with $\widehat{z}^{d,\varepsilon,K}_n(t^n_k,\cdot)\equiv \gamma^{n;d,\varepsilon,K}_k $ satisfying
    \[
          \big\| \widehat{z}^{d,\varepsilon,K}_n -   \widehat{Z}^{*,K,n;d} \big\|_{2;\mu^{K;d}_n } \le \varepsilon  ,\; \; \text{and} \; \; \Growth( \gamma^{n;d,\varepsilon,K}_k) ,\; \size(\gamma^{n;d,\varepsilon,K}_k) \le  \widehat{c}_n d^{\widehat{q}_n}  \varepsilon^{-\widehat{\tau}_n} (K)^{\widehat{m}_n} ,\; \forall \; k\in \K .
    \]
\end{theorem}

To realize the connecting splines as a single deep ReLU network, we use the following lemma.

\begin{lemma}\label{lem:time_indicator_realization}
    Given $ n\in \N^{-1} $, suppose $ h^{n;K}_k : \mathbb{R} \to \R_{+} ,\; k\in \K $ are tent (triangular) functions on $ \T^{n,K} $ ($   h^n_k \in [0,1] $ on $ [t_n,t_{n+1}] $) that connecting all observation time points (see Figure~\ref{fig:tent_function}, $ h^{n;K}_0 $ is ReLU). Then, $ h^{n;K}_k $ can be realized by deep ReLU networks with $ \size(h^{n;K}_0) =3 $ and $ \size(h^{n;K}_k) = 12 $ for all $ k\in \K \setminus \{ 0 \} $.
\end{lemma}

\begin{figure}[H]
\centering
\footnotesize
\begin{tikzpicture}
\begin{axis}[
    width=6.12cm, height=2.72cm,
    scale only axis,
    axis lines=left,
    axis line style={->},
    clip=false,
    xmin=-0.6, xmax=11.1, 
    ymin=-0.05, ymax=1.30,
    xlabel={$t$},
    ylabel={},
    xtick={0,1.5,3,4.5,6,7.5,9,10.5},
    xticklabels={$t_0^n(t_n)$,$t_1^n$,$t_2^n$,$t_3^n$,$t_4^n$,$t_5^n$,$t_6^n$,$t_7^n(t_{n+1})$},
    ytick={0,1},
    tick label style={font=\footnotesize},
    label style={font=\footnotesize},
    every axis plot/.append style={line join=round, line cap=round},
    grid=both,
    grid style={line width=0.2pt, draw=gray!25},
]

\node[font=\footnotesize, anchor=south]
    at (axis description cs:-0.16,0.50) {$h_k^{n;7}(t)$};

\addplot[very thick, red!75!black] coordinates {(-0.6,0) (0,0)};      
\addplot[very thick, red!75!black] coordinates {(-0.5,1.33) (1.5,0)};       
\addplot[very thick, red!75!black] coordinates {(1.5,0) (11.1,0)};    

\addplot[only marks, mark=*, mark size=2.0pt, draw=red!75!black, fill=red!75!black]
    coordinates {(0,1)};

\node[font=\footnotesize, anchor=south west, text=red!75!black]
    at (axis cs:0.15,1.02) {$h_0^{n;7}(t_0^n)=1$};

\node[
    font=\footnotesize,
    anchor=south west,
    text=red!75!black,
    fill=white,
    fill opacity=0.85,
    text opacity=1,
    inner sep=1pt
] at (axis cs:2.10,0.16) {$h_0^{n;7}(t)$};

\addplot[very thick, blue!70!black] coordinates {(-0.6,0) (4.5,0)};   
\addplot[very thick, blue!70!black] coordinates {(4.5,0) (6,1) (7.5,0)}; 
\addplot[very thick, blue!70!black] coordinates {(7.5,0) (11.1,0)};   

\addplot[only marks, mark=*, mark size=2.0pt, draw=blue!70!black, fill=blue!70!black]
    coordinates {(6,1)};

\node[font=\footnotesize, anchor=south west, text=blue!70!black]
    at (axis cs:6.05,1.02) {$h_4^{n;7}(t_4^n)=1$};

\node[font=\footnotesize, anchor=south east, text=blue!70!black]
    at (axis cs:10.95,0.12) {$h_4^{n;7}(t)$};

\addplot[
    only marks, mark=*, mark size=1.3pt,
    draw=black!40, fill=black!40
] coordinates {(0,0) (1.5,0) (3,0) (4.5,0) (6,0) (7.5,0) (9,0) (10.5,0)};

\end{axis}
\end{tikzpicture}
\caption{Tent function $ h^{n;K}_k $ on \( \T^{n,K} \) (\(K=7\)).}
\label{fig:tent_function}
\end{figure}

Using Lemma~\ref{lem:time_indicator_realization}, we obtain the expressivity result for approximating the integrand process by a single deep ReLU network under $L^2_{1+d,d}(\mu^{K;d}_n ) $ as follows.

\begin{theorem}[Single Integrand Network Realization \& Approximation]
    \label{theorem:joint nerual network realization_new}
    Under Assumption~\ref{ass:dynamic_ass} with order $p> 2 $ and Assumptions~\ref{ass:stronger_ass} with order $\tilde{p}> 4 $ and Assumption~\ref{ass:g_new}, for any $n \in \N^{-1} $, there exist constants $\bar{c}_n,\bar{q}_n,\bar{\tau}_n,\bar{m}_n \ge 1 $ independent of $d $, such that for any $\varepsilon \in (0,1]$, $K\in \mathbb{N}_{+} $, we have a deep ReLU network $\tilde{z}^{d,\varepsilon,K}_n : \mathbb{R}^{1+d}\rightarrow \mathbb{R}^d $ that satisfies for all \( t\in [t_n,t_{n+1}] \),
    \[
        \big\| \tilde{z}^{d,\varepsilon,K}_n - \widehat{Z}^{*,K,n;d} \big\|_{2;\mu^{K;d}_n } \le \varepsilon , \; \; \text{and} \; \; \Growth(\tilde{z}^{d,\varepsilon,K}_n(t,\cdot)) ,\; \size(\tilde{z}^{d,\varepsilon,K}_n ) \le \bar{c}_n d^{\bar{q}_n} \varepsilon^{-\bar{\tau}_n}(K)^{\bar{m}_n} .
    \]
\end{theorem}

\begin{remark}[Expressivity of Deep Delta Hedge]
    Theorem~\ref{theorem:joint nerual network realization_new} can also be interpreted as an expressivity result for the ``deep'' \textit{delta hedge} induced by \textit{DeepMartingale}. This is not restricted to Bermudan options, but also applies to European options. In particular, European options can also be delta hedged without the curse of dimensionality when the market dynamics and payoff exhibit sufficiently ``nice'' structures.
\end{remark}

\subsubsection{Expressivity of DeepMartingale}

We now present our main result in Theorem~\ref{thm:express_deep_mtg}. It shows that \textit{DeepMartingale} achieves an expressive solution to the dual formulation of the optimal stopping problem under discrete monitoring points with continuous-time observation: the required computational resources grow at most polynomially in the dimension $d$ and the target accuracy $\varepsilon$ of Problem~\ref{pb:backward_minimization}, as quantified by the size of the deep ReLU approximators. 

\begin{theorem}[Expressivity of DeepMartingale]\label{thm:express_deep_mtg}
    Under Assumption~\ref{ass:N_0_structural_ass_model_deter}, Assumption~\ref{ass:dynamic_ass} with order $p> 2 $, Assumption~\ref{ass:g_new}, and Assumption~\ref{ass:stronger_ass} with order $\tilde{p}> 4 $, then there exist constants $\tilde{c},\tilde{q},\tilde{r} > 0 $ independent of $d$, such that for any $\varepsilon \in (0,1] $, there exist $ K_{d,\varepsilon} \le \tilde{c}d^{\tilde{q}} \varepsilon^{-1} $, deep ReLU networks $\tilde{z}^{d,\varepsilon}_n :\mathbb{R}^{1+d} \rightarrow \mathbb{R}^d $, $n \in \N^{-1} $ and $\widetilde{M}^{d,\varepsilon}_n:=\sum_{m=0}^{n-1}\sum_{k=0}^{K_{d,\varepsilon} -1}\tilde{z}^{d,\varepsilon}_m(t^m_k,X^d_{t^m_k} )\cdot \Delta W^d_{t^m_k}$, $n \in \N^{-1} $ that satisfy for any $n \in \N^{-1} $, $ t\in [t_n,t_{n+1}] $, 
    \[
        \big\| \widetilde{U}^d_n(\widetilde{M}^{d,\varepsilon}  ) - Y^{*;d}_n \big\|_{L^2}  \le (N-n)\varepsilon , \; \; \text{and} \; \; \size(\tilde{z}^{d,\varepsilon}_n) ,\; \Growth(\tilde{z}^{d,\varepsilon}_n(t,\cdot))  \le \tilde{c}d^{\tilde{q}}\varepsilon^{-\tilde{r}} 
    \]
    and 
    Moreover, it directly implies $ 0 \le \mathbb{E}[\widetilde{U}^d_n(\widetilde{M}^{d,\varepsilon})] - \mathbb{E}[Y^{*;d}_n] \le (N-n)\varepsilon  $. If $ Y^{*;d}_n \ge 0 ,\; n \in \N $, $ \mathbb{P} $-a.s., then $ 0 \le  \big(\mathbb{E}|\widetilde{U}^d_n(\widetilde{M}^{d,\varepsilon})|^2 \big)^{\frac{1}{2}} - \big(\mathbb{E}|Y^{*;d}_n|^2 \big)^{\frac{1}{2}} \le (N-n)\varepsilon $. 
\end{theorem}

The proof of Theorem~\ref{thm:express_deep_mtg} follows by directly combining the approximation results established above; we provide the details in the Appendix.

\begin{remark}
    Theorem~\ref{thm:express_deep_mtg} justifies that Problem~\eqref{pb:backward_minimization} can be solved by \textit{DeepMartingale} without the curse of dimensionality:
    \begin{itemize}
        \item \emph{Theoretical implication:} key computational resources, such as the relevant neural network sizes, can be chosen to \textbf{scale at most polynomially in $d$}, while retaining the approximation capability for the dual upper bound/second moment and the solvability of Problem~\eqref{pb:backward_minimization};
        \item \emph{Practical implication:} one may first solve low-dimensional instances and estimate the polynomial growth order in $d$ (as a function of the error level $\varepsilon$) via regression, and then use this estimate to guide high-dimensional computations. We also call it \textit{dimension scaling law}. See Section~\ref{sec:Implement} for details.
    \end{itemize}
\end{remark}

\subsubsection{Example: Affine It\^o diffusion}\label{subsubsec:affine_ito}
We illustrate our structural framework and the resulting \textit{DeepMartingale} expressivity guarantees using a widely adopted class of models---affine It\^o diffusions. This class covers many dynamics used in practice (e.g., the Black--Scholes model and Ornstein--Uhlenbeck processes), thereby demonstrating that our main expressivity statement is relevant for realistic applications.

We first recall the affine It\^o diffusion setting from \citep{Jentzen23}.

\begin{definition}[Affine It\^o Diffusion]\label{def:AID}
If $X^d$ satisfies \eqref{eq:SDE} and the coefficient functions $a^d:\mathbb{R}^{d} \rightarrow \mathbb{R}^{d} $, $b^d:\mathbb{R}^d \rightarrow \mathbb{R}^{d\times d} $ satisfy for all $x,y\in \mathbb{R}^d$, $\lambda \in \mathbb{R} $,
\[
    a^d(\lambda x+y) + \lambda a^d(0)=\lambda a^d(x) + a^d(y) ,\; \;
        b^d(\lambda x+y) + \lambda b^d(0)=\lambda b^d(x) + b^d(y) ,
\]
   then we call $X^d$ an affine It\^o diffusion (AID). 
\end{definition}

To match our structural framework and derive \textit{DeepMartingale} expressivity, we impose the following expressivity conditions on AID models.

\begin{definition}[AID with $\frac{1}{2}$-$\log $ Growth ]\label{def:ADI_express}
If $X$ follows Definition~\ref{def:AID} and there exist constants $C^{*},Q^{*} $ such that for any $x\in \mathbb{R}^d $,  
\begin{equation*}
    \Lip a^d \le C^{*}(\log d)^{\frac{1}{2}}, \;  \ \LipH b^d \le C^{*}(\log d)^{\frac{1}{4}},\quad \|a^d(0)\|,\  \|b^d(0)\|_{\H} \le C^{*} d^{Q^{*}},
\end{equation*}
and in addition, $\|x^d_0\|\le C^{*}d^{Q^{*}} $, then we call $X^d$ AID with $\frac{1}{2}$-$\log $ Growth (AID-log).
\end{definition}

\begin{remark}\label{rem:discussion_of_expression_rate_model}
Definition~\ref{def:ADI_express} accommodates a broad class of models commonly used in practice. Indeed, AID dynamics admit the representation (see Lemma~\ref{lemma:equiv_affine_new})
\[
dX^d_t = (A^1_d X^d_t + b^1_d)\,dt + \sum_{i=1}^d ( A^{2,i}_d X^d_t    +   b^{2,i}_d)\, dW^{i;d}_t   ,   
\]
where $b^1_d \in \mathbb{R}^d$, $A^1_d, b^2_d = (b^{2,1}_d,\ldots,b^{2,d}_d) \in \mathbb{R}^{d\times d}$, and $A^2_d = (A^{2,1}_d,\ldots,A^{2,d}_d) \in \mathbb{R}^{d\times d\times d}$. The AID-log condition then requires $\|A^1_d\|_{2} , \|A^2_d\|^2_{2}:= \sum_{i=1}^d \|A^{2,i}_d\|^2_2  \le C^{*}(\log d)^{\frac{1}{2}} $ and $\|b^1_d\|, \|b^2_d\|_{\H} \le C^{*} d^{Q^{*}}$. Several widely used models satisfy Definition~\ref{def:ADI_express}:
\begin{enumerate}
    \item \textbf{Geometric Brownian motion:} $ dX^d_t = \mathrm{diag}(X^d_t)\bigl(\mu^d \, dt + \mathrm{diag}(\sigma^d)\, dW^d_t\bigr) $ with $\mu^d, \sigma^d \in \mathbb{R}^d$ satisfying $\|\mu^d\|_{\infty} , \|\sigma^d\|^2_{\infty}  \le  C (\log d)^{\frac{1}{2}} $. 
    This condition allows the components of $ \mu^d $ and $ \sigma^d $ to grow at rates $ \mathrm{O}( (\log d)^{\frac{1}{2}} ) $ and $ \mathrm{O}( (\log d)^{\frac{1}{4}} ) $, respectively. In particular, it covers all experimental GBM settings considered in \citep{Becker19}.
    \item \textbf{Ornstein--Uhlenbeck process:} $dX^d_t = (\mu^d - X^d_t)\, dt + \Sigma^d \, dW^d_t$ with $\mu^d \in \mathbb{R}^d$ and $\Sigma^d \in \mathbb{R}^{d\times d}$ satisfying $\|\mu^d\|_{\infty}, \|\Sigma^d\|_{\H} \le C d^Q$. For the constant volatility case $\Sigma^d = \sigma \mathbf{I}_d$ with $\sigma \ge 0 $, we have $\|\Sigma^d\|_{\H} = \sigma \|\mathbf{I}_d\|_{\H} = \sigma d^{\frac{1}{2}}$, which satisfies the above condition.
\end{enumerate}
\end{remark}

Under Definition~\ref{def:ADI_express}, the structural framework proposed above applies to AID-log models as follows.

\begin{lemma}\label{lemma:AID-log_ass1}
    If $X^d $ is an AID-log (Definition \ref{def:ADI_express}), then $X^d $ satisfies Assumption~\ref{ass:N_0_structural_ass_model_deter}, Assumption~\ref{ass:dynamic_ass} for any order $p>2 $ and Assumption~\ref{ass:stronger_ass} for any order $\tilde{p} > 4 $.
\end{lemma}

Consequently, the expressivity of \textit{DeepMartingale} for AID-log dynamics follows immediately.

\begin{theorem}[Expressivity for DeepMartingale: AID-log]\label{thm:AID-log_express}
    If $X^d $ is an AID-log and $g^d $ satisfies Assumption~\ref{ass:g_new}, then the expressivity for DeepMartingale (Theorem~\ref{thm:express_deep_mtg}) holds.
\end{theorem}

\section{Numerical implementation}\label{sec:Implement}

\subsection{DeepMartingale (pure dual) algorithm, delta hedge, and dimension-scaling approach}
In this section, we numerically illustrate the convergence, stability, and scalability (expressivity) of \textit{DeepMartingale} for the dual formulation of optimal stopping.

We include a comparison with the recent primal--dual hedging framework of \cite{guo2024simultaneousupperlowerbounds,puredual-mf,yang2024deepprimaldualbsdemethod}. In our experiments and under comparable computational budgets, \textit{DeepMartingale} delivers accurate upper bounds and stable hedging performance at moderate to high dimensions, while the primal--dual baseline may become less stable as the dimension increases. A plausible explanation is that our method is a \emph{pure dual} approach in the spirit of \cite{roger10,puredual-mf}: we directly optimize over a parametrized class of martingales and derive a hedging strategy without relying on any primal approximation. By contrast, primal--dual approaches may inherit approximation errors from the primal component, which can propagate to the dual bound and hedging performance in high dimension.

We first present the training procedure. We then describe our evaluation protocol and the construction of the corresponding ``deep'' delta hedge strategy. Finally, we propose a regression-based dimension-scaling law estimation, motivated by our expressivity theory.

All DNN computations are performed in single precision (float32) on an NVIDIA A100 GPU (1095 MHz core clock, 40 GB memory) with dual AMD Rome 7742 CPUs, running PyTorch 2.3.0 on Ubuntu 20.04.4 LTS. Other computations except for DNN are conducted on Apple Silicon M4 Pro CPU with 64 memory. 

\subsubsection{Training}
Let \(I\) denote the network depth and let all hidden layers have width \(q\). Let \(M\) be the number of training epochs. At each epoch, we generate \(J\) fresh sample paths to update the parameters, which empirically improves training efficiency and reduces memory usage in high dimension.
Given the SDE coefficients \((a,b)\) in \eqref{eq:SDE} and the payoff function \(g\), the training procedure is summarized in Algorithm~\ref{Alg:1}. We initialize all networks using Xavier normal initialization and train with the Adam optimizer.

\begin{algorithm}
\caption{DeepMartingale Training \label{Alg:1}}
\vskip6pt
\begin{algorithmic}
\Procedure{DeepPureDual}{$d,T,N,K,I,q,J,M,a,b,g,x_0$}
    \State Simulate \(J\) sample paths \((X,W)\) according to \eqref{eq:SDE} with step size \(\Delta t = T/(N K)\). Initialize \textit{DeepMartingale} \(M^{\theta,K}\), \(\theta=(\theta_n)_{n=0}^{N-1}\), defined in \eqref{deep_mart}, using depth \(I\) and width \(q\).
    \For{$m \gets 0$ \textbf{to} $M$}
        \State Set \(U_N \equiv g(t_N,X_{t_N})\) and initialize \(U_n\) for \(n\in \N^{-1} \).
        \For{$n \gets N-1$ \textbf{to} $0$}
            \If{$g$ non-negative}
                \State Update \(\theta_n\) by solving \eqref{eq_optim_second_moment} in Problem~\ref{pb:backward_minimization}.
            \Else
                \State Update \(\theta_n\) by solving \eqref{eq_optim} in Problem~\ref{pb:backward_minimization}.
            \EndIf
            \State Update \(U_n\) via the one-step recursion \eqref{eq:backward_recursion}.
        \EndFor
    \EndFor
    \State \textbf{return} \(\theta\)
\EndProcedure
\end{algorithmic}
\end{algorithm}

\subsubsection{Evaluation and delta hedge}\label{subsubsec:numerical_delta_hedge_evaluation}
After training, we estimate the mean and standard deviation of the upper bound using an independent evaluation set. Specifically, we generate \(J_1\) new sample paths and compute \(U_0\) by setting \(U_N \equiv g(t_N,X_{t_N})\) and applying the backward recursion \eqref{eq:backward_recursion}. To mitigate memory pressure, we compute the statistics by iterating this evaluation over \(M_1\) batches and aggregating the results.
We report the out-of-sample upper bound by
\[
\widehat{U}_0 := \mathbb{E}[U_0] .
\]

Under the complete-market assumption (i.e., invertible \(b(t,x)\)), Remark~\ref{rem:delta_hedge_deepMart} yields a ``deep'' delta hedge on each monitoring interval \([t_n,t_{n+1}]\), \(n=0,\ldots,N-1\), via
\[
\vartheta_n^{\theta_n}(t,x) := b^{-1}(t,x)\, z_n^{\theta_n}(t,x).
\]
Since in high-dimensional settings it is typically infeasible to obtain the true optimal stopping rule, we evaluate the hedging strategy using the distribution of the \emph{worst-case} (over monitoring times) hedging error under rebalancing frequency \(K\), following the metric also used in \cite{guo2024simultaneousupperlowerbounds}:
\begin{equation}
\mathrm{H}_{err} := V_0 + \min_{n=0,\ldots,N}\Big(M_n^{\theta,K} - g(t_n,X_{t_n})\Big),
\label{eq:hedging_error_worst}
\end{equation}
where \(V_0\) is the option price used to initialize the hedging portfolio.
By the dual representation \eqref{eq_strong_dual}--\eqref{eq_sure_optimal}, \(\mathrm{H}_{err}\) captures the (pathwise) slack of the hedging strategy against an adversarial exercise decision: \(\mathrm{H}_{err}\ge 0\) corresponds to a pathwise super-hedge, while \(\mathrm{H}_{err}<0\) indicates a shortfall. This worst-case criterion is robust and practically relevant, since real exercise behavior may deviate from the model-implied optimal stopping rule. Here we set $ V_0 $ as the benchmark price $ U_{\mathrm{ref}}^d $ in literature.

\subsubsection{Regression-based dimension scaling law: estimating an empirical expressivity order}\label{subsubsec:dimension-scaling-law}
To quantify how the accuracy of the learned \textit{DeepMartingale} \(M^{\theta,K}\) scales with the dimension \(d\), we conduct a regression-based estimation motivated by our expressivity result (Theorem~\ref{thm:express_deep_mtg}). For a set of moderate dimensions (e.g., \(2\le d\le 10\)), we compare the evaluated upper bound \(\widehat{U}_0^d\) against reference values \(U_{\mathrm{ref}}^d\) from the literature or high-accuracy benchmarks, while keeping the numerical and training configuration fixed across dimensions (e.g., \(K,I,q,M\) do not depend on \(d\)). We define the absolute error
\[
\mathrm{err}(d):=\widehat{U}_0^d - U_{\mathrm{ref}}^d,
\]
and postulate a power-law scaling \(\mathrm{err}(d)\approx C d^{p}\) for constants \(C>0\) and \(p>0\). Taking logarithms yields
\begin{equation}
\log(\mathrm{err}(d)) \approx p\log(d)+c , \; \;  c=\log C.
\label{eq:log-transform-scaling}
\end{equation}
We estimate \(p\) via ordinary least squares on the log-log data, and report the coefficient of determination \(R^2\) to quantify goodness of fit. The slope \(p\) is interpreted as an \emph{empirical expressivity order} with respect to the dimension: smaller \(p\) indicates better robustness to increasing \(d\), while larger \(p\) suggests faster degradation. In particular, \(p\approx 1\) corresponds to an approximately linear error growth in \(d\).



\subsection{Experiments - Bermudan max-call}\label{subsec:experiments}


We consider well-studied benchmarks to assess the performance of \textit{DeepMartingale} and the associated hedging strategy in high dimensions. We use unbounded ReLU activation, network depth $ I=3 $, training and evaluation batch sizes $ J = J_1 = 8,192 $, and the number of evaluation batches $ M_1 = 200 $. To study the scaling with respect to the dimension, we fix the number of integration points (rebalancing frequency) $ K_{\mathrm{ori}} =32 $ (as in \cite[Section~5.1.2]{guo2024simultaneousupperlowerbounds}), the network width $ q_{\mathrm{ori}} = 40 $, and the number of training epochs $ M_{\mathrm{ori}} = 5,000 $.   

Following \citep{andersen04,broadiecao08,schoen13,Becker19, guo2024simultaneousupperlowerbounds, puredual-mf}, we specify the dynamics by setting $ a(t,X_t) =(r - \delta  ) X_t $ and $b(t,X_t) =  \mathrm{diag}(X_t)  \Sigma  $, and we take the discounted payoff function to be of max-call type
\[
    g(t,X_t) = e^{-rt} \bigl(\max_{1\le i \le d}(X^i_t) - S \bigr)^{+} ,
\]
where $r,\delta, \Sigma := \mathrm{diag}(\sigma_1,\ldots,\sigma_d), S $ denote the risk-less interest rate, dividend rate, volatility matrix, and strike price, respectively.

We study\footnote{All the codes for replicating the experiments can be found in \url{https://github.com/GEOR-TS/DeepMartingale}.} the following \textbf{asymmetrical setting}:
$S=100$, $x_0=s_0 \mathbf{1}_d $, $r=0.05 $, $\delta=0.1 $, and
\[
    \sigma_i = 
    \begin{cases}
         0.08 + 0.32 \times \frac{i-1}{d-1}  & \text{if} \; \; d\le 5  \\
         0.1 + \frac{i}{2d} & \text{if} \; \; d>5 
    \end{cases}
     \; \; , \; i=1,\ldots,d .
\]
In addition, we also evaluated the symmetrical case and obtained qualitatively similar numerical results. To keep the presentation focused and to better highlight the advantages of \textit{DeepMartingale}, we omit the symmetrical case in this paper.

To illustrate the high-dimensional performance of our method for solving the dual optimal stopping problem and for constructing the associated delta hedging strategy, we also implement the algorithms in \cite[\textit{Guo-DeepPrimalDual}]{guo2024simultaneousupperlowerbounds} (source code is provided in their paper) and \cite[\textit{Alfonsi-PureDual}]{puredual-mf}.
We do not directly compare against the seminal primal--dual methods in \cite{andersen04}, since our primary focus is on the high-dimensional setting. We also implemented the recent deep primal--dual approach in \cite[\textit{Yang-DeepPrimalDual}]{yang2024deepprimaldualbsdemethod} (source code is provided in their paper). However, their method is mainly designed for high-frequency monitoring with aligned observations. In our more frequent observation setting, we adapted their code to incorporate frequent observations, but using the hyperparameter configuration reported in their paper leads to out-of-memory issues on our NVIDIA A100 (40G memory) due to the $K$-times larger data size, even for $d=2$. Hence, we are unable to include a numerical comparison with their algorithm. We also do not include a direct comparison against the upper bound methods in \cite{andersen04,Becker19}, since those approaches derive upper bounds via nested simulation after solving the primal problem; while highly accurate, such bounds are not explicitly realizable through a hedging strategy.

For \textit{Guo-DeepPrimalDual}, we first run the original code for all dimensions without modifying either the network architecture or the training procedure. After deriving the scaling law of \textit{DeepMartingale}, we further adjust their implementation by setting the number of substeps (i.e., the observation frequency) to match that of \textit{DeepMartingale}, and by choosing network widths with the same growth rate, so that the overall training setup is scaled to the same order as closely as possible. We train their method with the quadratic term, as in their code, which yields the best performance. We then evaluate the resulting upper bound under the full martingale (with the quadratic term) and under the delta martingale (delta hedging only), together with the corresponding worst-case delta hedging error. We do not consider the quadratic term in hedging since it is difficult to realize this direct quadratic hedging portfolio by the trained martingale.

For \textit{Alfonsi-PureDual}, we follow the algorithmic procedure in their paper. Since their best upper-bound performance is achieved by using polynomial basis functions and by enriching the basis with Vanilla options with respect to the underlying assets, we adopt the same setup, while modifying the observation frequency to $K$ so as to match the scaled configuration of \textit{DeepMartingale}. We only report the scaled performance, since in their implementation the polynomial basis dimension grows with the state dimension. We then compute the corresponding upper bound under the full portfolio (including Vanilla options) and under the delta portfolio (delta hedging only), as well as the associated worst-case delta hedging error.


Finally, for a fair comparison, all evaluation settings (in particular, evaluating the upper bound and hedging with $M_1 \times J_1$ samples) are chosen to be the same as those used for \textit{DeepMartingale}.

\subsubsection{Expressivity -- dimension scaling law estimation}


In this part, we fix the configuration to be constant $ K_{\mathrm{ori}}, q_{\mathrm{ori}}, M_{\mathrm{ori}} $ across dimensions from $ d=2 $ to $10$. The benchmarks for computing the absolute error are taken to be the reference values in \cite{Becker19}. We regress the error data $ \mathrm{err}(d) $ of \textit{DeepMartingale} collected for $ d=2 $ to $10$ to infer an implied expressivity order from low-dimensional results, and then use the resulting order to scale up the neural networks, training designs and rebalancing frequency. The estimated orders are as follows:
$$
    \begin{cases}
        \mathrm{err}(d) \propto d^{1.15} , \; R^2 = 0.8361 \;  & \text{if} \; \; s_0=90 \\
         \mathrm{err}(d) \propto d^{1.24} , \; R^2 = 0.9012  \; & \text{if} \; \; s_0=100 \\
         \mathrm{err}(d) \propto d^{1.18} ,\; R^2 = 0.9202  \; & \text{if} \; \; s_0=110 \\
    \end{cases}
$$
We observe that when $ s_0 =90 $, $ R^2 \le 0.9 $. This is mainly because the asymmetrical payoff/dynamics make the learning problem less regular, and the fixed baseline configuration can already be near-optimal in very low dimensions (e.g., $ d=2,3 $). This flattens the initial error growth and makes the implied scaling law inferred from $d\le 10$ less stable. As a result, the fitted power law may underestimate the constant terms and may not fully reflect the scaling behavior at larger $d$.

\subsubsection{Comparison after scaling}\label{subsubsec:experiment_compare_after_scaling}

When scaling to $d=100$, the regressed power can lead to memory pressure for both \textit{DeepMartingale} and \textit{Guo-DeepPrimalDual}. For this reason, we reduce the scaling exponent to $ 1.07 $, which keeps \textit{DeepMartingale} memory-feasible for all $ d\le 100 $. In contrast, \textit{Guo-DeepPrimalDual} still runs Out-of-Memory at $d=100$.

\begin{figure}[tbp] 
\begin{minipage}[t]{0.33 \linewidth} 
\centering
\includegraphics[height=3.3cm,width=5.1cm]{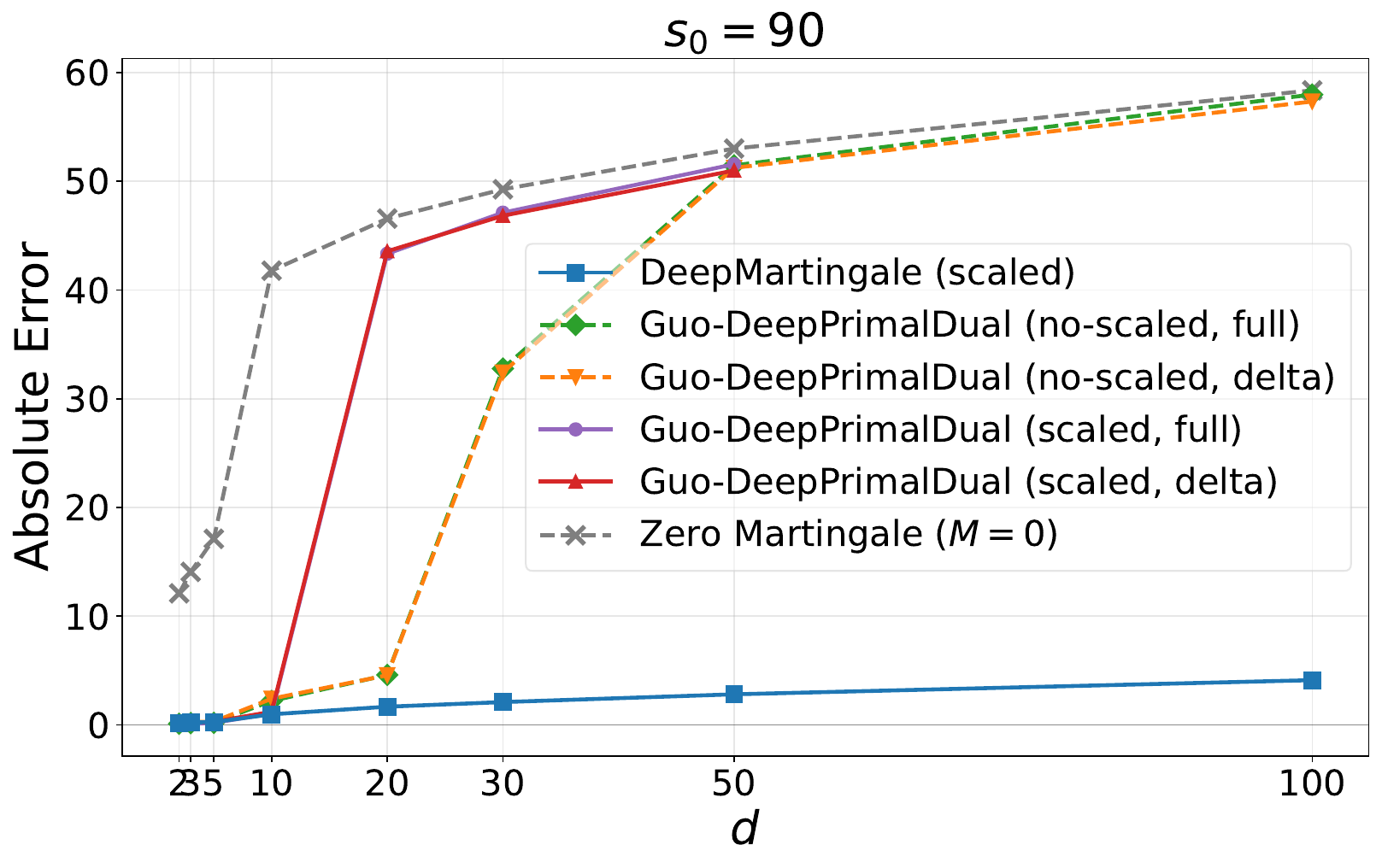}
\end{minipage}%
\begin{minipage}[t]{0.33 \linewidth} 
\centering
\includegraphics[height=3.3cm,width=5.1cm]{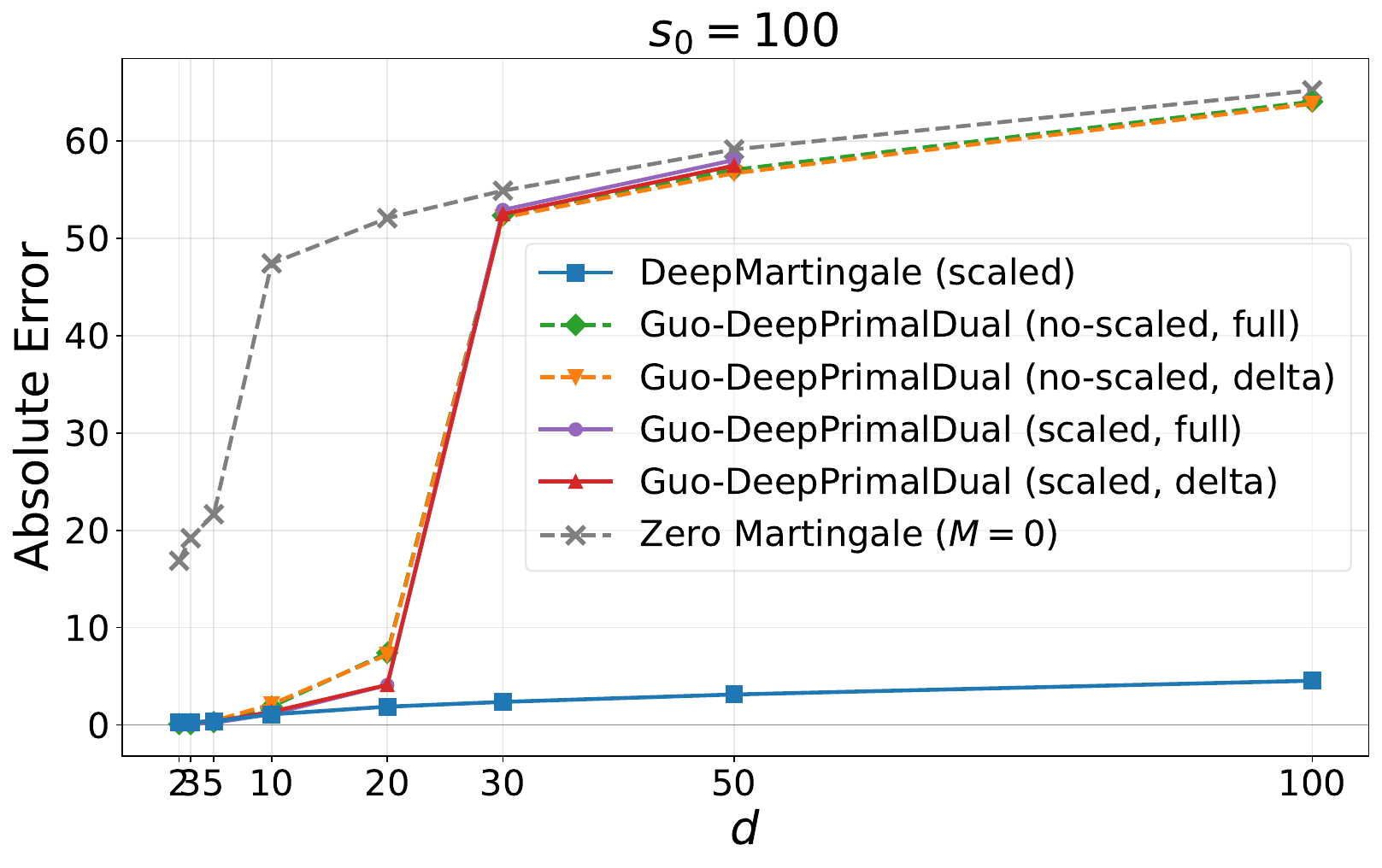}
\end{minipage}%
\begin{minipage}[t]{0.33 \linewidth} 
\centering
\includegraphics[height=3.3cm,width=5.1cm]{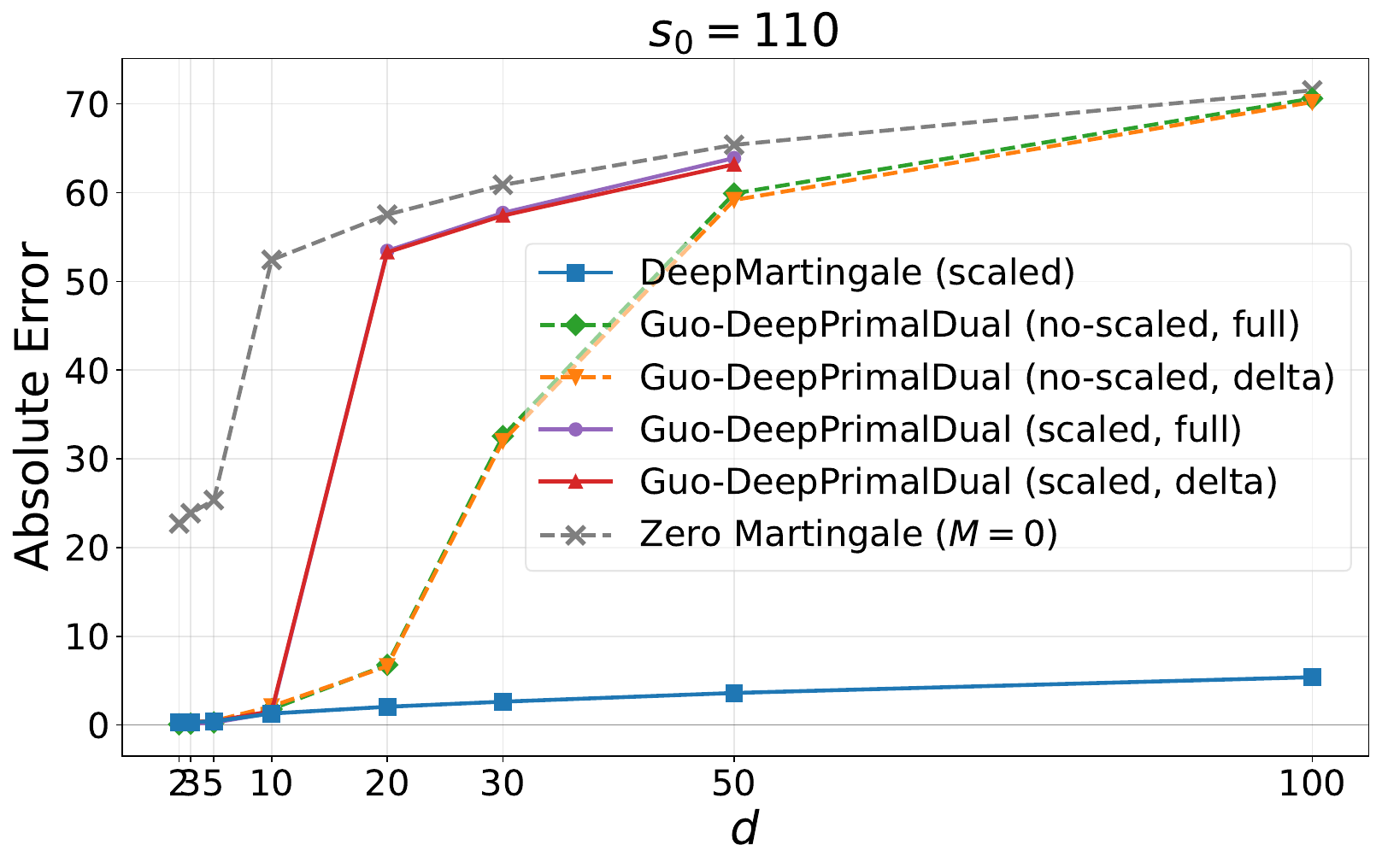}
\end{minipage}%
\caption{ $ \mathrm{err}(d) $ vs $d$ (Asymmetrical)}
\label{fig:express_asym_compare}
\end{figure}

\begin{figure}[tbp] 
\begin{minipage}[t]{0.33 \linewidth} 
\centering
\includegraphics[height=2.9cm,width=5cm]{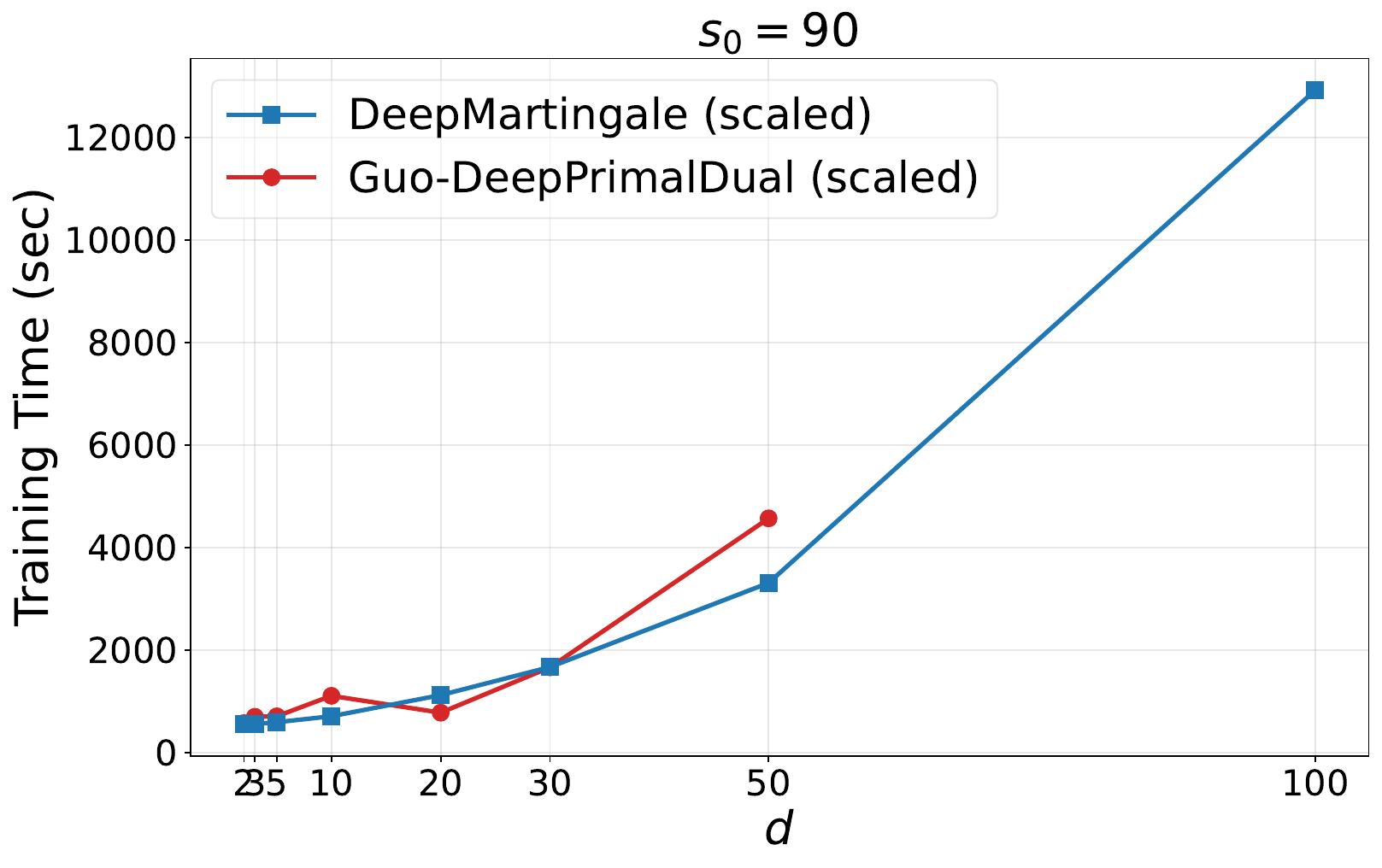}
\end{minipage}%
\begin{minipage}[t]{0.33 \linewidth} 
\centering
\includegraphics[height=2.9cm,width=5cm]{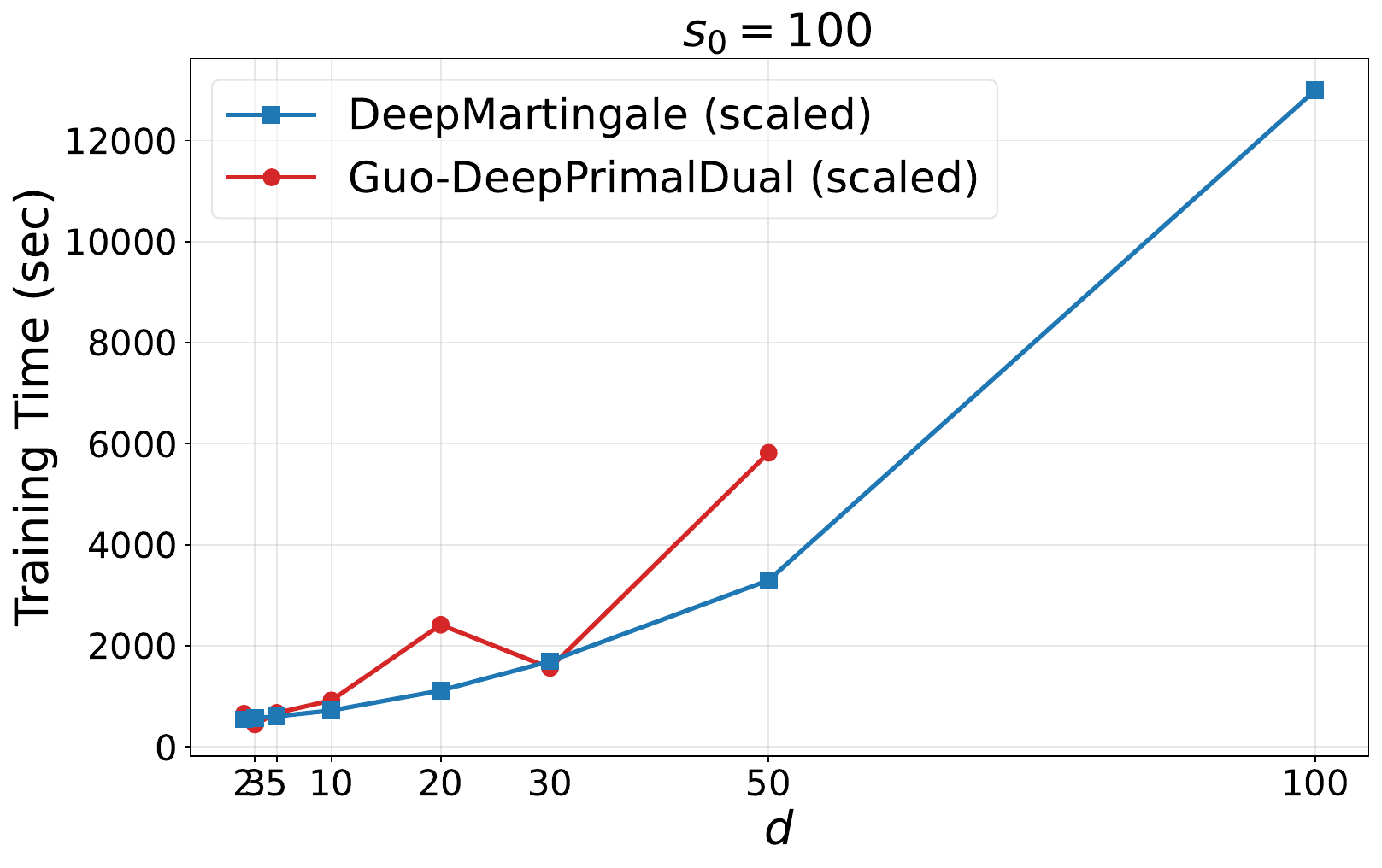}
\end{minipage}%
\begin{minipage}[t]{0.33 \linewidth} 
\centering
\includegraphics[height=2.9cm,width=5cm]{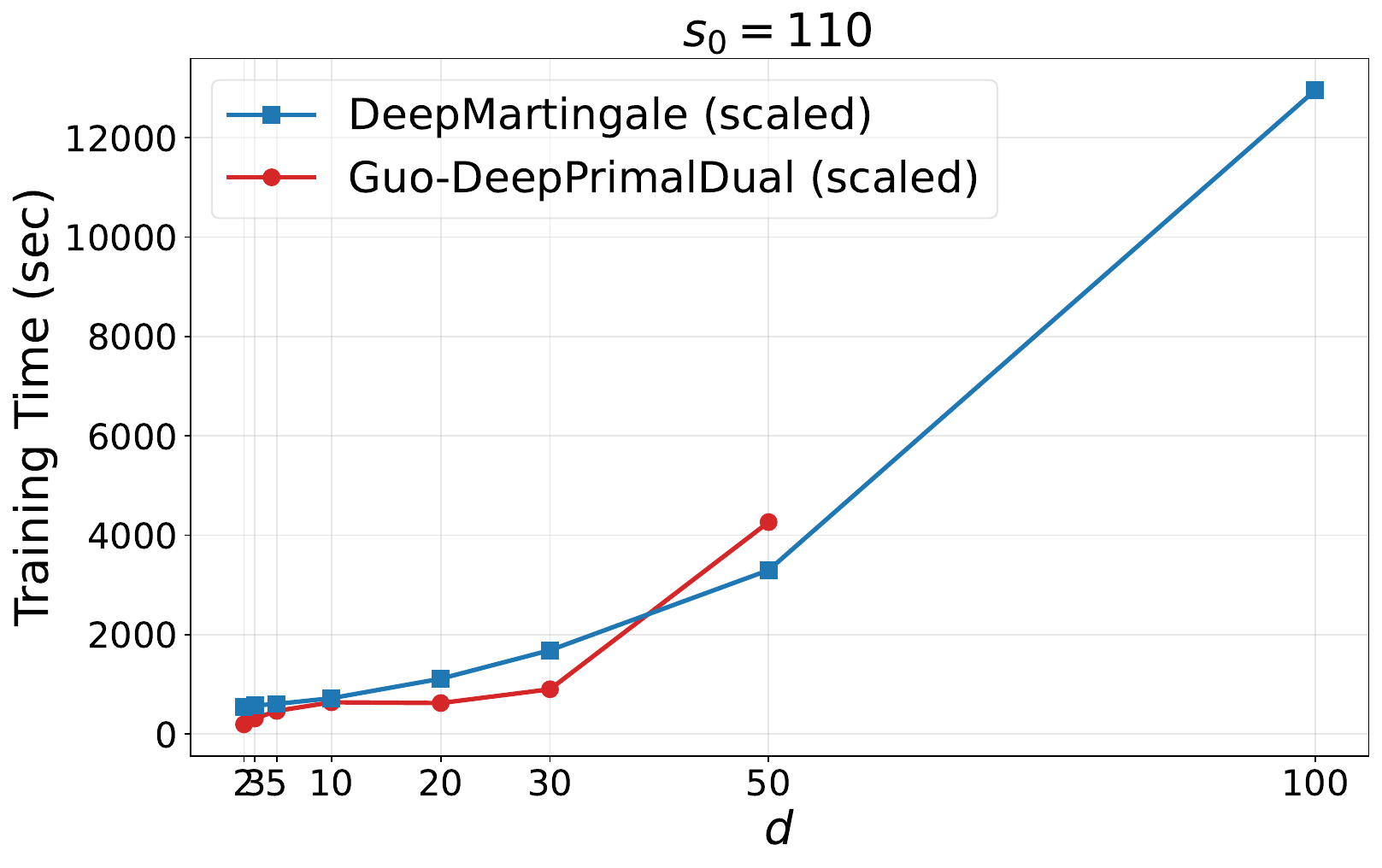}
\end{minipage}%
\caption{ Training times of \textit{DeepMartingale} vs \textit{Guo-DeepPrimalDual} (Asymmetrical)}
\label{fig:training_time_asym_DM}
\end{figure}

We compute the absolute error using the benchmark values in \cite{Becker19}, i.e., we set $ U^{d}_{\mathrm{ref}} $ to the point estimates in \cite[deep optimal stopping]{Becker19}. Figure~\ref{fig:express_asym_compare} reports the absolute-error curves versus the dimension. Here, ``full'' denotes the full martingale that includes the quadratic term, while ``delta'' denotes the delta martingale that only contains delta hedging. We observe that \textit{DeepMartingale} is substantially more accurate than \textit{Guo-DeepPrimalDual} in high dimensions. \textit{Guo-DeepPrimalDual} performs exceptionally well in low to mid dimensions (e.g., $2\le d\le 10$ or $20$). However, even when we scale \textit{Guo-DeepPrimalDual} as carefully as possible using the same orders, \textit{DeepMartingale} achieves the best performance among all curves as the dimension increases. In particular, \textit{Guo-DeepPrimalDual} tends to a degenerate martingale ($M\equiv 0$), suggesting that it becomes difficult to learn the high-dimensional structure and that the optimization may oscillate around the initialization. Together with the training-time comparison in Figure~\ref{fig:training_time_asym_DM}, this indicates that \textit{DeepMartingale} maintains good efficiency and accuracy as $d$ increases. Moreover, under the same scaling rule, \textit{DeepMartingale} remains memory-feasible at $d=100$, whereas \textit{Guo-DeepPrimalDual} fails, suggesting more efficient memory utilization for the pure-dual design.

Moreover, we compare against \textit{Alfonsi-PureDual} under scaling according to the same orders. The notations and results are reported in Table~\ref{Notation} and Table~\ref{table:max_call_asym}, respectively. We use \textbf{boldface} to highlight the best values among all methods, and denote by \textbf{N/A} the cases that run out of memory or exceed $24$ hours of runtime. Out-of-memory issues occur more frequently for the dual problem than for the primal problem, due to the nested substeps (observation points, rebalancing frequency) $K$ required for numerical integration, which in turn leads to $K$-times larger data sizes than in a pure primal approach. Therefore, we do not report results for $d>100$. Nevertheless, one can expect to scale to higher dimensions by using devices with larger memory or by adjusting hyperparameters such as network widths and observation frequency, or by employing more advanced implementation techniques.

Overall, these results further support the observation that \textit{Guo-DeepPrimalDual} achieves excellent performance in low to mid dimensions. This also suggests a possible direction for improving \textit{DeepMartingale} by incorporating a quadratic martingale term; however, we focus here on the high-dimensional setting and leave this extension for future research. As the dimension increases, \textit{DeepMartingale} dominates the performance among all considered approaches.

\begin{table}[ht]
\centering
\footnotesize
\begin{tabular}{llc}
		\toprule 
	 	Notation & \multicolumn{2}{c}{Description} \\
        \midrule 
        $U^{\text{dm}} $ & \multicolumn{2}{c}{Upper bound by \textit{DeepMartingale}.} \\
        $U^{\text{gd}} $ & \multicolumn{2}{c}{Upper bound by \textit{Guo-DeepPrimalDual}.} \\
        $U^{\text{af}} $ & \multicolumn{2}{c}{Upper bound by \textit{Alfonsi-PureDual}.} \\
        $U^{\text{zero}} $ & \multicolumn{2}{c}{Upper bound by zero martingale $M\equiv0$ (no hedging).} \\
      \midrule
        \multirow{1}{*}{$\text{Ref}_{2}$}  & $d=5$ & $95 \%$-CI from \citep[improved regression]{broadiecao08}. \\
        & $d\neq 5 $ & $95 \%$-CI from \citep[deep optimal stopping]{Becker19}. \\
		\bottomrule
	\end{tabular}
\caption{Notation List}
\label{Notation}
\end{table}

\begin{table}[ht]
\centering
\footnotesize
\begin{tabular}{cccccccccc}
	\toprule
		\multirow{2}{*}{$d$}   & \multirow{2}{*}{$s_0$}  & \multirow{2}{*}{$U^{\mathrm{dm}}$} & \multicolumn{2}{c}{$U^{\mathrm{gd}} $} & \multicolumn{2}{c}{$U^{\mathrm{af}} $} & \multirow{2}{*}{$U^{\mathrm{zero}} $} & \multirow{2}{*}{$\text{Ref}_{1}$} \\
		\cmidrule(lr){4-5} \cmidrule(lr){6-7}
		 &  &  & full & delta & vanilla & delta &  &  \\
		\midrule
		 & $90$  &  $14.505$ & $\mathbf{14.397}$ & $14.502$ & $ 14.561 $ & $ 15.577 $ &  $26.444 $ & $[14.299,14.367] $  \\
	    $2$    & $100$ &   $20.055$ & $\mathbf{19.876}$ & $20.022$ & $ 20.041 $ & $ 20.948 $ &  $36.676 $ & $[19.772,19.829] $ \\
		     &    $110$    &  $27.446$ & $\mathbf{27.277}$ & $27.452$ & $ 27.416 $ & $ 28.153 $ & $49.871 $ & $[27.138,27.163] $ \\
         \midrule
        & $90$  &  $19.319 $ & $\mathbf{19.177}$ & $19.276$ & $ 19.814 $ & $ 20.639 $ &  $33.165 $ & $[19.065 , 19.104] $  \\
	    $3$    & $100$ &   $26.957 $ & $\mathbf{26.821}$ & $26.954$ & $ 27.557 $ & $ 28.322 $ &  $45.879 $ & $[26.648 , 26.701] $ \\
		     &    $110$    &  $36.160 $ & $\mathbf{35.964}$ & $36.133$ & $ 36.883 $ & $ 37.569 $ &  $59.696 $ & $[35.806,35.835] $ \\
         \midrule
		 & $90$  &  $27.930 $ & $\mathbf{27.863}$ & $27.960$ & N/A & N/A &  $44.775 $  & $[27.468,27.686] $  \\
	     $5$   & $100$ &  $38.379 $ & $\mathbf{38.205}$ & $38.339$ & N/A & N/A &  $59.630 $ & $[37.730,38.020] $ \\
		     &    $110$    &  $49.911 $ & $\mathbf{49.765}$ & $49.910$ & N/A & N/A &  $74.852 $ & $[49.155 , 49.531] $ \\
         \midrule
		   & $90$ &  $\mathbf{86.970}$ & $87.037$ & $87.207$ & N/A & N/A &  $127.757 $ & $[85.610, 86.090] $ \\
	     $10$   & $100$ &  $\mathbf{105.839}$ & $105.884$ & $106.098$ & N/A & N/A &  $152.149 $  & $[104.830, 105.370] $ \\
		     &    $110$  &  $\mathbf{125.043}$ & $125.114$ & $125.314$ & N/A & N/A &  $176.167 $ & $[124.050, 124.660] $ \\
         \midrule
		   & $90$ &  $ \mathbf{127.776} $ & $169.427$ & $169.648$ & N/A & N/A &  $172.657 $ & $[125.819,126.383] $ \\
	     $20$   & $100$ &  $\mathbf{151.663}$ & $153.881$ & $153.900$ & N/A & N/A &  $201.820 $  & $[149.480,150.053] $ \\
		     &    $110$  &  $\mathbf{175.595}$ & $226.981$ & $226.811$ & N/A & N/A &  $231.061 $ & $[173.144,173.937] $ \\
         \midrule
		   & $90$ &  $\mathbf{156.806}$ & $201.803$ & $201.520$ & N/A & N/A &  $203.952 $ & $[153.800, 154.780] $ \\
	     $30$   & $100$ &  $\mathbf{183.960}$ & $234.481$ & $234.054$ & N/A & N/A &  $236.456 $  & $[180.640, 181.700] $ \\
		     &    $110$  &  $\mathbf{211.188}$ & $266.276$ & $265.967$ & N/A & N/A &  $269.430 $ & $[207.480, 208.620] $ \\
         \midrule
		   & $90$ &  $\mathbf{199.142} $ & $247.877$ & $247.300$ & N/A & N/A &  $249.318 $ & $[195.793, 196.963] $ \\
	     $50$   & $100$ &  $\mathbf{231.022}$ & $285.907$ & $285.324$ & N/A & N/A &  $286.993 $  & $[227.247, 228.605] $ \\
		     &    $110$  &  $\mathbf{262.936} $ & $323.218$ & $322.503$ & N/A & N/A &  $324.709 $ & $[258.661, 260.092] $  \\
         \midrule
		 & $90$   & $\mathbf{267.801}$  & N/A  & N/A & N/A & N/A & $322.028$ &   $[263.043, 264.425] $  \\
	     $100$   & $100$  & $\mathbf{307.315} $ & N/A & N/A & N/A & N/A & $367.939$ &   $[301.924, 303.843] $ \\
		     &    $110$   & $\mathbf{346.965}$ & N/A & N/A & N/A & N/A & $413.105$ &   $[340.580, 342.781] $   \\
		\bottomrule
	\end{tabular}
\caption{Bermudan max-call (Asymmetrical) after scaling }
\label{table:max_call_asym} 
\end{table}

\subsubsection{Hedging performance}\label{subsubsec:experiment_hedge}

\begin{figure}[tbp] 
\begin{minipage}[t]{0.33 \linewidth} 
\centering
\includegraphics[height=3.2cm,width=5.2cm]{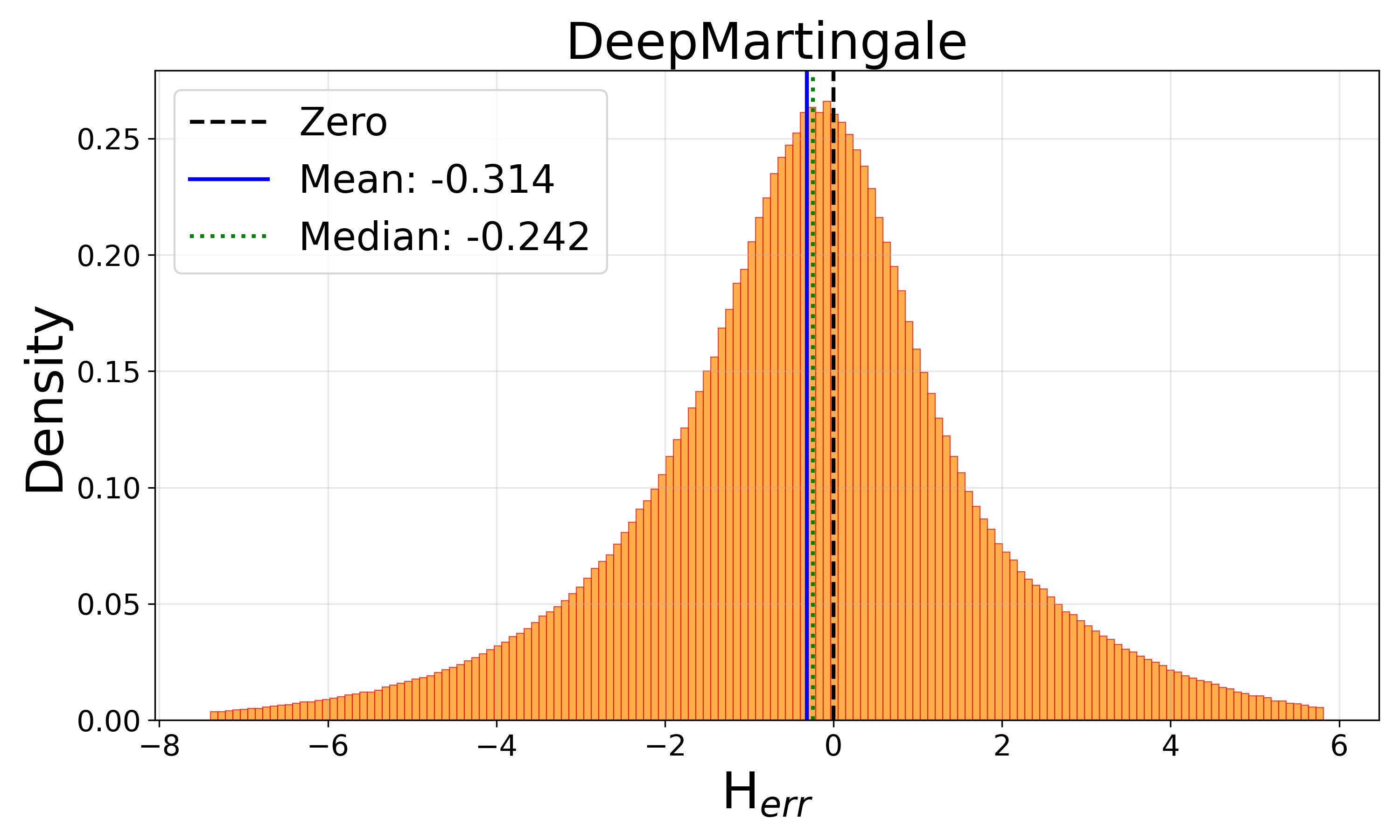}
\end{minipage}%
\begin{minipage}[t]{0.33 \linewidth} 
\centering
\includegraphics[height=3.2cm,width=5.2cm]{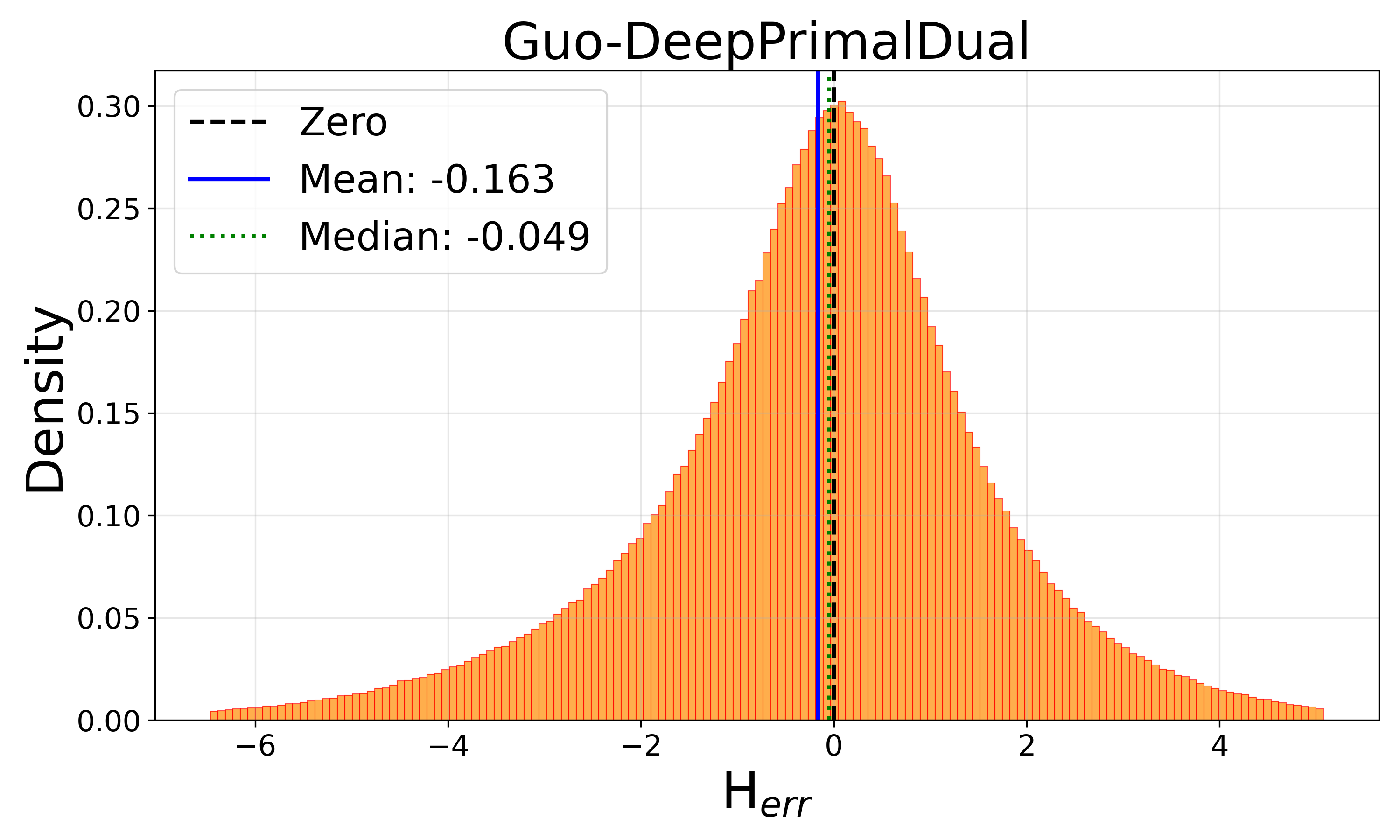}
\end{minipage}%
\begin{minipage}[t]{0.33 \linewidth} 
\centering
\includegraphics[height=3.2cm,width=5.2cm]{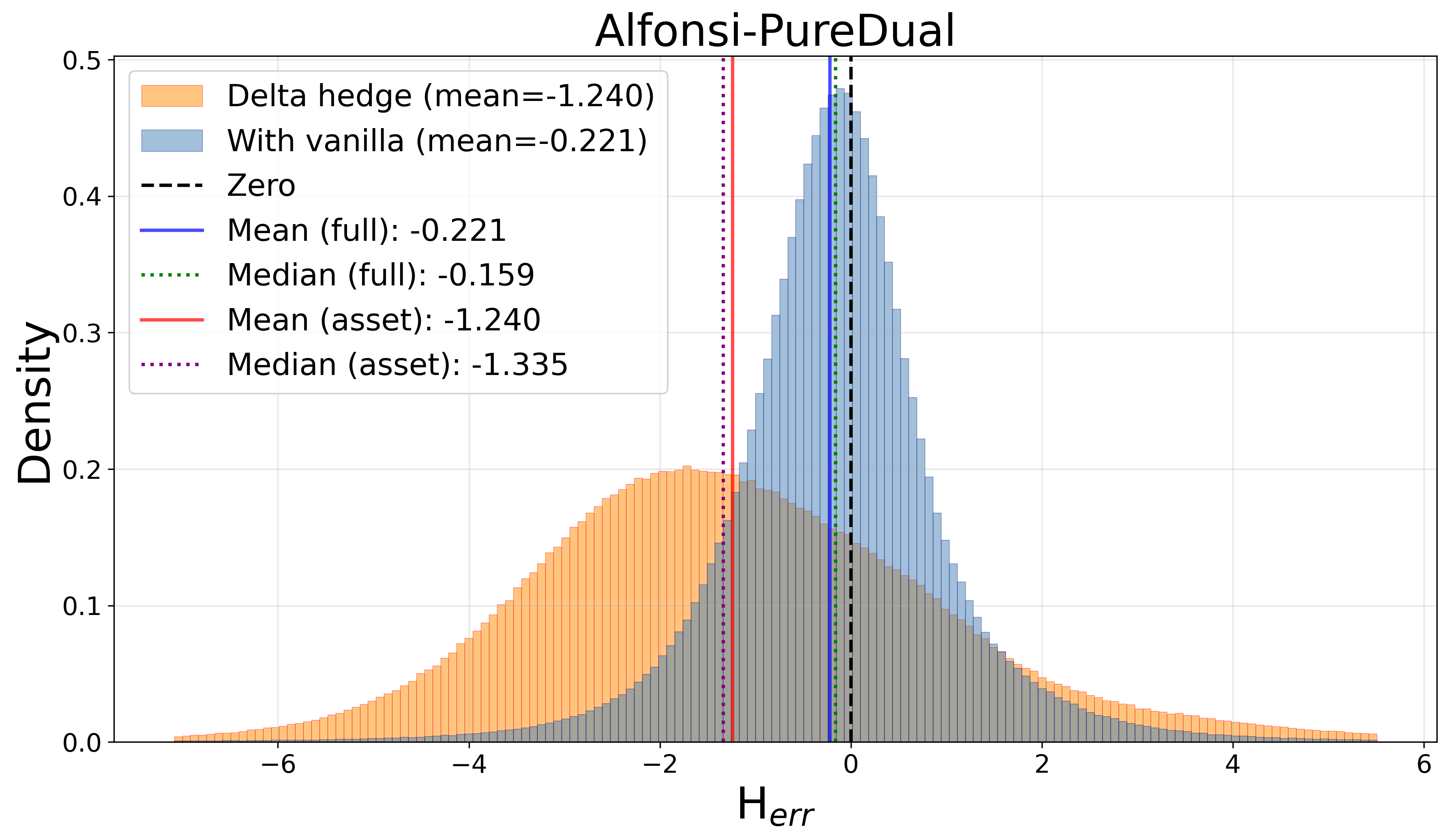}
\end{minipage}%
\caption{ Worst Case Hedging Error (Asymmetrical, $s_0=90, \; d=2$)}
\label{fig:hedging_compare_scaled_asym_d_2}
\end{figure}

In this part, we compare the hedging strategy at $s_0=90$ between the low-dimensional setting $d=2$ and the high-dimensional setting $d=50$. For $d=2$, we compare the ``deep'' delta hedging strategy produced by \textit{DeepMartingale} against the delta hedging strategy from \textit{Guo-DeepPrimalDual}, as well as hedging strategies with and without vanilla options written on the underlying assets from \textit{Alfonsi-PureDual}. As shown in Figure~\ref{fig:hedging_compare_scaled_asym_d_2}, \textit{Guo-DeepPrimalDual}'s delta hedge performs best among all delta-hedging portfolios. \textit{DeepMartingale} also achieves strong performance, although its distribution is slightly less peaked than that of \textit{Guo-DeepPrimalDual}. After adding at-the-money vanilla options, \textit{Alfonsi-PureDual} yields the best performance, exhibiting a more peaked shape and a smaller downside (i.e., less negative loss). This suggests a potential way to further improve hedging performance by incorporating vanilla options to mitigate second-order risk, which we leave for future research.

\begin{figure}[tbp] 
\begin{minipage}[t]{0.48 \linewidth} 
\centering
\includegraphics[height=4cm,width=7.2cm]{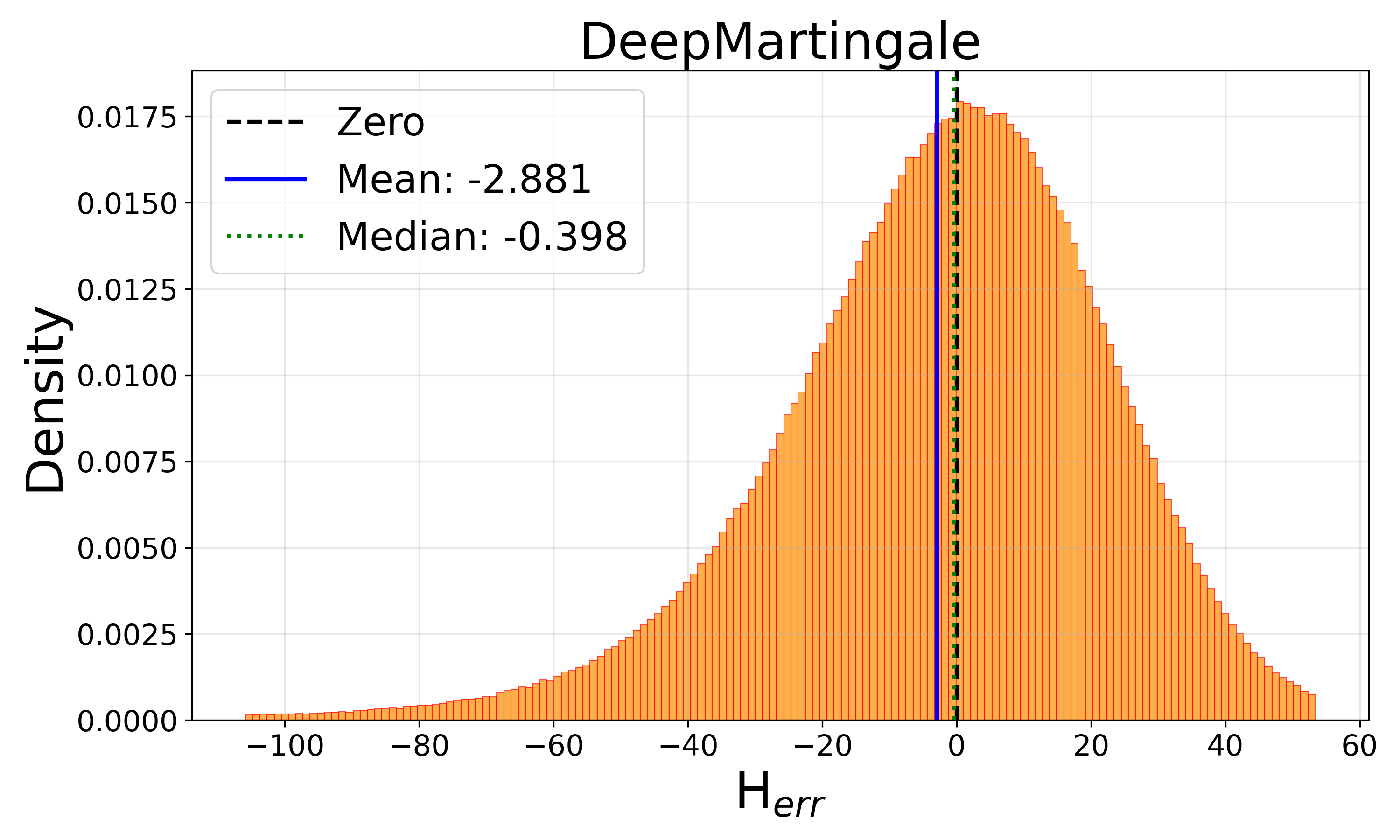}
\end{minipage}%
\begin{minipage}[t]{0.48 \linewidth} 
\centering
\includegraphics[height=4cm,width=7.2cm]{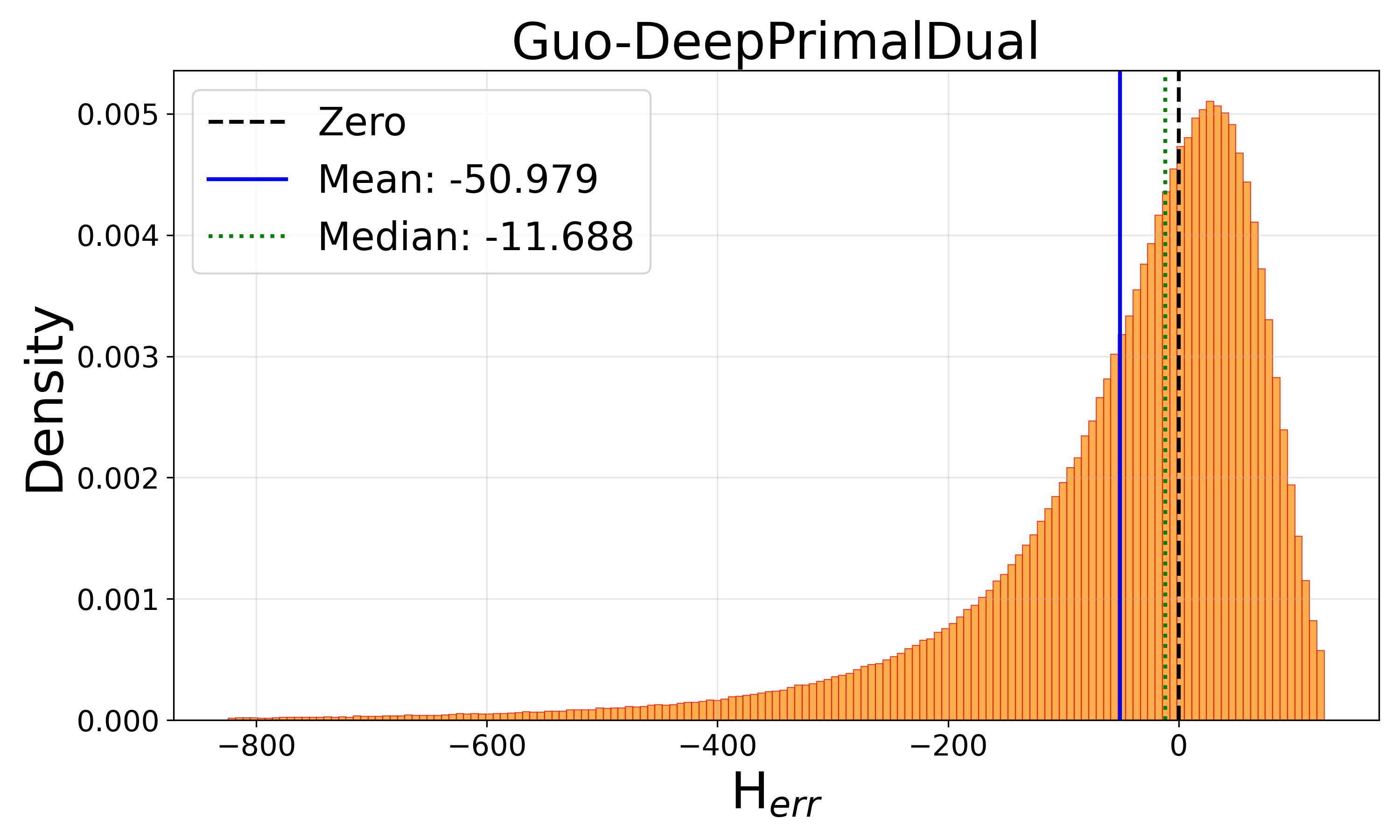}
\end{minipage}%
\caption{ Worst Case Hedging Error (Asymmetrical, $s_0=90, \; d=50$)}
\label{fig:hedging_compare_scaled_asym_d_50}
\end{figure}

For the high-dimensional case $d=50$, we only compare against the delta hedging strategy from \textit{Guo-DeepPrimalDual}, since the algorithm of \textit{Alfonsi-PureDual} does not easily scale to high dimensions. From Figure~\ref{fig:hedging_compare_scaled_asym_d_50}, it is apparent that \textit{DeepMartingale} delivers excellent performance for the delta-hedging portfolio, while \textit{Guo-DeepPrimalDual} fails to provide a reliable delta hedging strategy. This highlights the dimensional scalability of \textit{DeepMartingale} for high-dimensional hedging.

\bibliographystyle{abbrv}
\bibliography{main}   

\begin{thebibliography}{10}

\bibitem{puredual-mf}
A.~Alfonsi, A.~Kebaier, and J.~Lelong.
\newblock A pure dual approach for hedging bermudan options.
\newblock {\em Mathematical Finance}, 35(4):745--759, 2025.

\bibitem{andersen04}
L.~Andersen and M.~Broadie.
\newblock Primal-dual simulation algorithm for pricing multidimensional american options.
\newblock {\em Management Sci.}, 50(9):1222--1234, 2004.

\bibitem{BARTOLUCCI2023194}
F.~Bartolucci, E.~{De Vito}, L.~Rosasco, and S.~Vigogna.
\newblock Understanding neural networks with reproducing kernel banach spaces.
\newblock {\em Applied and Computational Harmonic Analysis}, 62:194--236, 2023.

\bibitem{RKBS2024}
F.~Bartolucci, E.~D. Vito, L.~Rosasco, and S.~Vigogna.
\newblock Neural reproducing kernel banach spaces and representer theorems for deep networks.
\newblock {\em Preprint}, 2024.
\newblock submitted March 13, https://arxiv.org/abs/2403.08750.

\bibitem{Becker19}
S.~Becker, P.~Cheridito, and A.~Jentzen.
\newblock Deep optimal stopping.
\newblock {\em J. Mach. Learn. Res.}, 20(74):1--25, 2019.

\bibitem{Becker20}
S.~Becker, P.~Cheridito, and A.~Jentzen.
\newblock Pricing and hedging american-style options with deep learning.
\newblock {\em J. Risk Financ. Manag.}, 13(7), 2020.

\bibitem{belome09}
D.~Belomestny, C.~Bender, and J.~Schoenmakers.
\newblock True upper bounds for {B}ermudan products via non-nested {M}onte {C}arlo.
\newblock {\em Math. Finance}, 19(1):53--71, 2009.

\bibitem{bender08}
C.~Bender, A.~Kolodko, and J.~Schoenmakers.
\newblock Enhanced policy iteration for {A}merican options via scenario selection.
\newblock {\em Quant. Finance}, 8(2):135--146, 2008.

\bibitem{broadiecao08}
M.~Broadie and M.~Cao.
\newblock Improved lower and upper bound algorithms for pricing {A}merican options by simulation.
\newblock {\em Quant. Finance}, 8(8):845--861, 2008.

\bibitem{brown10}
D.~B. Brown, J.~E. Smith, and P.~Sun.
\newblock Information relaxations and duality in stochastic dynamic programs.
\newblock {\em Oper. Res.}, 58(4-part-1):785--801, 2010.

\bibitem{carriere96}
J.~F. Carriere.
\newblock Valuation of the early-exercise price for options using simulations and nonparametric regression.
\newblock {\em Insur. Math. Econ.}, 19(1):19--30, 1996.

\bibitem{CSW2019}
J.~Chen, T.~Sit, and H.~Y. Wong.
\newblock Simulation-based value-at-risk for nonlinear portfolios.
\newblock {\em Quant. Finance}, 19(10):1639--1658, 2019.

\bibitem{Clement2002}
E.~Clement, D.~Lamberton, and P.~Protter.
\newblock An analysis of a least squares regression method for {A}merican option pricing.
\newblock {\em Finance \& Stochastics}, 6:449--471, 2002.

\bibitem{parto14}
G.~Da~Prato and J.~Zabczyk.
\newblock {\em Stochastic Equations in Infinite Dimensions}.
\newblock Cambridge University Press, 2014.

\bibitem{gonon23}
L.~Gonon.
\newblock Deep neural network expressivity for optimal stopping problems.
\newblock {\em Finance Stoch.}, 28:865–910, 2024.

\bibitem{GononRandomNN23}
L.~Gonon, L.~Grigoryeva, and J.-P. Ortega.
\newblock {Approximation bounds for random neural networks and reservoir systems}.
\newblock {\em Ann. Appl. Probab.}, 33(1):28 -- 69, 2023.

\bibitem{grohs21}
P.~Grohs and L.~Herrmann.
\newblock Deep neural network approximation for high-dimensional parabolic {H}amilton-{J}acobi-{B}ellman equations.
\newblock {\em Preprint}, 2021.
\newblock submitted March 9, https://arxiv.org/abs/2103.05744.

\bibitem{Jentzen23}
P.~Grohs, F.~Hornung, A.~Jentzen, and et~al.
\newblock A proof that artificial neural networks overcome the curse of dimensionality in the numerical approximation of black–scholes partial differential equations.
\newblock {\em Mem. Am. Math. Soc.}, 284(1410), 2023.

\bibitem{guo2024simultaneousupperlowerbounds}
I.~Guo, N.~Langrené, and J.~Wu.
\newblock Simultaneous upper and lower bounds of american-style option prices with hedging via neural networks.
\newblock {\em Quantitative Finance}, 25(4):509--525, 2025.

\bibitem{han18}
J.~Han, A.~Jentzen, and W.~E.
\newblock Solving high-dimensional partial differential equations using deep learning.
\newblock {\em PNAS}, 115(34):8505--8510, 2018.

\bibitem{Haugh04}
M.~B. Haugh and L.~Kogan.
\newblock Pricing {A}merican options: A duality approach.
\newblock {\em Oper. Res.}, 52(2):258--270, 2004.

\bibitem{DeepPrimalRandom}
C.~Herrera, F.~Krach, P.~Ruyssen, and J.~Teichmann.
\newblock Optimal stopping via randomized neural networks.
\newblock {\em Front. Math. Finance}, 3(1):31--77, 2024.

\bibitem{hornik91}
K.~Hornik.
\newblock Approximation capabilities of multilayer feedforward networks.
\newblock {\em Neural Networks}, 4(2):251--257, 1991.

\bibitem{Jentzen20}
M.~Hutzenthaler, A.~Jentzen, T.~Kruse, and T.~A. Nguyen.
\newblock A proof that rectified deep neural networks overcome the curse of dimensionality in the numerical approximation of semilinear heat equations.
\newblock {\em SN Partial Differ. Equations Appl.}, 1(10), 2020.

\bibitem{modern-Probability}
O.~Kallenberg.
\newblock {\em Foundations of Modern Probability}.
\newblock Springer Cham, 2021.

\bibitem{KohlerNN2010}
M.~Kohler, A.~Krzyżak, and N.~Todorovic.
\newblock Pricing of high-dimensional american options by neural networks.
\newblock {\em Mathematical Finance}, 20(3):383--410, 2010.

\bibitem{kolodko04}
A.~Kolodko and J.~Schoenmakers.
\newblock Upper bounds for {B}ermudan style derivatives.
\newblock {\em Monte Carlo Methods Appl.}, 10(3-4):331--343, 2004.

\bibitem{lapeyre2020neuralnetworkregressionbermudan}
B.~Lapeyre and J.~Lelong.
\newblock Neural network regression for bermudan option pricing, 2020.

\bibitem{longstaff01}
F.~A. Longstaff and E.~S. Schwartz.
\newblock Valuing {A}merican options by simulation: A simple least-squares approach.
\newblock {\em Rev. Financ. Stud.}, 14(1):113--147, 06 2001.

\bibitem{Ma-Zhang02}
J.~Ma and J.~Zhang.
\newblock Representation theorems for backward stochastic differential equations.
\newblock {\em Ann. Appl. Probab.}, 12(4):1390--1418, 2002.

\bibitem{MAO201191}
X.~Mao.
\newblock Linear stochastic differential equations.
\newblock In X.~Mao, editor, {\em Stochastic Differential Equations and Applications}, pages 91--106. Woodhead Publishing, second edition, 2011.

\bibitem{SDE03}
B.~{\O}ksendal.
\newblock {\em Stochastic Differential Equations}.
\newblock Springer Berlin, Heidelberg, 2003.

\bibitem{opschoor20}
J.~A.~A. Opschoor, P.~C. Petersen, and C.~Schwab.
\newblock Deep {R}e{LU} networks and high-order finite element methods.
\newblock {\em Anal. Appl.}, 18(05):715--770, 2020.

\bibitem{RAISSI2019}
M.~Raissi, P.~Perdikaris, and G.~Karniadakis.
\newblock Physics-informed neural networks: A deep learning framework for solving forward and inverse problems involving nonlinear partial differential equations.
\newblock {\em J. Comput. Phys.}, 378:686--707, 2019.

\bibitem{reppen23}
A.~M. Reppen, H.~M. Soner, and V.~Tissot-Daguette.
\newblock Neural optimal stopping boundary.
\newblock {\em Math. Finance}, 35(2):441--469, 2025.

\bibitem{Roger02}
L.~C.~G. Rogers.
\newblock Monte carlo valuation of {A}merican options.
\newblock {\em Math. Finance}, 12(3):271--286, 2002.

\bibitem{roger10}
L.~C.~G. Rogers.
\newblock Dual valuation and hedging of {B}ermudan options.
\newblock {\em SIAM J. Financ. Math.}, 1(1):604--608, 2010.

\bibitem{schoen13}
J.~Schoenmakers, J.~Zhang, and J.~Huang.
\newblock Optimal dual martingales, their analysis, and application to new algorithms for bermudan products.
\newblock {\em SIAM J. Financ. Math.}, 4(1):86--116, 2013.

\bibitem{vanroy01}
J.~Tsitsiklis and B.~Van~Roy.
\newblock Regression methods for pricing complex american-style options.
\newblock {\em IEEE Trans. Neural Networks}, 12(4):694--703, 2001.

\bibitem{yang2024deepprimaldualbsdemethod}
J.~Yang and G.~Li.
\newblock A deep primal-dual bsde method for optimal stopping problems, 2024.

\bibitem{zhangjianfeng17}
J.~Zhang.
\newblock {\em Backward Stochastic Differential Equations}.
\newblock Springer New York, New York, 2017.

\end{thebibliography}


\appendix

\begin{appendix}

\section{Detailed Theories \& Proofs for Section \ref{sec:Numerical_Approximation}}


To investigate numerical approximation for the non-driver decoupled FBSDE~\eqref{eq:general_FBSDE}, we consider the time horizon $[0,T]$ for simplicity. To preserve full theoretical generality, we first relax the coefficient functions $(a^d,b^d)$ of the It\^o process to random coefficients $a^d(\omega,t,x),b^d(\omega,t,x)$ (see \citep{zhangjianfeng17}).

\begin{assumption}\label{ass:N_0_structural_ass_model}
      $a^d(\omega,t,x),b^d(\omega,t,x) $ are $\mathcal{F}\otimes \mathcal{B}(R^{1+d})$-measurable functions that satisfy the following conditions: 
    \begin{enumerate}
        \item for any $x\in \mathbb{R}^d $, mappings $(\omega,t) \mapsto a^d(\omega,t,x),\; (\omega,t) \mapsto b^d(\omega,t,x) $ are $\mathbb{F} $-progressively measurable;
        \item $a^d,b^d $ are uniformly Lipschitz continuous in $x$ and for almost all $ (t,\omega) $,
        \begin{equation}\label{eq:N_0_model_Lip}
        \Lip a^d(\omega,t,\cdot), \Lip b^d(\omega,t,\cdot) \le C (\log d)^{\frac{1}{2}} ; \; \text{and}
    \end{equation}
    \item  $a^{0;d}_t(\omega):=a^d(\omega,t,0),b^{0;d}_t(\omega):=b^d(\omega,t,0) $ with
    \begin{equation}\label{eq:N_0_model_bound}
    \| a^{0;d}_t \|_{L^2} ,\;  \| b^{0;d}_t \|_{L^2} \le  C d^Q 
    \end{equation}
    \end{enumerate}
    for some positive constants $C,Q$ independent of $d$. For notational simplicity, we omit $\omega $ in $a,b$. 
\end{assumption}

Similarly, let $\bar{g}^d$ be a random terminal function $\bar{g}(\omega,x)$, and impose the following condition.

\begin{assumption}\label{ass:N_0_structural_ass_g}
      $\bar{g}^d$ satisfies
    \begin{equation*}
    \Lip \bar{g}^d(\omega,\cdot) ,\;      \mathbb{E}|\bar{g}^d(\cdot,0)| \le C d^{Q} ,\; \mathbb{P}\text{-a.s.} \; \omega
    \end{equation*}
    for some positive constants $C,Q$ independent of $d$. The constants $C,Q$ are the same as in Assumption \ref{ass:N_0_structural_ass_model}, which is ensured by taking their maximum. For notational simplicity, we omit $\omega $ in $\bar{g}^d$.
\end{assumption}

Based on the Lipschitz condition in Assumption \ref{ass:N_0_structural_ass_model}, we immediately obtain a linear growth bound for the coefficients.

\begin{proposition}[Coefficient Linear Growth]\label{prop:dynamic_bound_new}
    Under Equation \eqref{eq:N_0_model_Lip} in Assumption \ref{ass:N_0_structural_ass_model}, we have for $ dt\times d\mathbb{P} \text{-a.e.} \; (t,\omega)  $ and all $ x\in \R^d $,
    \begin{equation*}
        \|a^d(t,x)\| \le C(\log d)^{\frac{1}{2}}\|x\| + \|a^d(t,0)\| ,\quad \|b^d(t,x)\|_{\H} \le C(\log d)^{\frac{1}{2}}\|x\| + \|b^d(t,0)\|_{\H}
    \end{equation*}
\end{proposition}

\begin{proof}[Proof of Proposition \ref{prop:dynamic_bound_new}]
The argument is direct. Indeed,
    \begin{equation*}
        \|a^d(t,x)\| \le \|a^d(t,x)-a^d(t,0)\| + \|a^d(t,0)\| \le C(\log d)^{\frac{1}{2}}\|x\| + \|a^d(t,0)\| .
    \end{equation*}
    The same reasoning applies to $b^d$, which completes the proof.
\end{proof}

By Assumption \ref{ass:N_0_structural_ass_g} and the same Lipschitz-based argument as in Proposition \ref{prop:dynamic_bound_new}, we also obtain the following linear growth estimate for $\bar{g}^d$.

\begin{proposition}[$\bar{g}^d$ Linear Growth]\label{prop:g_bound_new}
    Under Assumption \ref{ass:N_0_structural_ass_g}, we have
    \begin{equation*}
        |\bar{g}^d(x)| \le Cd^Q\|x\| + |\bar{g}^d(0)| ,  \; \mathbb{P}\text{-a.s.} \; \omega.
    \end{equation*}
\end{proposition}

\subsubsection{Proof of expressivity for SDEs and a specific type of BSDE}

Let $X^{*,s}_t $ denote the pathwise maximum $X^{*,s}_t := \sup_{s\le u \le t }\|X_u \| $ (or $\|X_u \|_{\H}$). For a $ d\times d \times d $ tensor $X=(X_1,\ldots,X_d)$ where $X_i(i=1,\ldots,d)$ are $d\times d$ matrices, we recall that $\|X\|_{\H}^2 = \sum_{i=1}^d \|X_i\|^2_{\H} $. We introduce the following notations, which will be used throughout the numerical approximation expressivity theory.
\begin{itemize}
    \item[(i)] $L^0(\mathcal{F}, \mathbb{R}^n) : \mathbb{R}^n\text{-valued } \mathcal{F}\text{-measurable random variable}$.
 \item[(ii)] $L^0(\mathbb{F}, \mathbb{R}^n) : \mathbb{R}^n\text{-valued } \mathbb{F}\text{-progressively measurable process}
$.
\item[(iii)] $L^{p,q}(\mathbb{F}, \mathbb{P}, \mathbb{R}^n) := \{ Z \in L^0(\mathbb{F}, \mathbb{R}^n) : ( \int_0^T \|Z_t\|^p \, dt )^{\frac{1}{p}} \in L^q(\mathcal{F}_T ; \mathbb{R}) \}
$, and denote the norm by $\| Z \|^q_{p,q}:= \mathbb{E}( \int_0^T \|Z_t\|^p \, dt )^{\frac{q}{p}} $ or Hilber-Schmit norm $ \| Z \|^q_{p,q;\H}:= \mathbb{E}( \int_0^T \|Z_t\|_{\H}^p \, dt )^{\frac{q}{p}} $.
\item[(iv)] $S^p(\mathbb{F}, \mathbb{P}, \mathbb{R}^n) := \{ Y \in L^0(\mathbb{F}, \mathbb{R}^n) : Y \text{ continuous (in } t\text{) } \mathbb{P}\text{-a.s. and } Y^{*,0}_T \in L^p(\mathcal{F}_T ; \R) \}
$.
\end{itemize}

To bound the SDE solution under Assumption \ref{ass:N_0_structural_ass_model} and the FBSDE solution under Assumption \ref{ass:N_0_structural_ass_g} with the desired expression rate, we require an expressivity-oriented version of the BDG inequality.

\begin{lemma}[BDG inequality (one-sided)] \label{lem:BDG_new}
    For any $p>0$, there exists a universal constant $C_p >0 $ that depends only on $p$, such that for any $\sigma^d \in L^{2,p}(\mathbb{F},\mathbb{P},\mathbb{R}^d)$, $M^d_t:=\int_0^t\sigma^d_s \cdot dW^d_s $, we have \( \big\| M^{*,0;d}_T \big\|^p_{L^p} \le C_p \| \sigma^d \|^p_{2,p}  \).
    
    If $p\ge 2 $, for any $\sigma^d \in L^{2,p}(\mathbb{F},\mathbb{P},\mathbb{R}^{d\times d}) $ (or $L^{2,p}(\mathbb{F},\mathbb{P},\mathbb{R}^{d\times d\times d}),\sigma^d = (\sigma^d_1,\ldots,\sigma^d_d)$), $M^d_t:=\int_0^t\sigma^d_s  dW^d_s $ is a $d$-dimensional vector (or $d\times d $ matrix) martingale, then \(\big\| M^{*,0;d}_T \big\|^p_{L^p} \le C_p \|\sigma^d_t\|^p_{2,p;\H} \).
\end{lemma}

\begin{remark}
    The proof follows the same steps as Theorem 2.4.1 in \citep{zhangjianfeng17}. In our statement, we emphasize that the constant $C_p$ is universal in the sense that it depends only on $p$ and is independent of the dimension $d$. Moreover, we extend the result to the multi-dimensional tensor case; for $p\ge 2 $, the inequality coincides with the classical formulation. While \citep{parto14} extend the original BDG inequality to an infinite-dimensional setting (general Hilbert space), their framework typically requires predictability. For our purposes, we instead work with $L^{2,p}$ to retain the required generality.
\end{remark}

\begin{proof}[Proof of Lemma~\ref{lem:BDG_new}]
As in \citep{zhangjianfeng17}, we assume $M^{*,0;d}_T, \int_0^t \|\sigma^d_s\|^2 ds  $ are bounded, which can be easily extended to the unbounded case using the truncation method e.g., see \citep{zhangjianfeng17}. For $p\ge 2 $, we use the argument in \citep{zhangjianfeng17}, which contains no dimension-dependent constant, and we obtain
    \[
        \begin{aligned}
           \big\| M^{*,0;d}_T \big\|^p_{L^p} = p\int_{0}^{\infty}\lambda^{p-1}\mathbb{P}(M^{*,0;d}_T \ge \lambda ) d\lambda  & \le p\int_{0}^{\infty}\lambda^{p-2}\mathbb{E}[|M^d_T|1_{(M^{*,0;d}_T\ge \lambda)} ] d\lambda \\
            & = \frac{p}{p-1}\mathbb{E}[|M^d_T| \cdot |M^{*,0;d}_T|^{p-1} ].
        \end{aligned}
    \]
    Then, by the H\"older inequality, \( \big\| M^{*,0;d}_T \big\|^p_{L^p} \le \frac{p}{p-1} \big\| M^d_T \big\|_{L^p} \big\| M^{*,0;d}_T \big\|^{p-1}_{L^p}  \).
    As $M^{*,0;d}_T $ is bounded, we have
    \begin{equation}\label{eq:BDG_p_ge_2_1}
        \big\| M^{*,0;d}_T \big\|^p_{L^p} \le (\frac{p}{p-1})^p \big\| M^d_T \big\|^p_{L^p}.
    \end{equation}
    By applying the It\^o formula under the multi-dimension condition ($\sigma $ and $W^d$ are part of the $\mathbb{R}^d$-valued process), we obtain 
    \begin{align*}
        d|M^d_t|^2 &= \|\sigma^d_t\|^2 dt + 2M^d_t \sigma^d_t \cdot dW^d_t ,\; \text{and} \\
        d|M^d_t|^p &= d\big(|M^d_t |^2 \big)^{\frac{p}{2}} = \frac{1}{2}p(p-1)|M^d_t|^{p-2} \|\sigma^d_t\|^2 dt + p|M^d_t|^{p-2}M^d_t \sigma^d_t \cdot dW^d_t,
    \end{align*}
    where $\mathbb{E}[\int_0^T p|M^d_t|^{p-2}M^d_t \sigma^d_t \cdot dW^d_t ]=0  $. Thus, by the H\"older inequality,
    \begin{equation}\label{eq:BDG_p_ge_2_2}
        \begin{aligned}
             \big\| M^d_T \big\|^p_{L^p} & \le C_p\mathbb{E}\Big[\int_0^T |M^d_t|^{p-2} \|\sigma^d_t\|^2 dt \Big] \le C_p\mathbb{E}\Big[|M^{*,0;d}_T| \int_0^T\|\sigma^d_t\|^2 dt \Big]  \le C_p  \big\| M^{*,0;d}_T \big\|^{p-2}_{L^p}  \|  \sigma^d \|^2_{2,p}  ,
        \end{aligned}
    \end{equation}
    where $C_p = \frac{1}{2}p(p-1)$. Then, combining Equations \eqref{eq:BDG_p_ge_2_1} and \eqref{eq:BDG_p_ge_2_2}, we obtain \( \big\| M^{*,0;d}_T \big\|^{p}_{L^p} \le C_p \|  \sigma^d \|^p_{2,p} \),
    where $C_p=(\frac{p}{p-1})^{\frac{p^2}{2}}(\frac{p(p-1)}{2} )^{\frac{p}{2}} $.

    For the $d\times d$-matrix and $d\times d\times d $-tensor, the argument is the same. The $d\times d$-matrix scenario is obvious if one replaces the norm with $\|\cdot\|_{\H}$ and the inner product with trace $\Tr(\cdot), $ which is a Hilbert inner product. Similarly, the $d\times d\times d $-tensor scenario is
    \begin{align*}
        d\|M^d_t\|_{\H}^2 &= \|\sigma^d_t\|_{\H}^2 dt + 2 \sum\limits_{i=1}^d (M^{i;d}_t)^{\Transpose} \sigma_t^{i;d} \cdot dW^d_t ,\; \text{and} \\
        d\|M^d_t\|_{\H}^p &= d\big(\|M^d_t \|_{\H}^2 \big)^{\frac{p}{2}} = \frac{p}{2}\|M^d_t\|_{\H}^{p-4} \left(\frac{p}{2}\|M^d_t\|_{\H}^{2} \|\sigma^d_t\|_{\H}^2 + \frac{p(p-2)}{2}\|\sum\limits_{i=1}^d (M^{i;d}_t)^{\Transpose} \sigma_t^{i;d}\|^2 \right) dt \notag \\
        & + p\|M^d_t\|_{\H}^{p-2}(M^d_t)^{\Transpose} \sigma^d_t \cdot dW^d_t,
    \end{align*}
where $M^{j;d} $ denotes the $j$-th column of $M^d$ and $\sigma^{j;d} =(\sigma^{j;d}_1,\ldots,\sigma^{j;d}_d) $; $\sigma^{j;d}_i $ is the $j $-th column of $\sigma^d_i$, and the only difference would be filled by 
    \begin{equation*}
        \begin{aligned}
            \Big \| \sum_{j=1}^d [ (M^{j;d}_t)^{\Transpose} \sigma_t^j ] \Big \|^2_{\H} & = \sum_{i=1}^d \Big[ \sum_{j=1}^d \sum_{k=1}^d (M^{j;d}_t)^{\Transpose}_{k} (\sigma^d_t)^j_{k,i}  \Big]^2=\sum_{i=1}^d \Big[\sum\limits_{j=1}^d \sum\limits_{k=1}^d (M^d_t)^{\Transpose}_{j,k} \left((\sigma^d_t)_i \right)_{k,j}  \Big]^2 \\
            & = \sum_{i=1}^d \big| \Tr( (M^d_t)^{\Transpose} (\sigma^d_t)_i ) \big|^2
        \end{aligned}
    \end{equation*}
    and \( \sum_{i=1}^d \big| \Tr( (M^d_t)^{\Transpose} (\sigma^d_t)_i ) \big|^2 \le \|M^d_t\|_{\H}^2 \sum_{i=1}^d\|(\sigma^d_t)_i\|^2_{\H} = \|M^d_t\|_{\H}^2 \|\sigma^d_t\|_{\H}^2 \)
    via the Cauchy-Schwartz inequality. Then the subsequent analysis is the same as in the previous argument.

    For $0< p< 2$, as in \citep{zhangjianfeng17}, we denote the quadratic variation $\langle M^d \rangle_t = \int_0^t \|\sigma^d_s\|^2 ds $ and
    \begin{equation*}
        \big\| \langle M^d \rangle_t^{\frac{p-2}{4}} \sigma^d_t \big\|_{L^2} =\mathbb{E}\Big[\int_0^T \langle M^d \rangle_t^{\frac{p-2}{2}} \|\sigma^d_t\|^2 dt  \Big] = \mathbb{E}\Big[\int_0^T \langle M^d \rangle_t^{\frac{p-2}{2}} d\langle M^d \rangle_t \Big] = \frac{2}{p}\mathbb{E}[\langle M^d \rangle_T^{\frac{p}{2}} ] <\infty.
    \end{equation*}
    Define $N^d_t:=\int_0^t \langle M^d \rangle_s^{\frac{p-2}{4}} \sigma^d_s\cdot dW^d_t $, which is square-integrable and $\mathbb{E}[(N^d_T)^2] = \frac{2}{p}\mathbb{E}[\langle M^d \rangle_T^{\frac{p}{2}} ] $. Then, by the It\^o formula,
    \begin{equation*}
        M^d_t = \int_0^t \langle M^d \rangle_s^{\frac{2-p}{4}} dN^d_s = \langle M^d \rangle_t^{\frac{2-p}{4}} N^d_t - \int_0^t N^d_s d\langle M^d \rangle_s^{\frac{2-p}{4}}.
    \end{equation*}
    Given the monotonicity of $\langle M^d \rangle $ in $t$, 
    \begin{equation*}
        M^{*,0;d}_T \le \langle M^d \rangle_T^{\frac{2-p}{4}} N_T^{*,0;d} + \int_0^T|N^d_s|d\langle M^d \rangle_s^{\frac{2-p}{4}} \le 2 N^{*,0;d}_T\langle M^d \rangle_T^{\frac{2-p}{4}}.
    \end{equation*}
    Given the H\"older inequality,
    \begin{equation*}
        \big\|M^{*,0;d}_T \big\|^p_{L^p} \le 2^p \mathbb{E}[|N^{*,0;d}_T|^p \langle M^d\rangle_T^{\frac{p(2-p)}{4}} ] \le 2^p   \big\| N^{*,0;d}_T \big\|_{L^2}^p \big\|    \left(\mathbb{E}[\langle M^d \rangle_T^{\frac{p}{2}}] \right)^{\frac{2-p}{2}}.
    \end{equation*}
    Thus, according to the Doob maximum inequality in Lemma 2.2.4 in \citep{zhangjianfeng17}, we have
    \begin{equation*}
         \big\|M^{*,0;d}_T \big\|^p_{L^p} \le 4^p  \big\| N^d_T \big\|_{L^2}^p \left(\mathbb{E}[\langle M^d \rangle_T^{\frac{p}{2}}] \right)^{\frac{2-p}{2}} = 4^p (\frac{2}{p})^{\frac{p}{2}} \left(\mathbb{E}[\langle M^d \rangle_T^{\frac{p}{2}}] \right)^{\frac{p}{2}} \left(\mathbb{E}[\langle M^d \rangle_T^{\frac{p}{2}}] \right)^{\frac{2-p}{2}} = C_p \mathbb{E}[\langle M^d \rangle_T^{\frac{p}{2}}],
    \end{equation*}
    where $C_p=4^p (\frac{2}{p})^{\frac{p}{2}} $, which completes the proof. 
\end{proof}

Following \citep{zhangjianfeng17}, to establish the expressivity proof for $K$, we first provide a bound for the solution of an SDE under Assumption \ref{ass:N_0_structural_ass_model} with the expression rate, where the regularity is satisfied as in \citep{zhangjianfeng17}.
\begin{theorem} \label{theorem:dynamic_bound_new}
    Given $p\ge 2$, under Assumption \ref{ass:N_0_structural_ass_model}, and the further assumption that $a^d(t,0) \in L^{1,p}(\mathbb{F},\mathbb{P},\mathbb{R}^d)$, $b^d(t,0) \in L^{2,p}(\mathbb{F},\mathbb{P},\mathbb{R}^{d\times d}) $ is bounded by $ C d^Q $, which is the same as in Assumption \ref{ass:N_0_structural_ass_model}, the following properties of the SDE's solution holds:
    \begin{equation}
        \big\| X^{*,0;d}_T \big\|_{L^p} \le B_p d^{Q_p} (\log d)^{R_p} (1+\|x^d_0\|^p),
    \end{equation}
    where $B_p,Q_p,R_p $ are constants independent of $d$.
\end{theorem}
\begin{proof}[Proof of Theorem~\ref{theorem:dynamic_bound_new}]
This proof also follows \citep{zhangjianfeng17}, but it is first necessary to clarify the expression rate. For $p\ge 2 $, without loss of generality, we assume $X^d \in L^p(\mathbb{F},\mathbb{P},\mathbb{R}^d ) $ (the general case can be solved by the truncation method given in \citep{zhangjianfeng17}). First, we derive
    \begin{equation*}
        X^{*,0;d}_T \le \|x^d_0\| + \int_{0}^T \|a^d(s,X^d_s)\|ds + \sup\limits_{0\le t\le T} \Big\|\int_0^t b^d(s,X^d_s) dW^d_s \Big\|
    \end{equation*}
    according to \citep{zhangjianfeng17}. By Lemma \ref{lem:BDG_new} (BDG inequality), Proposition \ref{prop:dynamic_bound_new}, and the Jensen inequality, we have
    \begin{align}
        \big\| X^{*,0;d}_T \big\|_{L^p}  &  \le 3^{p-1} \left(\|x^d_0\|^p + \big\| a^d(\cdot,X^d_{\cdot}) \big\|^p_{1,p} + \big\| b^d(\cdot,X^d_{\cdot}) \big\|^p_{2,p;\H} \right) \notag \\
        & \le B_p (\log d )^{\frac{p}{2}} \Big(1+\|x^d_0\|^p + \E \int_0^T \|X^d_t\|^p dt  \Big) \label{eq:max_p_X}
    \end{align}
    for $B_p=3^{p-1}(2+T)^p C^p  $.
    As in \citep{zhangjianfeng17}, for $p>2$ according to the It\^o formula, we have
    \begin{align*}
        d\|X^d_t\|^2 & = \left(2 X^d_t \cdot a^d(t,X^d_t) + \|b^d(t,X^d_t)\|_{\H}^2 \right) dt + 2(X^d_t b^d(t,X^d_t)) \cdot dW^d_t \\
        d\|X^d_t\|^p & = d\big(\|X^d_t\|^2 \big)^{\frac{p}{2}} = \Big(p\|X^d_t\|^{p-2} X^d_t\cdot a^d(t,X^d_t) + \frac{p}{2}\|X^d_t\|^{p-2} \|b^d(t,X^d_t)\|^2_{\H}  \notag \\
       & + \frac{p(p-2)}{2}\|X^d_t\|^{p-4} \|X^d_t b^d(t,X^d_t)\|^2 \Big) dt  + p\|X^d_t\|^{p-2} (X^d_t b^d(t,X^d_t)) \cdot dW^d_t. 
    \end{align*}
    As in \citep{zhangjianfeng17}, $\int_0^t p\|X^d_s\|^{p-2} \sigma^d(s,X^d_s) d W^d_s $ is a martingale. Therefore, by the property of the Hilbert-Schmit norm, Proposition \ref{prop:dynamic_bound_new}, and the Jensen inequality,
    \begin{align*}
        \big\|X^d_t \big\|^p_{L^p} & \le \|x^d_0\|^p + \mathbb{E}\int_0^t\|X^d_s\|^{p-4} \Big(p\|X^d_s\|^2 |X^d_s\cdot a^d(s,X^d_s)| + \frac{p}{2}\|X^d_s\|^2\|b^d(s,X^d_s)\|^2_\H \\
        & \quad \quad \quad \quad + \frac{p(p-2)}{2}\|X^d_s b^d(s,X^d_s)\|^2 \Big) ds \\
        & \le \|x^d_0\|^p + B^{'}_p (\log d ) \Big(\mathbb{E}\Big[|X^{*,0;d}_T|^{p-1}\int_0^T\|a^d(t,0)\|dt \Big] + \mathbb{E}\Big[|X^{*,0;d}_T|^{p-2}\int_0^T\|b^d(t,0)\|_{\H}^2dt \Big] \\
        & \quad \quad \quad \quad + \int_0^T \big\| X^d_s \big\|^p_{L^p} ds  \Big) \\
    \end{align*}
    for $B^{'}_p =2(p+\frac{p}{2}+ \frac{p(p-2)}{2 }+ 4)C^2  $. By the Gronwall inequality and Young inequality, we have
    \begin{align*}
        \big\|X^d_t \big\|^p_{L^p} &\le \exp(B^{'}_p T \log d ) \Big[\|x^d_0\|^p + B^{'}_p(\log d) \Big(\mathbb{E}[|X^{*,0;d}_T|^{p-1}\int_0^T\|a^d(t,0)\|dt] \\
        & \quad \quad \quad \quad + \mathbb{E}\big[|X^{*,0;d}_T|^{p-2}\int_0^T\|b^d(t,0)\|_{\H}^2dt \big] \Big) \Big] \\ 
        & \le B^{''}_p d^{Q_p} (\log d )\|x^d_0\|^p + \frac{1}{p} \Big( \varepsilon^{-1} (2p-3) \mathbb{E}|X^{*,0;d}_T|^p + [B^{''}_p d^{Q_p} (\log d )]^p \varepsilon^{p-1} \| a^d(\cdot,0) \|^p_{1,p}  \\
        & \quad \quad \quad \quad  + 2[B^{''}_p d^{Q_p} (\log d )]^{\frac{p}{2}} \varepsilon^{\frac{p-2}{2}}  \| b^d(\cdot,0) \|^p_{2,p;\H}   \Big)  \\
        & \le B^{''}_p d^{Q_p} (\log d )\|x^d_0\|^p + \frac{1}{p} \varepsilon^{-1} (2p-3)  \big\|X^{*,0;d}_T \big\|^p_{L^p}  + \Big(2[B^{''}_p d^{Q_p} (\log d )]^{\frac{p}{2}} \varepsilon^{\frac{p-2}{2}} \\
        & \quad \quad \quad \quad + [B^{''}_p d^{Q_p} (\log d )]^p \varepsilon^{p-1} \Big) C^p d^{pQ}
    \end{align*}
    for $B^{''}_p = 1+B^{'}_p, Q_p = B^{'}_p T $. Combining this with Equation \eqref{eq:max_p_X}, we derive
    \begin{align*}
        \big\|X^{*,0;d}_T \big\|^p_{L^p} & \le (B_p+B^{''}_p) d^{Q_p} T(\log d)^{\frac{p}{2}}(1+\|x^d_0\|^p) + \frac{2}{p}B_p  T(\log d)^{\frac{p}{2}} \Big(2[B^{''}_p d^{Q_p} (\log d )]^{\frac{p}{2}}\varepsilon^{\frac{p-2}{2}} \notag \\
        & + [B^{''}_p d^{Q_p} (\log d )]^p \varepsilon^{p-1} \Big) C^p d^{pQ} + \frac{2p-3}{p}B_p  T(\log d)^{\frac{p}{2}} \varepsilon^{-1}  \big\|X^{*,0;d}_T \big\|^p_{L^p} .
    \end{align*}
    By taking $\varepsilon = \frac{2(2p-3)}{p}B_p  T(\log d)^{\frac{p}{2}} $, we obtain
    \begin{align*}
        \big\|X^{*,0;d}_T \big\|^p_{L^p}  
        & \le B^{\circ}_p d^{Q^{\circ}_p} (\log d)^{R^{\circ}_p} (1+\|x^d_0\|^p )
    \end{align*}
    for some positive constants $B^{\circ}_p,Q^{\circ}_p,R^{\circ}_p $. For $p=2$, the argument is much easier following the same procedure of the proof in \citep{zhangjianfeng17}, so we do not show it here. 
\end{proof}

Based on Theorem \ref{theorem:dynamic_bound_new}, we provide a corollary of the $d\times d$-matrix for further reference. Let $X^d$ satisfy
\begin{equation}\label{eq:matrix_SDE}
    X^d_t = x^d_0 + \int_0^t a^d(s,X^d_s) ds +   \int_{0}^t b^d(s,X^d_s) dW^d_s, \; \forall t\in [0,T],
\end{equation}
where $a^d,b^d$ are $\mathbb{R}^{d\times d}, \mathbb{R}^{d\times d\times d} $-valued functions, respectively, with $b^d=(b^d_1,\ldots,b^d_d)$ for $b^d_i\in \mathbb{R}^{d\times d} $.
Then, we have the following theorem.
\begin{theorem}\label{thm:dynamic_X_multidim_tensor}
 Given $p\ge 2$, under Assumption \ref{ass:N_0_structural_ass_model} for the $d\times d$-matrix version for $a^d$ and the $d\times d\times d $-tensor version for $b^d$, 
and the further assumption that $a^d(t,0) \in L^{1,p}(\mathbb{F},\mathbb{P},\mathbb{R}^{d\times d})$, $b^d(t,0) \in L^{2,p}(\mathbb{F},\mathbb{P},\mathbb{R}^{d\times d \times d}) $ is bounded by $ C d^Q $, which is the same as in Assumption \ref{ass:N_0_structural_ass_model},
the following properties of the SDE's solution \eqref{eq:matrix_SDE} holds:
    \begin{equation}
        \mathbb{E}|X^{*,0;d}_T |^p \le B_p d^{Q_p} (\log d)^{R_p} (1+\|x^d_0\|^p_\H ),
    \end{equation}    
    where $B_p,Q_p,R_p $ are constants independent of $d$.
\end{theorem}
\begin{remark}
   The solution $X^d$ in Theorem \ref{thm:dynamic_X_multidim_tensor} can be modified to \textbf{a.s.}-continuous version; thus, $X^d$ can further be predictable, which implies that we can apply the Hilbert space's BDG inequality for $X^d$ developed by \citep{parto14}. Here we simply apply Lemma \ref{lem:BDG_new}. 
\end{remark}
\begin{proof}[Proof of Theorem~\ref{thm:dynamic_X_multidim_tensor}]
The procedure for the proof is similar to the proof used in Theorem \ref{theorem:dynamic_bound_new}. Without loss of generality (by the truncation method), assume $X^d\in L^p(\mathbb{F},\mathbb{P},\mathbb{R}^{d\times d}) $.
    As $\|X^d_t\|_{\H} \le \|x^d_0\|_{\H} + \int_0^t\|a^d(s,X^d_s)\|_{\H} ds + \sum\limits_{i=1}^d \| \int_0^t b^d_i(s,X^d_s) dW^{i;d}_s \|_{\H}  , $
    then, for $p\ge 2$,
    \begin{equation*}
        |X^{*,0;d}_T|^p \le (d+2)^{p-1} \Big(\|x^d_0\|_{\H}^p + \big(\int_0^t\|a^d(s,X^d_s)\|_{\H} ds \big)^p +  \sum_{i=1}^d \sup_{0\le t\le T} \big\| \int_0^t b^d_i(s,X^d_s) dW^{i;d}_s \big\|_{\H}^p \Big).
    \end{equation*}
    As $\| \int_0^t b^d_i(s,X^d_s) dW^{i;d}_s \|_{\H}^p = \big(\sum_{j=1}^d \|\int_0^t b^{j;d}_i(s,X_s)dW^{i;d}_s\|^2 \big)^\frac{p}{2} \le d^{\frac{p-2}{2}} \sum_{j=1}^d \|\int_0^t b^{j;d}_i(s,X^d_s)dW^{i;d}_s\|^p   $ ($b^{j;d}_i $ is the $j$-th column of $b^d_i$), 
    then according to Lemma \ref{lem:BDG_new}, 
    \begin{equation*}
        \mathbb{E}\Big[\sup\limits_{0\le t\le T} \|\int_0^t b^{j;d}_i(s,X^d_s)dW^{i;d}_s\|^p \Big] \le C_p  \| b^{j;d}_i(\cdot,X^d_{\cdot}) \|^p_{2,p} .
    \end{equation*}
    Thus, by Proposition \ref{prop:dynamic_bound_new} for the $d\times d\times d $-tensor version,
    \begin{equation*}
      \big\| X^{*,0;d}_T \big\|^p_{L^p}  \le C_1 d^{Q_1} (\log d )^{\frac{p}{2}} \Big(1+\|x^d_0\|_{\H} + \mathbb{E}\int_0^T\|X^d_t\|_{\H}^p dt \Big)
    \end{equation*}
    for $C_1=2^{2p-2}C_pC^p T^{p-1} , Q_1=\frac{p+2}{2} + p-1  $. By the Hilbert space version of the It\^o formula (e.g., in \citep{parto14}),
    \begin{align*}
        d\|X^d_t\|_{\H}^2 &  = \left[ 2 \Tr( (X^d_t)^{\Transpose} a^d(t,X^d_t)) + \|b^d(t,X^d_t)\|_{\H}^2  \right]dt + 2\sum\limits_{j=1}^d [ (X^{j;d}_t)^{\Transpose} b^{j;d}(t,X^d_t) ] \cdot dW^d_t ,\; \text{and} \\
        d\|X^d_t\|_{\H}^p & = d\big(\|X^d_t\|_{\H}^2 \big)^{\frac{p}{2}} = \Big(p\|X^d_t\|_{\H}^{p-2} \Tr( (X^d_t)^{\Transpose} a^d(t,X^d_t)) + \frac{p}{2}\|X^d_t\|_{\H}^{p-2} \|b^d(t,X^d_t)\|^2_\H  \notag \\
        + & \frac{p(p-2)}{2}\|X^d_t\|_{\H}^{p-4} \Big\|\sum\limits_{j=1}^d [ (X^{j;d}_t)^{\Transpose} b^{j;d}(t,X^d_t) ] \Big\|^2 \Big) dt  + p\|X^d_t\|_{\H}^{p-2} \sum\limits_{j=1}^d [ (X^{j;d}_t)^{\Transpose} b^{j;d}(t,X^d_t) ]\cdot dW^d_t,
    \end{align*}
    where $X^{j;d} $ and $ b^{j;d}_i $ are the $j $-th column of $X^d $ and $b^d_i(i=1,\ldots,d) $, respectively, and $b^{j;d}=(b^{j;d}_1,\ldots,b^{j;d}_d ) $.
    Note that 
    \begin{equation*}
        \Big\|\sum\limits_{j=1}^d [ (X^{j;d}_t)^{\Transpose} b^{j;d} ] \Big\|^2_\H = \sum\limits_{i=1}^d \Big[ \sum\limits_{j=1}^d\sum\limits_{k=1}^d (X^{j;d}_t)^{\Transpose}_{k} b^{j;d}_{k,i}  \Big]^2=\sum\limits_{i=1}^d \Big[\sum\limits_{j=1}^d \sum\limits_{k=1}^d (X^d_t)^{\Transpose}_{j,k} (b^d_i)_{k,j}  \Big]^2 = \sum\limits_{i=1}^d[\Tr( (X^d_t)^{\Transpose} b^d_i )]^2 .
    \end{equation*}
    Since $\Tr(\cdot) $ is a Hilbert inner product of the matrix space, then, by the Cauchy-Schwartz inequality, \(  \big\|\sum_{j=1}^d [ (X^{j;d}_t)^{\Transpose} b^{j;d} ] \big\|^2_\H \le \|X^d \|_{\H}^2 \sum_{i=1}^d \|b^d_i\|^2_\H = \|X\|_{\H}^2 \|b^d(t,X^d_t)\|^2_\H \).
    Therefore, using an argument similar to that in the proof of Theorem \ref{theorem:dynamic_bound_new}, we can calculate the result.  
\end{proof}

We also use the following Lipschitz continuous theorem to map $x \mapsto \mathbb{E}X^{x;d}_t $ for any given $t\in [0,T]$, where $X^{x;d}$ denotes the process starting with the initial value $x$.
\begin{theorem}[Lipschitz continuous for SDE]\label{theorem:Lip_SDE}
    Under Assumption \ref{ass:N_0_structural_ass_model}, given $p\ge 2 $, for any $x_i \in \mathbb{R}^d ,\; i=1,2 $, there exists a positive constant $Q_p $ (chosen to be the same as in Theorem \ref{theorem:dynamic_bound_new}), such that \( \| X^{x_1;d}_t - X^{x_2;d}_t  \|_{L^p} \le  d^{Q_p} \|x_1 - x_2\| \) for all $ t\in [0,T] $. 
    For $x_i \in \mathbb{R}^{d\times d},\; i=1,2 $ of the matrix version's SDE, we replace $\|\cdot \| $ with $\|\cdot\|_{\H} $.
\end{theorem}
\begin{proof}[Proof of Theorem~\ref{theorem:Lip_SDE}]
For the $\mathbb{R}^d $-valued scenario, we have 
     \begin{align*}
        d\|X^{x_1;d}_t - X^{x_2;d}_t\|^2 & = \left(2(X^{x_1;d}_t - X^{x_2;d}_t ) \cdot \Delta a^d_t + \|\Delta b^d_t\|_{\H}^2 \right) dt + 2[(X^{x_1;d}_t - X^{x_2;d}_t) \Delta b^d_t ] \cdot dW^d_t \\
        d\|X^{x_1;d}_t - X^{x_2;d}_t\|^p &  = \Big(p\|X^{x_1;d}_t - X^{x_2;d}_t\|^{p-2} (X^{x_1;d}_t - X^{x_2;d}_t)\cdot \Delta a^d_t + \frac{p}{2}\|X^{x_1;d}_t - X^{x_2;d}_t\|^{p-2} \|\Delta b^d_t\|^2_\H  \notag \\
         + \frac{p(p-2)}{2}\|X^{x_1;d}_t & - X^{x_2;d}_t\|^{p-4} \|(X^{x_1;d}_t - X^{x_2;d}_t) \Delta b^d_t\|^2 \Big) dt  \\
         + p\|X^{x_1;d}_t & -  X^{x_2;d}_t\|^{p-2}  [(X^{x_1;d}_t - X^{x_2;d}_t)\Delta b^d_t] \cdot dW^d_t. 
    \end{align*}
    where $\Delta a^d_t = a^d(t,X^{x_1;d}_t) - a^d(t,X^{x_2;d}_t) $ and $ \Delta b^d_t =b^d(t,X^{x_1;d}_t) - b^d(t,X^{x_2;d}) $. Using the Lipschitz assumption Equation \eqref{eq:N_0_model_Lip} in Assumption \ref{ass:N_0_structural_ass_model}, we have
    \begin{align*}
        \mathbb{E}\|X^{x_1;d}_t - X^{x_2;d}_t\|^p 
        & \le \|x_1 - x_2\|^p + B^{'}_p (\log d)\int_0^t \mathbb{E}\|X^{x_1;d}_s - X^{x_2;d}_s\|^p ds
    \end{align*}
    for $B^{'}_p = pC+\frac{p(p-1)}{2}C^2 $. By the Gronwall inequality,
    \begin{equation*}
        \mathbb{E}\|X^{x_1;d}_t - X^{x_2;d}_t\|^p \le e^{B^{'}_p T (\log d)} \|x_1 - x_2 \|^p \le d^{B^{'}_p T} \|x_1 - x_2 \|^p .
    \end{equation*}
    Taking $Q_p = B^{'}_p T $, we obtain the result for the $\mathbb{R}^d $ case. For the $\mathbb{R}^{d\times d} $-valued scenario, the argument is the same. Thus, we complete the proof by replacing $\|\cdot\| $ with $\|\cdot\|_{\H} $. 
\end{proof}

To prepare the proof of expressivity for $K$, we first consider the following related BSDE:
\begin{equation}\label{eq:BSDE_new}
    Y^d_t = \xi^d - \int_t^T Z^d_s \cdot dW^d_s ,\; 0\le t \le T
\end{equation}
for some random terminal value $\xi^d $. The following theorem bounds the solution $(Y^d,Z^d) $ without dependence on $d$; the well-posedness can be found in \citep{zhangjianfeng17}.
\begin{theorem}[BSDE bound]\label{theorem:Priori_est_new}
    Given $ p\ge 2$, for any $\xi \in L^p(\mathcal{F}_T ; \mathbb{R}) $, let $(Y^d,Z^d) \in S^2(\mathbb{F},\mathbb{P},\mathbb{R}) \times L^{2}(\mathbb{F},\mathbb{P},\mathbb{R}^d )  $ be the unique $\mathbb{F} $-progressively measurable solution of BSDE \eqref{eq:BSDE_new}. Then the following holds: there exists a positive constant $B^{*}_p $ independent of $d$ such that
\begin{equation}
    \|(Y^d,Z^d)\|^p := \big\| Y^{*,0;d}_T \big\|^p_{L^p} + \big\| Z^d \big\|^p_{2,p}   \le B^{*}_p \|\xi^d \|^p_{L^p} .
\end{equation}
If $p\ge 4$, then for any $\xi^d \in L^p(\mathcal{F}_T,\mathbb{P},\mathbb{R}^{d }) $, let $(Y^d,Z^d) \in S^2(\mathbb{F},\mathbb{P},\mathbb{R}^d) \times L^{2,p}(\mathbb{F},\mathbb{P},\mathbb{R}^{d\times d} ) $ be the unique $\mathbb{F} $-progressively measurable solution of BSDE \eqref{eq:BSDE_new}. Then we also have
\begin{equation}
    \|(Y^d,Z^d)\|^p_\H:= \big\| Y^{*,0;d}_T \big\|^p_{L^p} + \big\| Z^d \big\|^p_{2,p;\H}  \le B^{*}_p \|\xi^d \|^p_{L^p} .
\end{equation}
    
\end{theorem}
\begin{proof}[Proof of Theorem~\ref{theorem:Priori_est_new}]
We follow the procedure in \citep{zhangjianfeng17} without considering the Gronwall inequality. According to \citep{zhangjianfeng17}, for $p=2$, it is easy to deduce that 
    \begin{equation*}
    Y^{*,0;d}_T  \le  |\xi^d| + \sup\limits_{0\le t\le T} \Big|\int_0^t Z^d_s \cdot dW^d_s \Big| .
    \end{equation*}
    By applying the BDG inequality Lemma~\ref{lem:BDG_new}, we obtain
    \begin{equation} \label{eq:BDG_apply_new}
        \big\| Y^{*,0;d}_T \big\|^2_{L^2} \le  C_2  (  \|\xi^d \|^2_{L^2} + \| Z \|^2_{\L^2}  ) .
    \end{equation}
    By the It\^o formula, 
    \begin{equation} \label{eq:2_Ito_formula_new}
        d|Y^d_t |^2 = \|Z^d_t \|^2 dt + 2Y^d_t Z^d_t\cdot dW^d_t.
    \end{equation}
    As $\int_0^t Y^d_s Z^d_s \cdot dW^d_s $ is a martingale, we have
    \begin{equation*}
        \mathbb{E}\Big[|Y_t|^2 + \int_t^T\|Z_s\|^2 ds \Big] = \mathbb{E}|\xi|^2, \; \forall t \in [0,T].
    \end{equation*}
    Accordingly, \(  \| Z \|^2_{\L^2} \le 2 C_2 \|\xi^d \|^2_{L^2}  \).
    Therefore, by Equation~\eqref{eq:BDG_apply_new}, \( \big\| Y^{*,0;d}_T \big\|^2_{L^2} \le 2 C_2 \|\xi^d \|^2_{L^2} \),
    which immediately allows us to deduce the final result. 
    For the more general parameter $p>2 $, we first assume that $Y^d $ is bounded and $(\int_0^T \|Z^d_s\|^2 ds)^{\frac{p}{2}} < \infty  $. Then by applying It\^o formula, we obtain
    \begin{equation}\label{eq:p_Ito_formula_new}
        d|Y^d_t |^p = d\big(|Y^d_t|^2 \big)^{\frac{p}{2}} = \frac{p(p-1)}{2} | Y^d_t |^{p-2} \|Z^d_t\|^2 dt + p|Y^d_t|^{p-2} Y^d_t Z^d_t \cdot dW^d_t.
    \end{equation}
    Therefore, by Lemma~\ref{lem:BDG_new} and the Young inequality, we can derive
    \begin{align*}
       \big\| Y^{*,0;d}_T \big\|^p_{L^p}  
        & \le \|\xi^d \|^p_{L^p}  + \frac{p(p-1)}{2}\mathbb{E}\Big[\int_0^T |Y^d_t|^{p-2} \|Z^d_t\|^2 dt \Big] + p\mathbb{E}\Big[\sup\limits_{0\le t \le T} \big|\int_0^t |Y^d_t|^{p-2} Y^dt Z^d_t \cdot dW^d_t  \big| \Big] \\
        & \le \|\xi^d \|^p_{L^p} + \frac{p(p-1)}{2}\mathbb{E}\Big[\int_0^T |Y^d_t|^{p-2} \|Z^d_t\|^2 dt \Big] + \frac{1}{2(pC_2)}pC_2 \big\| Y^{*,0;d}_T \big\|^p_{L^p}  \\
        & + \frac{1}{2}(pC_2)^2 \mathbb{E}\Big[\int_0^T |Y^d_t|^{p-2} \|Z^d_t\|^2 ds   \Big].
    \end{align*}
    Thus, \(         \big\| Y^{*,0;d}_T \big\|^p_{L^p}  \le 2 \|\xi^d \|^p_{L^p} + [p(p-1)+(pC_2)^2 ] \mathbb{E}\big[\int_0^T |Y^d_t|^{p-2} \|Z^d_t\|^2 ds   \big] \).
    Also, by directly using the expectation of the integral version of Equation \eqref{eq:p_Ito_formula_new}, which is similar to the case $p=2$, one can easily show that
    \begin{equation*}
        \mathbb{E}|Y^d_t|^p  + \frac{p(p-1)}{2} \mathbb{E}\Big[\int_t^T |Y^d_t|^{p-2} \|Z^d_t\|^2 ds  \Big] \le \mathbb{E}|\xi^d |^p ,\; \forall t\in [0,T].
    \end{equation*}
    Hence, we have \( \mathbb{E}\big[\int_t^T |Y^d_t|^{p-2} \|Z^d_t\|^2 ds   \big] \le \frac{2}{p(p-1)} \| \xi^d \|_{L^p}^p \).
    We immediately deduce that \( \big\| Y^{*,0;d}_T \big\|^p_{L^p} \le (4+2C_2^2) \| \xi^d \|_{L^p}^p \).
    Then, by Equation~\eqref{eq:2_Ito_formula_new}, we have
    \begin{equation}\label{eq:p_2_BDG_solve}
        \Big(\int_0^T \|Z^d_s\|^2 dt \Big)^{\frac{p}{2}} \le 2^{\frac{p-2}{2}} |\xi^d|^p + 2^{\frac{p}{2}}\Big|\int_0^T Y^d_t Z^d_t \cdot dW^d_t \Big|^{\frac{p}{2}}.
    \end{equation}
    Thus, by Lemma~\ref{lem:BDG_new} and the Young inequality,
    \begin{align*}
         \| Z^d \|_{2,p}^p
        & \le  2^{\frac{p-2}{2}} \| \xi^d \|_{L^p}^p + 2^{\frac{p}{2}} C_{\frac{p}{2}}  \mathbb{E}\Big[ |Y^{*,0}_T |^{\frac{p}{2}} \Big(\int_0^T  \|Z_t\|^2 d_t \Big)^{\frac{p}{4}} \Big] \\
        & \le 2^{\frac{p-2}{2}} \| \xi^d \|_{L^p}^p + 2^{p-1} C_{\frac{p}{2}} ^2 \big\| Y^{*,0;d}_T \big\|^p_{L^p} + \frac{1}{2}  \| Z^d \|_{2,p}^p .
    \end{align*}
    Then, we immediately have \( \| Z^d \|_{2,p}^p \le 2^{\frac{p}{2}} \mathbb{E}|\xi|^p + 2^{p} C_{\frac{p}{2}}^2 \mathbb{E}|Y^{*,0}_T |^{p}\le \left[2^{\frac{p}{2}}+ 2^{p} C_{\frac{p}{2}}^2 (4+2(B^{*,1})^2) \right] \| \xi^d \|_{L^p}^p  \).
    Hence, there exists a positive constant $B^{*}_p$ such that \(  \|(Y^d,Z^d)\|^p \le B^{*}_p \|\xi^d \|^p_{L^p} \).
    
    For the second argument with the parameter $p\ge 4 $, we follow a  similar argument, replacing the norm with $\|\cdot\|_{\H} $ as in the proof of Lemma~\ref{lem:BDG_new}. Note that $p\ge 4 $ guarantees the application of Lemma~\ref{lem:BDG_new} for the term \( \big(\int_0^T (Y_t Z_t)\cdot dW_t \big)^{\frac{p}{2}} \),
     which completes the proof.
\end{proof}

\subsubsection{Expressivity of the focused FBSDE without driver term}


We now return to the equation, non-driver decoupled FBSDE, formulated as follows:
\begin{align}\label{eq:FBSDE_new}
    X^d_t&=x+\int_0^t a^d(s,X^d_s)ds + \int_0^t b^d(s,X^d_s)dW^d_s \notag ,\; \text{and}\\
    Y^d_t&=g^d(X^d_T)  - \int_{t}^T Z^d_s dW^d_s.
\end{align}
Under some regularity conditions, for any $x\in \mathbb{R}^d $ (or $\mathbb{R}^{d\times d}$), the above FBSDE~\eqref{eq:FBSDE_new} has a unique solution $(X^d,Y^d,Z^d)\in L^2(\mathbb{F},\mathbb{P},\mathbb{R}^d)\times S^2(\mathbb{F},\mathbb{P},\mathbb{R}) \times L^2(\mathbb{F},\mathbb{P},\mathbb{R}^d) $ (or $L^2(\mathbb{F},\mathbb{P},\mathbb{R}^{d\times d})\times S^2(\mathbb{F},\mathbb{P},\mathbb{R}^d) \times L^2(\mathbb{F},\mathbb{P},\mathbb{R}^{d\times d}) $) that is $\mathbb{F} $-progressively measurable. By our previous estimation of the SDE's and BSDE's solution, we immediately derive the following estimation theorem for $(Y^d,Z^d) $ of FBSDE~\eqref{eq:FBSDE_new}.
\begin{theorem}[Estimation for focused FBSDE]\label{thm:est_focused_FBSDE}
    Given $p\ge 2$, under Assumptions \ref{ass:N_0_structural_ass_model} and \ref{ass:N_0_structural_ass_g} and the further assumption that $g^d(x,\cdot) \in L^p( \mathcal{F}_T ; \R ) $ with the same bound $C d^Q $ as in Assumption \ref{ass:N_0_structural_ass_g}, then for our focused decoupled FBSDE~\eqref{eq:FBSDE_new} solution $(X^d,Y^d,Z^d)\in L^2(\mathbb{F},\mathbb{P},\mathbb{R}^d)\times S^2(\mathbb{F},\mathbb{P},\mathbb{R}) \times L^2(\mathbb{F},\mathbb{P},\mathbb{R}^d)$ ,
    \begin{equation}
         \|(Y^d,Z^d)\|^p \le C^{*}_p d^{Q^{*}_p} (1+\|x\|^p)
    \end{equation}
    for some positive constants $C^{*}_p,Q^{*}_p $. For solution $(X^d,Y^d,Z^d)\in L^2(\mathbb{F},\mathbb{P},\mathbb{R}^{d\times d})\times S^2(\mathbb{F},\mathbb{P},\mathbb{R}^d) \times L^2(\mathbb{F},\mathbb{P},\mathbb{R}^{d\times d})$, the argument is the same when we replace $\|\cdot\| $ with $\|\cdot\|_{\H} $.
\end{theorem}
\begin{proof}[Proof of Theorem~\ref{thm:est_focused_FBSDE}]
The proof is quite simple. By Theorem~\ref{theorem:Priori_est_new}, we have \(   \|(Y^d,Z^d)\|^p \le B^{*}_p\mathbb{E}|g^d(X^d_T)|^p \).
    According to Proposition~\ref{prop:g_bound_new} and Theorem~\ref{theorem:dynamic_bound_new},
    \begin{equation*}
        \|(Y^d,Z^d)\|^p \le 2^{p-1}B^{*}_p C^p d^{pQ}(1+\mathbb{E}\|X^d_T\|^p) \le C^{*}_p d^{Q^{*}_p} (1+\|x\|^p), 
    \end{equation*}
    where $C^{*}_p =2^{p}B^{*}_p C^pB_p $ and $ Q^{*}_p=pQ+Q_p+R_p  $. The argument for solution $(X^d,Y^d,Z^d)\in L^2(\mathbb{F},\mathbb{P},\mathbb{R}^{d\times d})\times S^2(\mathbb{F},\mathbb{P},\mathbb{R}^d) \times L^2(\mathbb{F},\mathbb{P},\mathbb{R}^{d\times d})$ is the same, which completes the proof. 
\end{proof}

To simplify the analysis of the expressivity of the numerical integration framework for $Z^d$, we clarify the structure of $Z^d$ using Feynman-Kac representation, which requires further regularity assumptions for the FBSDE. First, the coefficient functions must be deterministic. Here, we denote $u^d(t,X^d_t)=Y^d_t$ as in \citep{zhangjianfeng17}, and then provide subsequent FBSDE propositions with expression rates under more regular assumptions.
\begin{proposition} \label{prop:bound_Lip_u_new}
    Under Assumption~\ref{ass:N_0_structural_ass_model_deter} and Assumption~\ref{ass:N_0_structural_ass_g_deter}, there exist constants $b,r > 0$ independent of $d$, such that in $L^2(\mathbb{F},\mathbb{P},\mathbb{R}^d)\times S^2(\mathbb{F},\mathbb{P},\mathbb{R}) \times L^2(\mathbb{F},\mathbb{P},\mathbb{R}^d)$ case, for all $ x_1,x_2 \in \R^d $, 
    \begin{equation} \label{eq:Lip_u_new}
           |u^d(t,x_1)| \le b d^{r} (1+\|x_1\|) , \; \;  | u^d(t,x_1) - u^d(t,x_2) | \le b d^r \|x_1 - x_2\|  . 
    \end{equation}

    Similarly, in $L^2(\mathbb{F},\mathbb{P},\mathbb{R}^{d\times d})\times S^2(\mathbb{F},\mathbb{P},\mathbb{R}^d) \times L^2(\mathbb{F},\mathbb{P},\mathbb{R}^{d\times d})$ case, for all $ x_1,x_2 \in \R^{d\times d} $,
    \begin{equation} \label{eq:Lip_u_matrix_new}
           \|u^d(t,x)\| \le b d^{r} (1+\|x\|_{\H}) ,\; \;  \| u^d(t,x_1) - u^d(t,x_2) \| \le b d^r \|x_1 - x_2\|_{\H} . 
    \end{equation}

  
\end{proposition}
\begin{proof}[Proof of Proposition~\ref{prop:bound_Lip_u_new}]
For the $L^2(\mathbb{F},\mathbb{P},\mathbb{R}^d)\times S^2(\mathbb{F},\mathbb{P},\mathbb{R}) \times L^2(\mathbb{F},\mathbb{P},\mathbb{R}^d)$ scenario,
we denote $(X^{t,x;d}_s,Y^{t,x;d}_s,Z^{t,x;d}_s),\; t\le s\le T $ as the solution of FBSDE~\eqref{eq:FBSDE_new} when the dynamic $X^d$ starts at $t$ with the value $x$. It is easy to verify that $(a^d,b^d) $ in Assumption~\ref{ass:N_0_structural_ass_model_deter} satisfies Assumption~\ref{ass:N_0_structural_ass_model}. Note that $u^d(t,x)=Y^{t,x;d}_t$. 
The first part is proved by Theorem~\ref{thm:est_focused_FBSDE} when $p=2$; we use the monotonicity of the norm with respect to expectation.

For the Lipschitz part, by Theorem~\ref{theorem:Lip_SDE},
\begin{equation}\label{eq:Lip_X}
    \big\| X^{t,x_1;d}_T - X^{t,x_2;d}_T \big\|_{L^2} .
\end{equation}
As $(Y^{t,x_1;d} - Y^{t,x_2;d}, Z^{t,x_1;d}-Z^{t,x_2;d}) $ satisfies the following non-driver linear BSDE,
\begin{equation*}
    \bar{Y}^{d}_l = g^d(X^{t,x_1;d}_T ) - g^d(X^{t,x_2;d}_T ) + \int_l^T \bar{Z}^d_s dW^d_s, \; l\in [t,T],
\end{equation*}
then by Theorem~\ref{theorem:Priori_est_new},
\begin{equation*}
    |u^d(t,x_1)-u^d(t,x_2) | = |\bar{Y}^d_t |\le 
    \| \bar{Y}^{*,t;d}_T \|_{L^1} \le \| \bar{Y}^{*,t;d}_T \|_{L^2}  \le  (B^{*}_2)^{\frac{1}{2}} \| g^d(X^{t,x_1;d}_T) - g^d(X^{t,x_2;d}_T) \|_{L^2} .
\end{equation*}
By Assumption~\ref{ass:N_0_structural_ass_g} and Equation~\eqref{eq:Lip_X}, we have \( \| g^d(X^{t,x_1;d}_T) - g^d(X^{t,x_2;d}_T) \|_{L^2} \le  b_1 d^{r_1} \|x_1-x_2\|  \)
for $b_1=C $ and $ r_1=Q+Q_2$,
and thus can deduce that \(  |u^d(t,x_1)-u^d(t,x_2) | \le b_1 d^{r_1} \|x_1-x_2\| \).
By choosing $b,r$ to be the maximum of $(C^{*}_2)^{\frac{1}{2}},\frac{1}{2}Q^{*}_2$ and $b_1,r_1 $, we complete the proof. The $L^2(\mathbb{F},\mathbb{P},\mathbb{R}^{d\times d})\times S^2(\mathbb{F},\mathbb{P},\mathbb{R}^d) \times L^2(\mathbb{F},\mathbb{P},\mathbb{R}^{d\times d})$ scenario can be solved with the procedure by replacing $\|x\| $ in the above argument with $\|x\|_{\H} $. 
\end{proof}

\begin{proposition} \label{prop:bound_z_direct_new}
    Under Assumption~\ref{ass:N_0_structural_ass_model_deter} and Assumption~\ref{ass:N_0_structural_ass_g_deter}, for the $L^2(\mathbb{F},\mathbb{P},\mathbb{R}^d)\times S^2(\mathbb{F},\mathbb{P},\mathbb{R}) \times L^2(\mathbb{F},\mathbb{P},\mathbb{R}^d) $ scenario, we have
    \begin{equation*}
        \|Z^d_t\| \le b d^r \|b^d(t,X^d_t) \|_{\H} \le bC d^r (\log d)^{\frac{1}{2}} (1+\|X^d_t\|).
    \end{equation*}
    For the $L^2(\mathbb{F},\mathbb{P},\mathbb{R}^{d\times d})\times S^2(\mathbb{F},\mathbb{P},\mathbb{R}^d) \times L^2(\mathbb{F},\mathbb{P},\mathbb{R}^{d\times d})$ scenario, we replace $\|\cdot\| $ with $\|\cdot\|_{\H} $ in the above proposition.
\end{proposition}
\begin{proof}[Proof of Proposition~\ref{prop:bound_z_direct_new}]
For the $L^2(\mathbb{F},\mathbb{P},\mathbb{R}^d)\times S^2(\mathbb{F},\mathbb{P},\mathbb{R}) \times L^2(\mathbb{F},\mathbb{P},\mathbb{R}^d) $ scenario, if $(a^d,b^d,g^d) $ is further continuous differentiable, then according to the Feynman-Kac formula for BSDEs (\citep[Theorem 5.1.4]{zhangjianfeng17}), we know that $Z^d_t=\partial_x u^d(t,X^d_t)b^d(t,X^d_t) $. Then, according to Proposition~\ref{prop:bound_Lip_u_new}, \(  \|\partial_x u^d \| \le b d^r  \).
    By Proposition~\ref{prop:dynamic_bound_new}, for all $ (t,x) \in [0,T] \times \R^d $, \( \|b^d(t,x)\|_{\H} \le C ( \log d)^{\frac{1}{2}}(1+\|x\|) \),
    which immediately leads to the deduction that
    \begin{equation*}
        \|Z^d_t\|\le bC d^r (\log d)^{\frac{1}{2}} (1+\|X^d_t\|).
    \end{equation*}
    For the general $(a^d,b^d,g^d)$, if we choose smooth mollifiers $(a^{\eta;d},b^{\eta;d},g^{\eta;d}),\; 0<\eta<1 $, and denote the related FBSDE solution $(X^{\eta;d},Y^{\eta;d},Z^{\eta;d}) $ ~\eqref{eq:FBSDE_new}, then all of the the previous statements hold for $Z^{\eta;d} $. By using kernel $ K(x) = \exp{(-\frac{1}{1-\|x\|^2})} \mathbf{1}_{(\| x\| < 1) } $, 
    one can easily verify that the growth rate and Lipschitz constants for $(a^{\eta;d},b^{\eta;d},g^{\eta;d}) $ are dominated by $(a^d,b^d,g^d)$'s, therefore
    \begin{equation*}
        \|Z^{\eta;d}_t\|\le bC d^r (\log d)^{\frac{1}{2}} (1+\|X^{\eta;d}_t\|), \; \forall \eta \in (0,1).
    \end{equation*}
    As $X^{\eta;d}=X $ and $\mathbb{E}\int_0^T\| Z^{\eta;d}_t - Z^d_t\|^2 dt \rightarrow 0\; (\eta \rightarrow 0^{+})  $, there exists a $\mathbb{P} $-a.s convergence subsequence $Z^{\eta_n;d} \rightarrow Z^d \; (n\rightarrow \infty) $, therefore, by letting $n\rightarrow \infty $, we have
    \begin{equation*}
        \|Z^d_t\|\le bC d^r (\log d)^{\frac{1}{2}} (1+\|X^d_t\|).
    \end{equation*}
    The argument for the $L^2(\mathbb{F},\mathbb{P},\mathbb{R}^{d\times d})\times S^2(\mathbb{F},\mathbb{P},\mathbb{R}^d) \times L^2(\mathbb{F},\mathbb{P},\mathbb{R}^{d\times d})$ scenario is the same if we replace $\|\cdot\| $ with $\|\cdot\|_{\H} $, which completes the proof. 
\end{proof}


\begin{lemma}[Representation of $Z^d$ by smooth solution]\label{lem:Representation_z_smooth_solution_new}
         Under Assumption~\ref{ass:N_0_structural_ass_model} and Assumption~\ref{ass:N_0_structural_ass_g} with $(a^d,b^d,g^d)$ being continuously differentiable in $(x,y,z) $, we denote the solution $(X,Y,Z) \in L^2(\mathbb{F},\mathbb{P},\mathbb{R}^d)\times S^2(\mathbb{F},\mathbb{P},\mathbb{R}) \times L^2(\mathbb{F},\mathbb{P},\mathbb{R}^d) $ with the related function $u^d(t,X^d_t)=Y_t$ (Markovian). Then $u^d $ is continuous differentiable in $x$ with bounded derivatives that have the expression rate $b d^r$, and
    \begin{equation*}
        \partial_x u^d(t,X^d_t)=\nabla Y^d_t \big(\nabla X^d_t \big)^{-1},\; Z^d_t=\nabla Y^d_t \big(\nabla X^d_t \big)^{-1} b^d(t,X^d_t),
    \end{equation*}
    where $(\nabla X, \nabla Y, \nabla Z) \in L^2(\mathbb{F},\mathbb{P},\mathbb{R}^{d\times d})\times S^2(\mathbb{F},\mathbb{P},\mathbb{R}^d) \times L^2(\mathbb{F},\mathbb{P},\mathbb{R}^{d\times d}) $ is the unique $\mathbb{F}$-progressively measurable solution of the following non-driver decoupled linear FBSDE:
    \begin{align} \label{eq:nabla_FBSDE_new}
        \nabla X^d_t &= \mathbf{I}_d + \int_0^t \partial_x a^d(s,X^d_s) \nabla X^d_s ds +  \sum_{i=1}^d \int_0^t \partial_x b^{i;d}(s,X^d_s) \nabla X^d_s dW^{i;d}_s , \notag\\
        \nabla Y^d_t &=\partial_x g^d(X^d_T) \nabla X^d_T - \int_t^T \nabla Z^d_s dW^d_s ,
    \end{align}
    where $\partial_x b^d = (\partial_x b^{1;d} ,\ldots,\partial_x b^{d;d} ) \in \mathbb{R}^{d\times d\times d} $, $b^{i;d} \in \mathbb{R}^d $ is the $i $-th column of $b^d$. Denote the linear mapping (w.r.t $G$) $ \widetilde{b}^{d}(t,G) := \big(\partial_x b^{1;d}(s,X^d_s) G ,\ldots,\partial_x b^{d;d}(s,X^d_s) G \big) $ for any $ G\in \R^{d\times d} $.

    
\end{lemma}
\begin{proof}[Proof of Lemma~\ref{lem:Representation_z_smooth_solution_new}]
Except for the expression rate, the proof is given in \citep{zhangjianfeng17}. For the expression rate, by directly applying Equation~\eqref{eq:Lip_u_new} in Proposition~\ref{prop:bound_Lip_u_new}, it is easy to verify that \(  \| \partial_x u^d \| \le b d^r  \),
which completes the proof. 
\end{proof}

\subsubsection{Proof of Theorems \ref{thm:Numerical_integration_est} and \ref{theorem:express_N_0_new}}

\begin{proof}[Proof of Theorem~\ref{thm:Numerical_integration_est}]
We first further assume that $(a^d,b^d,g^d) $ are continuously differentiable in $x$, then denote $(X^d,Y^d,Z^d) \in L^2(\mathbb{F},\mathbb{P},\mathbb{R}^d)\times S^2(\mathbb{F},\mathbb{P},\mathbb{R}) \times L^2(\mathbb{F},\mathbb{P},\mathbb{R}^d)$ as the solution of FBSDE~\eqref{eq:FBSDE_new}. This satisfies Lemma~\ref{lem:Representation_z_smooth_solution_new}. Similar to the procedures in \citep{zhangjianfeng17} and \citep{Ma-Zhang02}, we have
    \begin{equation} \label{eq:SDE_derivative_fraction_new}
        (\nabla X^d_t)^{-1} = \mathbf{I}_d - \int_0^t (\nabla X^d_s)^{-1} \Big(\partial_x a^d -  \sum\limits_{i=1}^d (\partial_x b^{i;d} )^2  \Big)(s,X^d_s) ds - 
        \int_{0}^t (\nabla X^d_s)^{-1} \partial_x b^{d}(s,X^d_s) dW^{d}_s .
    \end{equation}

    For further analysis, we first prove some important estimates. Through direct manipulation, we obtain for all $t\in [t_i,t_{i+1}] $,
    \begin{equation*}
        X^d_{t_i,t}:=X^d_{t} - X^d_{t_i} = \int_{t_i}^t a^d(s,X^d_s) ds + \int_{t_i}^t b^d(s,X^d_s) dW^d_s .
    \end{equation*}
    Then, according to Proposition~\ref{prop:dynamic_bound_new}, Theorem~\ref{theorem:dynamic_bound_new}, the Fubini theorem, and the Jensen inequality, 
    \begin{align*}
        \| X^d_{t_i,t} \|_{L^6}
        & \le C (\log d)^{\frac{1}{2}} \Big[2h + h^{\frac{5}{6} } \Big( \int_{t_i}^{t_{i+1}} \| X^d_{s} \|^6_{L^6}  ds \Big)^{\frac{1}{6}} + h^{ \frac{1}{3}} \Big( \int_{t_i}^{t_{i+1}} \| X^d_{s} \|^6_{L^6} ds \Big)^{\frac{1}{6}} \Big] \\
        & \le C (B_6)^{\frac{1}{6}} d^{\frac{1}{6}Q_6} (\log d)^{\frac{1}{2}+\frac{1}{6} R_6} (3h + h^{\frac{1}{2}}) (1 + \|x\| )  \le  B_0 d^{Q_0} h^{\frac{1}{2}} (1+\|x\|)
    \end{align*} 
    for the positive constants $B_0=4 C (B_6)^{\frac{1}{6}} $ and $ Q_0=\frac{1}{6}Q_6 + \frac{1}{2}+\frac{1}{6}R_6 $.

    By the Lipschitz condition in Assumption~\ref{ass:N_0_structural_ass_model_deter}, the coefficient functions in the linear decoupled FBSDE~\eqref{eq:nabla_FBSDE_new} admit the following bounds:
    \begin{equation*}
        \|\partial_x a^d(s,X^d_s) G \|_{\H} \le \| \partial_x a^d(s,X^d_s)\|_{2}  \|G\|_{\H} \le C(\log d)^{\frac{1}{2}} \|G\|_{\H},\; \forall \, G\in \mathbb{R}^{d\times d}, \; \textbf{a.s.},
    \end{equation*}
    \[
        \| \widetilde{b}^{d}(t,G) \|^2_{\H} = \sum_{i=1}^d \| \partial_x {b}^{i;d} G\|^2_{\H}  \le  \| G\|^2_{\H} \sum_{i=1}^d \| \partial_x {b}^{i;d} \|^2_{2}   \le  C^2 (\log d)^{\frac{1}{2}} \| G\|^2_{\H}  ,\; \forall \, G\in \mathbb{R}^{d\times d}, \; \textbf{a.s.} 
    \]
Indeed, \( \| \partial_x a^d(s,X^d_s)\|_{2} \le \Lip a^d   \) and \(  \sum_{i=1}^d \| \partial_x {b}^{i;d} \|^2_{2} \le  \sum_{i=1}^d (\mathrm{Lip}^i (b^{d}))^2 = (\LipH b^d)^2  \).
Consequently, the random affine coefficients $(\omega,t,G) \mapsto \partial_x a^d(t,X^d_t(\omega)) G,\; \widetilde{b}^d(\omega, t, G) $ satisfy Assumption~\ref{ass:N_0_structural_ass_model}. 

In addition, by the Lipschitz condition in Assumption~\ref{ass:N_0_structural_ass_g_deter}, the terminal function in FBSDE~\eqref{eq:nabla_FBSDE_new} can be bounded analogously:
    \begin{equation*}
        \|\partial_x g^d(X^d_T) G \| \le \|\partial_x g^d(X^d_T)\| \|G\|_{\H} \le C d^Q\|G\|_{\H} ,\; \forall \, G\in \mathbb{R}^{d\times d},\; \textbf{a.s.};
    \end{equation*}
therefore, the mapping $(\omega,G) \mapsto \partial_x g^d(X^d_T(\omega)) G $ satisfies Assumption~\ref{ass:N_0_structural_ass_g}. 

    By directly applying Theorems \ref{thm:dynamic_X_multidim_tensor} and \ref{theorem:Priori_est_new}, we obtain
    \begin{equation} \label{eq:nabla_all_est_new}
        \big\| (\nabla X)^{*,0;d}_T \big\|^6_{L^6} + \big\| (\nabla Y_t)^{*,0;d}_T  \big\|^6_{L^6} + \big\| \nabla Z^d \big\|^6_{2,6;\H}  \le  (\bar{B}^{*} d^{\bar{Q}^{*}})^6 
    \end{equation}
    for $\bar{B}^{*}=1+(B_6 + C^{*}_6)^{\frac{1}{6}} ,\;  \bar{Q}^{*}=\frac{1}{6}(Q_6+Q^{*}_6)+1  $. 
    
    Similarly, 
    since \( \big\| \partial_x a^d - \sum_{i=1}^d \partial_x (b^{i;d})^2 \big\|_{2} \le \big\| \partial_x a^d \big\|_{2} + \sum_{i=1}^d \big\| \partial_x b^{i;d} \big\|^2_{2} = \Lip a^d + ( \LipH b^d )^2 \lesssim (\log d)^{\frac{1}{2}} \)   
for \eqref{eq:SDE_derivative_fraction_new}, we can obtain a similar result for $(\nabla X^d)^{-1} $, where we use the same constants as above:
    \( \big\| \big( (\nabla X)^{-1} \big)^{*,0;d}_T \big\|_{L^6}  \le \bar{B}^{*} d^{\bar{Q}^{*}} \) and \( \big\|  (\nabla X^d_t)^{-1}-(\nabla X^d_{t_i})^{-1}\big\|_{L^6;\H} \le \bar{B}^{*} d^{\bar{Q}^{*}}  h^{\frac{1}{2}}  \) for all \( t\in [t_i,t_{i+1}] \).

    According to Lemma~\ref{lem:Representation_z_smooth_solution_new}, we have \(   Z^d_{t_i,t} = \nabla Y^d_t (\nabla X^d_t)^{-1} b^d(t,X^d_t) - \nabla Y^d_{t_i} (\nabla X^d_{t_i})^{-1} b^d(t_i,X^d_{t_i}) \).
    By direct manipulation, this becomes
    \begin{align*}
        \|Z^d_{t_i,t} \| & \le \| \nabla Y^d_t \| \cdot \|  (\nabla X^d_t)^{-1} \|_{\H} \cdot \|b^d(t,X^d_{t}) - b^d(t_i,X^d_{t_i})  \|_{\H}  \\
        &  + \|\nabla Y^d_t\| \cdot \|b^d(t_i,X^d_{t_i})\|_{\H} \cdot \| (\nabla X^d_t)^{-1} - (\nabla X^d_{t_i})^{-1} \|_{\H} + \| \nabla Y^d_t - \nabla Y^d_{t_i} \| \cdot \|(\nabla X^d_{t_i})^{-1} b^d(t_i,X^d_{t_i} )\|_{\H} \\
        & =: I^d_1(t) + I^d_2(t) + I^d_3(t) .
    \end{align*}
    Note that, by the uniformly $\frac{1}{2} $-H\"older continuous assumption in Assumption \ref{ass:N_0_structural_ass_model_deter},
    \begin{equation*}
        \|b^d(t,x) - b^d(s,x)\|_{\H} \le C d^Q |t-s|^{\frac{1}{2}}, \; \forall x\in \mathbb{R}^d,\; t,s\in [0,T].
    \end{equation*}

    We now focus on $t\in [t_i,t_{i+1}] $. For $I^d_1(t)$, by the H\"{o}lder inequality, we deduce that
    \begin{align*}
        \| I^d_1(t) \|_{L^2}^2
        & \le 2 \mathbb{E}\left[ \| \nabla Y^d_t \|^2 \|(\nabla X^d_t)^{-1}\|_{\H}^2 \big(\|b^d(t,X^d_t) - b^d(t_i,X^d_t)\|_{\H}^2 + \|b^d(t_i,X^d_t) - b^d(t_i,X^d_{t_i})\|_{\H}^2   \big) \right] \\
        & \le 2C^2D^{2Q}(\log d ) \big\| \nabla Y^d_t \big\|^2_{L^6} \big\| (\nabla X^d_t)^{-1} \big\|_{L^6;\H}^2
       \big(h + \| X^d_{t} - X^d_{t_i} \|_{L^6}^2 \big) \\
        & \le 4 C^2 (\bar{B}^{*})^{4}(B_0)^2  d^{2Q_0+4\bar{Q}^{*} + 2Q + 1}  h (1+\|x\|^2)  \le \Tilde{B}_1 d^{\Tilde{Q}_1}  (1+\|x\|^2) h
    \end{align*}
    for $\Tilde{B}_1=4 C^2 (\bar{B}^{*})^{4}(B_0)^2 ,\Tilde{Q}_1=2Q_0+4\bar{Q}^{*} + 2Q + 1 $.

    For $I^d_2(t) $, by the H\"{o}lder inequality, we deduce that
    \begin{align*}
        \| I^d_2(t) \|_{L^2}^2   & \le C^2 (\log d) \mathbb{E}\left[ \|\nabla Y_t \|^2 \left(1 + \|X_{t_i} \| \right)^2 \| (\nabla X_t)^{-1} - (\nabla X_{t_i})^{-1}\|_{\H}^2 \right] \\
        & \le 2C^2 (\log d) \big\| \nabla Y^d_t \big\|^2_{L^6}
         \big(1 + \| X^d_{t_i} \|_{L^6}^2  \big) \big\| (\nabla X^d_t)^{-1} - (\nabla X^d_{t_i})^{-1} \big\|_{L^6;\H}^2 
         \le \Tilde{B}_2 d^{\Tilde{Q}_2} (1+\|x\|^2) h 
    \end{align*}
    for $\Tilde{B}_2=4 C^2 (B_6)^{\frac{1}{3}} (\bar{B}^{*})^{4} $ and $ \Tilde{Q}_2=4\bar{Q}^{*}+\frac{1}{3}Q_6 + 1+\frac{1}{3}R_6  $.


    For $I^d_3(t) $, by the H\"{o}lder inequality, Fubini Theorem, and Jensen inequality, 
    \begin{align*}
         \| I^d_3(t) \|_{L^2}^2   & = \mathbb{E}\Big[ \Big\|\int_{t_i}^t \nabla Z^d_s dW^d_s \Big\|^2 \|(\nabla X^d_{t_i})^{-1} b^d(t_i,X^d_{t_i} )\|_{\H}^2   \Big] \\
        & \le 2C^2 (\log d) \mathbb{E}\Big[ \mathbb{E}_{t_i}\Big[ \Big\|\int_{t_i}^t \nabla Z^d_s dW^d_s \Big\|^2  \Big] \|(\nabla X^d_{t_i})^{-1}\|_{\H}^2 \left(1 + \|X^d_{t_i} \|^2 \right)   \Big] \\
        & \le 2C^2 (\log d) \mathbb{E}\Big[ \Big(\int_{t_i}^{t_{i+1}} \|\nabla Z^d_s \|_{\H}^2 ds \Big)  \|(\nabla X^d_{t_i})^{-1}\|_{\H}^2 \left(1 + \|X^d_{t_i} \|^2 \right)   \Big]. 
    \end{align*}
    By combining the results of $I_1,I_2,I_3 $, and \eqref{eq:nabla_all_est_new}, we have 
    \begin{align*}
        & \| Z^d - \overline{Z}^{h;d}  \|^2_{\L^2}
         = \sum\limits_{i=0}^{n-1} \int_{t_i}^{t_{i+1}} \mathbb{{E}}\|Z^d_{t_i,t} \|^2 dt   \le 3 \sum\limits_{i=0}^{n-1} \int_{t_i}^{t_{i+1}} \big(  \| I^d_1(t) \|_{L^2}^2 +  \| I^d_2(t) \|_{L^2}^2 +  \| I^d_3(t) \|_{L^2}^2  \big) dt \\
        & \le 
        3T(\Tilde{B}_1 +\Tilde{B}_2 )d^{\Tilde{Q}_1+\Tilde{Q}_2}(1+\|x\|^2)h + 6T C^2 d h \mathbb{E}\Big[ \big| \big((\nabla X)^{-1} \big)_T^{*,0;d} \big|^2 \big(1 + \| X^{*,0;d}_T \|^2 \big) \int_{0}^{T} \|\nabla Z^d_s\|_{\H}^2 ds \Big] \\
        & \le 3T(\Tilde{B}_1 +\Tilde{B}_2 )d^{\Tilde{Q}_1+\Tilde{Q}_2}(1+\|x\|^2)h + 6TC^2 d h  \big\| \big( (\nabla X)^{-1} \big)^{*,0;d}_T \big\|^2_{L^6;\H} \big( 1 + \| X^{*,0;d}_T \|_{L^6}^2 \big)  \big\| \nabla Z^d  \big\|_{2,6;\H}^2  \\
        & \le \Tilde{B} d^{\Tilde{Q}} (1+\|x\|^2) h
    \end{align*}
    for $\Tilde{B}=6T(1+C^2)(\Tilde{B}_1 +\Tilde{B}_2 +(B_6)^{\frac{1}{3}} (\bar{B}^{*})^{4} ) $ and $\Tilde{Q}= 1 + \Tilde{Q}_1+\Tilde{Q}_2 + 4\bar{Q}^{*}+\frac{1}{3}Q_6+ \frac{1}{3}R_6$, which completes the proof for the smooth $(a^d,b^d,g^d)$ solution. 
    
    For a more general $(a^d,b^d,g^d)$ solution, under Assumptions~\ref{ass:N_0_structural_ass_model_deter} and \ref{ass:N_0_structural_ass_g_deter}, we choose mollifiers $(a^{\eta;d},b^{\eta;d},g^{\eta;d}) , \; 0<\eta<1  $ that are continuously differentiable in $x$; then, according to the previous argument, there is a solution $(X^{\eta;d},Y^{\eta;d},Z^{\eta;d}) $ for the FBSDE~\eqref{eq:FBSDE_new} that satisfies the previous argument. Thus, the theorem holds for $Z^{\eta;d},\; 0\le \eta \le 1 $, where constants $\bar{B},\bar{Q} $ only depend on $(a^d,b^d,g^d)$. Indeed, the smooth mollifers $(a^{\eta;d},b^{\eta;d},g^{\eta;d}) $ can be generated by the kernel \( K(x) = \exp{(-\frac{1}{1-\|x\|^2})} \mathbf{1}_{(\|x\| < 1 )} \), 
    where $\int_{\mathbb{R}^d} K(x) dx = \int_{[-1,1]^d } K(x) dx =1 $. Then, by setting the mollifiers as \( a^{\eta;d}(t,x) = \int_{\mathbb{R}^d } a^d(t,x-\eta y) K(y) dy \) (mollifiers for $ b^d ,g^d  $ are similar), 
    one can easily verify that all of the bounding constants of $(a^{\eta;d},b^{\eta;d},g^{\eta;d}) $ can be controlled $(a^d,b^d,g^d)$'s. By Assumptions~\ref{ass:N_0_structural_ass_model_deter} and \ref{ass:N_0_structural_ass_g_deter}, the constants $(a^{\eta;d},b^{\eta;d},g^{\eta;d}) $ can be bounded by those in Assumptions~\ref{ass:N_0_structural_ass_model_deter} and \ref{ass:N_0_structural_ass_g_deter}. Therefore the expressivity result for $Z^{\eta;d},\; 0<\eta<1 $ holds for the bounding constants independent of $\eta $. 
    Similarly, the $L^2$ error $\mathbb{E}[\int_0^T \|Z^{\eta;d}_t-Z^d_t \|^2 dt]$ can be bounded by the error $\mathbb{E}|g^{\eta;d}(X^{\eta;d}_T)- g^d(X^d_T)  |^2 $ according to Theorem~\ref{theorem:Priori_est_new}, which goes to $0$ as $\eta \to 0^{+} $ as $X^{\eta;d} \to X^d $ in $L^2 $ and $g^{\eta;d} \to g^d (\eta \to 0^{+}) $ is uniform ($g^d$ Lipschitz continuous). By denoting the $\eta$-independent expression rate bound for all $Z^{\eta;d},0<\eta<1 $ as
    \begin{equation*}
        \mathbb{{E}}\Big[\sum_{i=0}^{n-1} \int_{t_i}^{t_{i+1}}\|Z^{\eta}_t - Z^{\eta}_{t_i} \|^2 dt \Big] \le \Tilde{B}d^{\Tilde{Q}} (1+\|x\|^2) h, \; \forall \, \eta \in (0,1)
    \end{equation*}
    and letting $\eta \to 0^{+}$, we have \( \big\| Z^d - \overline{Z}^{h;d}  \big\|_{\L^2}^2 \le \Tilde{B}d^{\Tilde{Q}} (1+\|x\|^2) h \),
    which completes the proof. 
\end{proof}

To prove Theorem~\ref{theorem:express_N_0_new}, we need to clarify that $V^d_{n+1},\; n \in \N^{-1} $ satisfy Assumption~\ref{ass:N_0_structural_ass_g_deter}, which is guaranteed by the following proposition.
\begin{proposition}\label{prop:V_Lip_growth}
    If $(a^d,b^d) $ satisfies Assumption~\ref{ass:N_0_structural_ass_model_deter} and $g^d(t_n,\cdot),\; n \in \N $ satisfies Assumption~\ref{ass:N_0_structural_ass_g_deter}, then for all $n \in \N^{-1} $, $V^d_{n+1}  $ satisfies Assumption~\ref{ass:N_0_structural_ass_g_deter}.
\end{proposition}
\begin{proof}[Proof of Proposition~\ref{prop:V_Lip_growth}]
We use backward induction for this proof. The base case $k=N $ is obvious by $V^d_N(x)=g^d(T,x) $. Suppose  Assumption~\ref{ass:N_0_structural_ass_g_deter} holds for $V^d_{n+2} $. Let $X^{t_{n+1}, x;d}_{t_{n+2}} $ denote the dynamic process starting at $t_{n+1} $ with value $x $. For $k=n+1$, by the relationship
    \[ V^d_{n+1}(x)=\max\big(g^d(t_{n+1},x), \mathbb{E}[V^d_{n+2}(X^{t_{n+1}, x ;d}_{t_{n+2}}) ] \big) ,\]
    we have
    \begin{align*}
        |V^d_{n+1}(x) | & \le |g^d(t_{n+1},x) | + \big| \mathbb{E}[V^d_{n+2}(X^{t_{n+1}, x;d}_{t_{n+2}}) ] \big|  \le C d^Q (1+ \|x\| ) + C d^Q  \| X^{t_{n+1}, x;d }_{t_{n+2}} \|_{L^2} 
    \end{align*}
    as well as 
    \begin{align*}
        |V^d_{n+1}(x) - V^d_{n+1}(y) | 
        & \le 2|g^d(t_{n+1},x) -g^d(t_{n+1},y) | +  \left| \mathbb{E}[V^d_{n+2}(X^{t_{n+1}, x;d }_{t_{n+2}}) ] - \mathbb{E}[V^d_{n+2}(X^{t_{n+1}, y;d }_{t_{n+2}}) ]   \right| \\
        & \le 2 C d^Q\|x-y\| + Cd^Q \big\| X^{t_{n+1}, x;d }_{t_{n+2}} - X^{t_{n+1}, y;d }_{t_{n+2}} \big\|_{L^2}  .
    \end{align*}
    By Theorems~\ref{theorem:dynamic_bound_new} and ~\ref{theorem:Lip_SDE}, we then have \( |V^d_{n+1}(x) | \le \bar{c}_{n+1} d^{\bar{q}_{n+1}}(1+\|x\|) \) and \( |V^d_{n+1}(x) - V^d_{n+1}(y) | \le \bar{c}_{n+1} d^{\bar{q}_{n+1}}\|x-y\| \) 
    for $\bar{c}_{n+1}=3 C(B_2)^{\frac{1}{2}} $ and $ {\bar{q}_{n+1}}=Q+Q_2+\frac{1}{2}R_2 $. By induction, after choosing the same constants as in Assumption~\ref{ass:N_0_structural_ass_g_deter} (e.g., taking the maximum), we complete the proof. 
\end{proof}

\begin{proof}[Proof of Theorem~\ref{theorem:express_N_0_new}]
Let $ \overline{Z}^{*,K,n;d} := \sum_{k=0}^{K-1} {Z}^{*;d}_{t^n_k} \mathbf{1}_{[t^n_k,t^n_{k+1})} (t) $ for $ t\in [t_n,t_{n+1}] $. 
By Theorem~\ref{thm:Numerical_integration_est}, there exist constants $\bar{B}_n, \bar{Q}_n > 0,\; n \in \N^{-1} $ independent of $d$, such that \( \big\| {Z}^{*;d} - \overline{Z}^{*,K,n;d} \big\|^2_{[t_n,t_{n+1}],2} \le  \bar{B}_n d^{\bar{Q}_n}  \frac{1}{K} \), 
which immediately becomes \( \big\| {Z}^{*;d} - \widehat{Z}^{*,K,n;d} \big\|^2_{[t_n,t_{n+1}],2} \le  \bar{B}_n d^{\bar{Q}_n}  \frac{1}{K} \) 
    through the minimization property of the conditional expectation. Therefore, let ${B}^{*}=\max_{n\in \N^{-1} }(\bar{B}_n) $ and $ {Q}^{*}=\max_{n\in \N^{-1} }(\bar{Q}_n) $ for any $\varepsilon>0 $. Then, by taking $ K_{d,\varepsilon} =\lceil {B}^{*} d^{{Q}^{*}} \varepsilon^{-1} \rceil  $, we have,
    \begin{equation*}
        \big\| \widehat{Z}^{*,K_{d,\varepsilon} ,n;d} - {Z}^{*;d} \big\|^2_{[t_n,t_{n+1}],2}  \le \varepsilon ,\; \;   K_{d,\varepsilon} \le {B}^{*} d^{{Q}^{*}}\varepsilon^{-1} + 1 ,\; \forall \; n\in \N^{-1} .
    \end{equation*}
    By choosing the same constants as above, we complete the proof.
\end{proof}

\section{Detailed Proofs for Section~\ref{sec:deep_mtg}}


\subsection{Detailed proof of the representation of $Z^{*}_{t^n_k}$}
The proof of Lemma~\ref{lemma:represent_z_x} requires the following lemma.

\begin{lemma} \label{lemma:independent_measurable_relation}
    In a probability space $(\Omega,\mathcal{F},\mathbb{P}) $, for the $\mathcal{B}(\mathbb{R}^d) \otimes \mathcal{F} $-measurable function $f: \mathbb{R}^d \times \Omega \rightarrow \mathbb{R} $, if for any given $x\in \mathbb{R}^d$, $\omega \mapsto f(x,\omega) $ is independent of $\sigma $-field $\mathcal{G} \subset \mathcal{F}$ and $X\in \mathbb{R}^d $ is measurable w.r.t. $\mathcal{G}$ with $\mathbb{E}|f(X,\cdot) |<\infty $, then (a). \( \mathbb{E}[ f(X,\cdot) | \mathcal{G} ] = \mathbb{E}[ f(X,\cdot) | X ],\; \mathbb{P} \)-a.s., and (b). \(  \mathbb{E}[ f(X,\cdot) | \mathcal{G} ] = \mathbb{E}[ f(y,\cdot) | \mathcal{G} ]\big|_{y=X} = g(X,*_{\omega} ) ,\; \mathbb{P}  \)-a.s. for some $\mathcal{B}(\mathbb{R}^d) \otimes \mathcal{F} $-measurable function $g$.
\end{lemma}
\begin{proof}[Proof of Lemma~\ref{lemma:independent_measurable_relation}]
For any $\mathcal{B}(\mathbb{R}^d) \otimes \mathcal{F} $-measurable non-negative function $f $, let
\[
f_n = \sum_{k=0}^{n\cdot2^n - 1} \frac{k}{2^n} 1_{(k\le f< k+1 ) } + n 1_{(f\ge n )}, \; n=1,2,\ldots ,
\]
then $0 \le f_n \le n $ and $f_n \uparrow f $ pointwisely. As $f_n $ is bounded, following the argument in \citep{SDE03} Theorem 7.1.2 and the boundedness of $f_n$, 
\begin{equation}\label{eq:conditional_expectation_relationship_n}
    \mathbb{E}[f_n(X,\cdot)| \mathcal{G}] = \mathbb{E}[f_n(X,\cdot)| X] = \mathbb{E}[f_n(y,\cdot)]\big|_{y=X},\; \forall \; n .
\end{equation}
According to the standard conditional expectation argument (e.g., \citep{modern-Probability}), $\mathbb{E}[f_n(X,\cdot)| \mathcal{G}] \uparrow \mathbb{E}[f(X,\cdot)| \mathcal{G}]$ and $ \mathbb{E}[f_n(X,\cdot)| X] \uparrow \mathbb{E}[f(X,\cdot)| X] \; \text{a.s.} $ Furthermore, $\mathbb{E}[f_n(y,\cdot)] \uparrow \mathbb{E}[f(y,\cdot)] , \forall y\in \mathbb{R}^d$, which implies $\mathbb{E}[f_n(y,\cdot)]\big|_{y=X} \uparrow \mathbb{E}[f(y,\cdot)]\big|_{y=X} $. Thus, by taking the limit of \eqref{eq:conditional_expectation_relationship_n}, we obtain
\begin{equation}\label{eq:conditional_expectation_relationship}
    \mathbb{E}[f(X,\cdot)| \mathcal{G}] = \mathbb{E}[f(X,\cdot)| X] = \mathbb{E}[f(y,\cdot)]\big|_{y=X} ,\; \P\text{-a.s.}
\end{equation}
For $f$, which has the integrable condition in Lemma~\ref{lemma:independent_measurable_relation}, by $f = f^{+} - f^{-} $, where $f^{+},f^{-} \ge 0 $, \eqref{eq:conditional_expectation_relationship} also holds for $f^{+},f^{-}$. Then, by the linearity of the conditional expectation and integrable condition, \eqref{eq:conditional_expectation_relationship} holds for $f$.
\end{proof}

\subsection{Detailed proof of convergence}



%

\begin{proof}[Proof of Theorem~\ref{theorem:inte_L2_approx}]
Since $ \mu_n^K ,\; n\in \N^{-1} $ are finite Borel measure on $ \R^{1+d} $, 
by applying the universal approximation theorem in $L^2(\mu_n) $ \citep{hornik91} with bounded and non-constant activation function, the approximation statement for $ \widehat{Z}^{*,K,n} $ holds, which can be constructed as a shallow network.

    For $ \xi^{*,K}_n $, we have the following relationship by It\^o isometry and Lemma~\ref{lemma: inte_relation_mu}: for all $n \in \N^{-1} $,
    \[
        \big\| \xi^{\theta_n,K}_n  - \xi^{*,K}_n \big\|_{L^2} = \Big( \E \Big| \sum\limits_{k=0}^{K-1} \big( \widehat{Z}^{*,K,n}(t^n_k,X_{t^n_k})-z^{\theta_n}_n(t^n_k,X_{t^n_k}) \big) \cdot \Delta W_{t^n_k} \Big|^2 \Big)^{\frac{1}{2}} = \big\| z^{\theta_n}_n - \widehat{Z}^{*,K,n} \big\|_{2;\mu^K_n} ,
    \]
    which can easily lead to our desired result.
\end{proof}

\begin{proof}[Proof of Theorem~\ref{theorem:theta_approx}]
For $n=N$, it is obvious that $\widetilde{U}_{N}(M^{\theta^K_{\varepsilon} })\equiv g(t_N,X_{t_N})\equiv \widetilde{U}_N(\widehat{M}^{*,K}) $. By backward induction, suppose 
\( \big\| \widetilde{U}_{n+1}(M^{\theta^K_{\varepsilon},K} ) -  \widetilde{U}_{n+1}(\widehat{M}^{*,K})\big\|_{L^2} \le (N-n-1)\varepsilon \).
    By Theorem~\ref{theorem:inte_L2_approx}, there exists a parameter $\theta^{K}_{n,\varepsilon} \in \Theta $ such that \( \big\| \xi^{\theta^K_{n,\varepsilon},K}_n - \xi^{*,K}_n \big\|_{L^2} < \varepsilon \). 
    Then, by applying Lemma~\ref{lem-error-propagate} and the Cauchy-Schwartz inequality, we have
    \[
        \big\| \widetilde{U}_{n}(M^{\theta^K_{\varepsilon},K} ) -  \widetilde{U}_{n}(\widehat{M}^{*,K})\big\|_{L^2} \le (N-n-1)\varepsilon + \big\| \xi^{\theta^K_{n,\varepsilon},K}_n - \xi^{*,K}_n \big\|_{L^2} \le (N-n)\varepsilon
    \]
    as claims.
\end{proof}

\begin{proof}[Proof of Corollary~\ref{coro:tight_upper}]
By Theorem~\ref{theorem:discrete_convergence}, for any $\varepsilon >0 $, there exists an $ K(\varepsilon) \in \mathbb{N} $ such that, for any $ K \ge K(\varepsilon) $, we have \( \big\| Y^{*}_n - \widehat{Y}^{*,K}_n \big\|_{L^1} \le  \big\| Y^{*}_n - \widehat{Y}^{*,K}_n \big\|_{L^2} \le \frac{1}{2} \varepsilon \) for all $ n\in \N^{-1} $.
Moreover, by Theorem~\ref{theorem:theta_approx} and Jensen's inequality, there exists $\theta^{K}_{\varepsilon} \in \Theta^N $ such that, for all $ n\in \N^{-1} $,
    \[
        \big\|  \widehat{Y}^{*,K}_n - \widetilde{U}_n(M^{\theta^{K}_{\varepsilon},K}) \big\|_{L^p} 
        \le  \big\| \widetilde{U}_n(\widehat{M}^{*,K}) - \widetilde{U}_n(M^{\theta^{K}_{\varepsilon},K}) \big\|_{L^2} \le \frac{1}{2} \varepsilon ,\; p=1,2 .
    \]
Combining the above with the fact that $\mathbb{E}[\widetilde{U}_n(M^{\theta^{K},K}) ]=\mathbb{E}\big[\mathbb{E}[\widetilde{U}_n(M^{\theta^{K},K}) | \mathcal{F}_{t_n}] \big] $ is an upper bound for $ \mathbb{E}[Y^{*}_n] $, we obtain, for all $ n\in \N^{-1} $,
     \[
         0\le \mathbb{E}[\widetilde{U}_n(M^{\theta^{K}_{\varepsilon},K}) - Y^{*}_n ] = \left| \mathbb{E}[\widetilde{U}_n(M^{\theta^{K},K}) - Y^{*}_n ] \right| \le \big\| Y^{*}_n - \widehat{Y}^{*,K}_n \big\|_{L^1} + \big\|  \widehat{Y}^{*,K}_n - \widetilde{U}_n(M^{\theta^{K}_{\varepsilon},K}) \big\|_{L^1} \le \varepsilon ,
     \]
which yields $ \mathbb{E}[Y^{*}_n] \le \inf_{\theta \in \Theta^N } \mathbb{E}[\widetilde{U}_n(M^{\theta,K}) ] \le \mathbb{E}[\widetilde{U}_n(M^{\theta^{K}_{\varepsilon},K}) ] \le \mathbb{E}[Y^{*}_n] + \varepsilon ,\; K \ge K(\varepsilon)  $. Hence, $ \lim_{K \to \infty} \inf_{\theta \in \Theta^N } \mathbb{E}[\widetilde{U}_n(M^{\theta,K}) ] = \mathbb{E}[Y^{*}_n] $. 

Moreover, if $ Y^{*}_n \ge 0 $, then $ \mathbb{E}|\widetilde{U}_n(M^{\theta^{K},K}) |^2  $ is an upper bound for $ \mathbb{E}|Y^{*}_n|^2 $. Therefore, for all $ n\in \N^{-1} $,
    \[
         0 \le (\mathbb{E}|\widetilde{U}_n(M^{\theta^{K},K})|^2)^{\frac{1}{2}} - (\mathbb{E} |Y^{*}_n |^2)^{\frac{1}{2}} \le   \big\|  {Y}^{*,K}_n - \widetilde{U}_n(M^{\theta^{K}_{\varepsilon},K}) \big\|_{L^2} \le \varepsilon , 
    \]
and thus for all $ K \ge K(\varepsilon) $, 
\[
    \begin{aligned}
        \big(\mathbb{E}\big|Y^{*}_n\big|^2 \big)^{\frac{1}{2}} & \le  \big( \inf_{\theta \in \Theta^N } \mathbb{E}\big|\widetilde{U}_n(M^{\theta, K}) \big|^2 \big)^{\frac{1}{2}} = \inf_{\theta \in \Theta^N} \big(\mathbb{E}\big|\widetilde{U}_n(M^{\theta, K}) \big|^2 \big)^{\frac{1}{2}} \\
        & \le  \big(\mathbb{E}\big| \widetilde{U}_n(M^{\theta^{K},K}) \big|^2 \big)^{\frac{1}{2}} \le \big(\mathbb{E}\big|Y^{*}_n\big|^2 \big)^{\frac{1}{2}} + \varepsilon .
    \end{aligned}
\]
Consequently,
    \[
        \Big( \lim_{K\to \infty} \inf_{\theta \in \Theta^N} \mathbb{E}\big|\widetilde{U}_n(M^{\theta, K}) \big|^2 \Big)^{\frac{1}{2}} = \lim_{K\to \infty} \Big( \inf_{\theta \in \Theta^N} \mathbb{E}\big|\widetilde{U}_n(M^{\theta, K}) \big|^2 \Big)^{\frac{1}{2}} = \big(\mathbb{E}\big|Y^{*}_n\big|^2 \big)^{\frac{1}{2}} ,
    \]
which implies $ \lim_{K\to \infty} \inf_{\theta \in \Theta^N} \mathbb{E}\big|\widetilde{U}_n(M^{\theta, K}) \big|^2 = \mathbb{E}\big|Y^{*}_n\big|^2 $ as claimed.
\end{proof}


\begin{proof}[Proof of Proposition~\ref{pro:L2-to-L1-loss}]
Since $ \theta^{K} \in \argmin_{\varepsilon_{K}} \mathbb{E}\big|\widetilde{U}_n(M^{\theta, K}) \big|^2  $, we have
\[
    \mathbb{E}\big|\widetilde{U}_n(M^{\theta^{K}, K}) \big|^2 \le \inf_{\theta\in \Theta^N } \mathbb{E}\big|\widetilde{U}_n(M^{\theta, K}) \big|^2 + \varepsilon_{K} .
\]
Since $ \varepsilon_{K} \downarrow 0 $, it follows that \(  \limsup_{K \to \infty} \mathbb{E}\big|\widetilde{U}_n(M^{\theta^{K}, K}) \big|^2 = \mathbb{E}\big|Y^{*}_n\big|^2 \).
Together with
\[
    \liminf_{K \to \infty} \mathbb{E}\big|\widetilde{U}_n(M^{\theta^{K}, K}) \big|^2 \ge \mathbb{E}\big|Y^{*}_n\big|^2 ,
\]
we obtain \(  \lim_{K \to \infty} \mathbb{E}\big|\widetilde{U}_n(M^{\theta^{K}, K}) \big|^2 = \mathbb{E}\big|Y^{*}_n\big|^2  \).
Furthermore, using
    \[
        0 \le   \mathbb{E}\big[ \mathrm{Var}_n\big( \widetilde{U}_n(M^{\theta^{K},K}) \big) \big]  = \mathbb{E}\big[ \mathbb{E}_n\big| \widetilde{U}_n(M^{\theta^{K},K}) \big|^2 \big] - \mathbb{E}\big| \mathbb{E}_n [\widetilde{U}_n(M^{\theta^{K},K})] \big|^2 \le \mathbb{E}\big| \widetilde{U}_n(M^{\theta^{K},K}) \big|^2 - \mathbb{E}\big| Y^{*}_n \big|^2 ,
    \]
we deduce $ \lim_{K \to \infty} \mathbb{E}\big[ \mathrm{Var}_n\big( \widetilde{U}_n(M^{\theta^{K},K}) \big) \big] = 0 $.
Therefore, by \citep[Corollary~5.4]{schoen13} and the uniform integrability of $ (M^{\theta^{K},K})_{K\ge 1} $, we conclude that
    \[
       \lim_{K\to \infty}  \mathbb{E}\big[ \widetilde{U}_n(M^{\theta^{K},K}) ] = \mathbb{E} [Y^{*}_n]  
    \]
as claimed.
\end{proof}

\subsection{Detailed proof of the expressivity of the value function approximation}

\subsubsection{Detailed proof of the infinite-width neural network and RanNN}

\begin{proof}[Proof of Proposition~\ref{prop:contin_RKBS}]
It is sufficient to show that each $f_i$ is continuous. By \citep{RKBS2024}, 
    we have
    \begin{equation*}
        f_i(x,\mu_{i-1}) = \phi_i(x)(\mu_{i-1}),
    \end{equation*}
    where $\phi_i(x) $ is the bounded linear operator from $\mathcal{M}(\Theta_{i-1},\mathcal{X}_i) $ to $\mathcal{X}_i $.
    For $x_1,x_2 \in \mathcal{X}_{i-1}$, $ \mu_{i-1}^1, \mu_{i-1}^2 \in \mathcal{M}(\Theta_{i-1}, \mathcal{X}_i ) $,
    \begin{align*}
    \|f_i(x_1, \mu_{i-1}^1 ) - f_i(x_2, \mu_{i-1}^2)\|_{\mathcal{X}_{i}} 
    &\leq \|f_i(x_1, \mu_{i-1}^1 ) - f_i(x_2, \mu_{i-1}^1 )\|_{\mathcal{X}_{i}} + \|f_i(x_2, \mu_{i-1}^1 ) - f_i(x_2, \mu_{i-1}^2 )\|_{\mathcal{X}_{i}} \\
    &\leq \|\mu_{i-1}^1 \|_{\text{TV}} \|x_1-x_2\|_{\mathcal{X}_{i-1}} + \|\phi_i(x_2)(\mu_{i-1}^{1}-\mu_{i-1}^{2})\|_{\mathcal{X}_{i}} ,\; \text{and} \\
    \|\phi_i(x_2)(\mu_{i-1}^{1}-\mu_{i-1}^{2}))\|_{\mathcal{X}_{i}} & = \Big\|\sum_{n\geq 0} \rho(x_2, n)(w^{i+1,1}_n-w^{i+1,2}_n) \Big\|_{\mathcal{X}_{i}} \\
    &\leq \|b^{i+1}_1 - b^{i+1}_2\|_{\mathcal{X}_{i}} + \left\|W^{i+1}_1-W^{i+1}_2 \right\|_{\mathcal{B}(\mathcal{X}_{i-1}, \mathcal{X}_{i})} \|x_2\|_{\mathcal{X}_{i-1}} \\
    &\leq \left(1 + \|x_2\|_{\mathcal{X}_{i-1}} \right)\|\mu_{i-1}^{1}-\mu_{i-1}^{2}\|_{\text{TV}},
    \end{align*}    
where $\mathcal{B}(\mathcal{X}_{i-1}, \mathcal{X}_{i}) $ denotes the bounded linear operator Banach Space. Thus, 
\begin{equation*}
    \|f_i(x_1, \mu_{i-1}^1 ) - f_i(x_2, \mu_{i-1}^2)\|_{\mathcal{X}_{i}} \le \|\mu_{i-1}^1 \|_{\text{TV}} \|x_1-x_2\|_{\mathcal{X}_{i-1}} + \left(1 + \|x_2\|_{\mathcal{X}_{i-1}} \right)\|\mu_{i-1}^{1}-\mu_{i-1}^{2}\|_{\text{TV}},
\end{equation*}
which implies the continuity.
\end{proof}

\begin{proof}[Proof of Proposition~\ref{prop:size_measurable}]
The measurability of $\Growth(\Tilde{f}) and \Lip(\Tilde{f}) $ is obvious due to the continuity of $f $ w.r.t. $x\in \mathbb{R}^{d_1}$ and $\mu\in \mathcal{U}$. For $\size(\Tilde{f})$, it is sufficient to check $\size\left(f_i\left(*,\mu_{i-1}(\cdot)\right)\right) $. Given parameter $\mu^{'}_{i-1}=\sum_{m=0}^{M} w^{i+1}_m \delta_m \in \mathcal{M}(\Theta_{i-1},\mathcal{X}_i) $, \(   \size\left(f_i(*,\mu^{'}_{i-1})\right) = \sum_{m=0}^{M}\sum_{n=1}^{M^{'}} \mathbf{1}_{(w^{i+1}_{mn} \neq 0)} \),
where $M^{'}=d_3 $ if $i=I+1 $, and $M^{'}= \infty $ otherwise, and $w^{i+1}_{mn} $ is the $n $-th component of $w^{i+1}_m $. The measurability is then clear by the indicator function $ \mathbf{1}_{(w^{i+1}_{mn}\neq 0)} $.
\end{proof}

\subsubsection{Preliminary result for the value function approximation}

By Assumption~\ref{ass:g_new}, we immediately have the following proposition.
\begin{corollary}[Growth Rate and Lipschitz Rate for $g^d$]\label{coro:g_growth_bound_new}
    Under Assumption~\ref{ass:g_new}, for any $n \in \N $, $\varepsilon>0 $, the following inequalities hold: \(  \Growth(g^d(t_n,\cdot)),\; \Lip(g^d(t_n,\cdot)) ,\;  \Growth(\widehat{g}^d_n) \le c d^{q}   \).
\end{corollary}


\begin{proof}[Proof of Corollary~\ref{coro:g_growth_bound_new}]
The proof of the linear growth rate can be found in \citep{gonon23}. For the Lipschitz growth rate, for any $n\in \N, \; x,y \in \mathbb{R}^d,\; \varepsilon >0 $, there exists $\widehat{g}^d_n $ such that $|\widehat{g}^d_n(x) - g^d(t_n,x)|\le \varepsilon c d^q (1+\|x\|) $ with $|\widehat{g}^d_n(x) - \widehat{g}^d_n(y)| \le c d^q\|x-y\| $. Then, 
\begin{align*}
     |g^d(t_n,x) - g^d(t_n,y)| & \le |g^d(t_n,x) - \widehat{g}^d_n(x)| + |\widehat{g}^d_n(x) - \widehat{g}^d_n(y)| + |g^d(t_n,y) - \widehat{g}^d_n(y)| \\
     & \le \varepsilon c d^q (1+\|x\|+\|y\|) + c d^q\|x-y\|.
\end{align*}
    To obtain the result, we let $\varepsilon \rightarrow 0^{+} $ and choose the same constants as above.
\end{proof}

Then, following the proof in \citep{gonon23}, we can bound the growth rate of the value function by the following corollary.
\begin{corollary}[Linear and Lipschitz Growth for $V^d_n $] \label{coro:V_growth_bound_new}
    Under Assumption~\ref{ass:dynamic_ass} with $p \ge 1 $ and Assumption~\ref{ass:g_new}, the following linear growth rate properties hold for any $n \in \N \setminus \{ 0 \} $: \( \Lip(V^d_n) ,\; \Growth(V^d_n) \le cd^q   \).
\end{corollary}
\begin{proof}[Proof of Corollary~\ref{coro:V_growth_bound_new}]
The proof for the linear growth rate of $V^d_n $ is the same as Lemma 4.8 in \citep{gonon23} via induction. 
The only thing we need to do to obtain the proof is to show the linear growth rate expressivity of $f^d_n$, which is guaranteed by the Jensen inequality: 
\begin{align*}
    \mathbb{E}\|f^d_n(x,\cdot)\| & \le \big(\mathbb{E}\|f^d_n(x,\cdot)\|^p \big)^{\frac{1}{p}} 
     \le \big(\mathbb{E}\big|\Growth(f^d_n(*,\cdot)) \big|^p \big)^{\frac{1}{p}} (1+\| x \|)   \le c d^q(1+\|x\|),
\end{align*}
    and to combine this with the Lipschitz rate, which has been proved by Proposition~\ref{prop:V_Lip_growth}.
    Finally, we choose the same constants, $c,q$, to complete the proof.
\end{proof}

\begin{proof}[Proof of Lemma~\ref{lemma:prob_new}]
The proof $ \P(U\le M_0) > 0 $ is simple, so we omit it (see, e.g., \citep{Jentzen23}). 
    
    Similar to \citep{gonon23}, our analysis is based on the following important relationship: for any events $A_n,\; n=0,\ldots,N_1$, 
    \( \mathbb{P}\big(\prod_{n=0}^{N_1}A_n \big) \ge  \mathbb{P}\big( \prod_{n=0}^{N_1-1}A_n \big) +\mathbb{P}(A_{N_1}) - 1 \ge \sum_{n=0}^{N_1} \mathbb{P}(A_n)  - N_1 \). 
    The proof is simply obtained using the basic probability formula.
    
    By applying the Markov inequality and Bernoulli inequality, we obtain $\frac{1}{(N_1+2)J_n}\le 1-(\frac{N_1+1}{N_1+2})^{\frac{1}{J_n}},\; n=1,\ldots,N_1 $, \( \mathbb{P}\left(U>(N_1+2)M_0 \right)\le \frac{\mathbb{E}[U]}{(N_1+2)M_0} \le \frac{1}{N_1+2} \), and 
    \[ \mathbb{P}\big( |X^1_n| > (N_1+2))J_nM_n \big) \le \frac{\mathbb{E}|X^1_n|}{(N_1+2)J_n M_n}\le \frac{1}{(N_1+2)J_n}\le 1 - (\frac{N_1+1}{N_1+2})^{\frac{1}{J_n}}, \; n=1,\ldots,N_1. \]
    Then, \( \mathbb{P}\big( U\le (N_1 + 2)M_0 \big)  \ge  \frac{N_1+1}{N_1+2} \), and
    \[ \mathbb{P}\big( \max_{i=1,\ldots,J_n} |X^i_n| \le (N_1+2)J_n M_n \big) = \left[1-\mathbb{P}\left(|X^1_n| > (N_1+2)J_nM_n \right) \right]^{J_n}\ge \frac{N_1+1}{N_1+2} ,\; n=1,\ldots,N_1 .\]
    Hence, \(  \mathbb{P}\big(  U\le (N_1+2) M_0 ,\; B_{N_1} \big) > 0 \), 
    which completes the proof.
\end{proof}

\subsubsection{Detailed proof of the value function approximation}

\begin{proof}[Proof of Lemma~\ref{lem:prob_measure_bound}]
By H\"older's inequality, for any fixed $\bar{p}\in (0,p) $,
    \[
        \M_{\bar{p}}(\bar{\rho}_{n+1}^{K;d}) 
        \le
        (\frac{N}{T})^{\frac{1}{\bar{p}}}d^{-\frac{1}{\bar{p}}} \Big(\sum\limits_{k=0}^{K-1} \big( \mathbb{E} \big\|\Delta W^d_{t^n_k} \big\|^{\frac{2p}{p-\bar{p}}}  \big)^{\frac{p-\bar{p}}{p}} \left(\mathbb{E} \left\|X^d_{t_{n+1}} \right\|^p \right)^{\frac{\bar{p}}{p}} \Big)^{\frac{1}{\bar{p}}} .
    \]
Since $\Delta^i W^d_{t^n_k} \sim N(0,\Delta t), \; \forall \; k\in \K ,\; i=1,\ldots,d $, the Minkowski inequality yields
\[
    \big(\mathbb{E} \|\Delta W^d_{t^n_k} \|^{\frac{2p}{p-\bar{p}}} \big)^{\frac{p-\bar{p}}{p}} 
\le 
\sum_{i=1}^d \big( \mathbb{E} | \Delta^i W^d_{t^n_k} |^{\frac{2p}{p-\bar{p}}} \big)^{\frac{p-\bar{p}}{p}} = C_{p,\bar{p}} \Delta t d 
\]
for all $ k\in \K $,
where $C_{p,\bar{p}} =  \big(\frac{2^{\frac{p}{p-\bar{p}}}}{\sqrt{\pi}}\tilde{\Gamma} (\frac{2p-\bar{p}}{p-\bar{p}}) \big)^{\frac{p-\bar{p}}{p}} $ and $\tilde{\Gamma}(x) $ denotes the Gamma function.
By Assumption~\ref{ass:dynamic_ass}, we have
\( \big(\mathbb{E} \|f^d_n(x,\cdot)\|^p  \big)^{\frac{1}{p}} \le \big(\mathbb{E} |\Growth(f^d_n(*,\cdot)) |^p  \big)^{\frac{1}{p}}(1+\|x\|) \le cd^q(1+\|x\|) \).
Consequently, for all $ n\in \N^{-1} $,
\begin{align*}
    \mathbb{E}\|X^d_{t_{n+1}}\|^p & 
    = \mathbb{E}\big[\mathbb{E} \left[\|f^d_n(x,\cdot)\|^p \right]\big|_{x=X^d_{t_n}} \big] 
     \le 2^{p-1}c^p d^{pq} \left(1+\mathbb{E}\|X^d_{t_n}\|^p \right) \\
    & \le \ldots 
     \le (1 + 2^{p-1}c^p d^{pq})^{n+1}(1+\|x^d_0\|^p )  
     \le (1 + 2^{p-1}c^p )^{N} d^{Npq}(1+\|x^d_0\|^p) .
\end{align*}
Therefore, using $\Delta t = \frac{T}{NN_0} $, we obtain \( \M_{\bar{p}}(\bar{\rho}_{n+1}^{K;d})  \le \widehat{k}_{n+1} d^{\widehat{p}_{n+1}} \) for any $ K\in \bbN_{+} $,
where $\widehat{k}_{n+1}=(C_{p,\bar{p}})^{\frac{1}{\bar{p}}} (1 + 2^{p-1}c^p)^{\frac{N}{p}}(1+c^p)^{\frac{1}{p}} $ and $ \widehat{p}_{n+1}=qN+q $. Note that $\widehat{k}_{n+1} $ and $ \widehat{p}_{n+1} $ are independent of $K$.
\end{proof}

\begin{proof}[Proof of Theorem~\ref{theorem:neural_approx_V_new}]



Applying Lemma~\ref{lem:prob_measure_bound} and Theorem~\ref{theorem:recursive_express_new}, choose $k_1,p_1 \ge 1 $ sufficiently large such that the sequences $(k_n,p_n),n \in \N \setminus \{0\} $ generated from $(k_1,p_1)$ in Theorem~\ref{theorem:recursive_express_new} satisfy $k_n\ge \widehat{k}_n,\; p_n \ge \widehat{p}_n ,\; n\in \N \setminus \{0\} $. This can be achieved by taking the maximum over these requirements, and the resulting choice remains independent of $K$. Hence, there exist constants $c_{n+1}, q_{n+1}, \tau_{n+1}  \ge 1,\; n\in \N^{-1}  $ independent of $d $ such that, for any $\varepsilon >0,K\in \mathbb{N}_{+} $, there exist deep ReLU networks $\widehat{V}^{d,\varepsilon}_{n+1} $ satisfying, for all $ K\in \bbN_{+} $,
\[
    \big\| \widehat{V}^{d,\varepsilon}_{n+1} - V^d_{n+1} \big\|_{2;\bar{\rho}^{K;d}_{n+1}} \le \varepsilon ,\; \text{and} ,\; \Growth(\widehat{V}^{d,\varepsilon}_{n+1}) , \; \size(\widehat{V}^{d,\varepsilon}_{n+1}) \le c_{n+1}d^{q_{n+1}}\varepsilon^{-\tau_{n+1}}.
\]
Consequently, \(  \big\| \widehat{V}^{d,\varepsilon}_{n+1} - V^d_{n+1} \big\|_{2;\tilde{\rho}^{K;d}_{n+1}} \le 
     (\frac{T}{N} d)^{\frac{1}{2}} \varepsilon \).
The proof is completed by choosing the same constants as above while preserving the expressivity estimates.
\end{proof}


\subsubsection{Detailed proof of integrand approximation}

\begin{proof}[Proof of Theorem~\ref{theorem:z_express_approx_new}]

For any $\varepsilon \in (0,1],\; n\in \N^{-1} $, the neural network $\widehat{V}_{n+1} $ satisfies 
\[
    \big\| \widehat{V}^{d,\varepsilon}_{n+1} - V^d_{n+1} \big\|_{2;\tilde{\rho}^{K;d}_{n+1}} \le \varepsilon ,\; \text{and} ,\; \Growth(\widehat{V}^{d,\varepsilon}_{n+1}) , \; \size(\widehat{V}^{d,\varepsilon}_{n+1}) \le c_{n+1}d^{q_{n+1}}\varepsilon^{-\tau_{n+1}} , \; \forall \; K\in \bbN_{+},
\]
and $\widehat{f}^{t_{n+1};d}_{t^n_k} $ satisfies $\widehat{f}^{t_{n+1};d}_{t^n_k} = f^{t_{n+1};d}_{t^n_k}  ,\; k\in \K $
with the properties stated in Assumption~\ref{ass:stronger_ass}. Let $X^{t^n_k,x;d}_{t_{n+1}}=f^{n+1;d}_k(x,\cdot) := f_{t^n_k}^{t_{n+1};d}(x,\cdot) $ and $\theta^{n+1}_k $ be the random parameter of RanNN $\widehat{f}^{n+1;d}_{k}:=\widehat{f}^{t_{n+1};d}_{t^n_k} $. For any $k\in \K $, let the $(\theta^{i,n+1}_k,W^{i;d}_{t^n_k}),i=1,\ldots,J $ i.i.d version of $(\theta^{n+1}_k,W^d_{t^n_k}) $, $\widehat{f}^{i,n+1;d}_k $ be the RanNN w.r.t. $\theta^{i,n+1}_k $ and for $\omega \in \Omega,\; x\in \mathbb{R}^{d}, \; $ and $ J\in \mathbb{N}_{+} $, 
\[\Gamma^{n;d,\varepsilon}_k(x)(\omega) := \frac{1}{J\Delta t}\sum\limits_{i=1}^J \widehat{V}^{d,\varepsilon}_{n+1}\left(\widehat{f}^{i,n+1;d}_{k}(x,\omega) \right)\Delta W^{i;d}_{t^n_k}(\omega) , \] 
and let 
\( \Gamma^{n;d,\varepsilon,K}(t,x)(\omega):= \sum_{k=0}^{K-1} \Gamma^{n;d,\varepsilon}_k(x)(\omega) 1_{t^n_k}(t) \).
Let 
\[ \bar{z}^{n;d,\varepsilon,K}(t,x):= \sum_{k=0}^{K-1} \frac{1}{\Delta t}\mathbb{E}\left[\widehat{V}^{d,\varepsilon}_{n+1}\left(\widehat{f}^{n+1;d}_k(x,\cdot) \right)\Delta W^d_{t^n_k} \right] 1_{t^n_k}(t) ,\]
and \(  I^{d,\varepsilon,K}_n(\omega) := \| \bar{z}^{n;d,\varepsilon,K} -  \Gamma^{n;d,\varepsilon,K}(\omega) \|_{2; \mu^{K;d}_n}   
\) for $ \omega \in \Omega $.
By direct estimation, we obtain
\begin{equation}
    \big\| \widehat{Z}^{*,K,n;d} - \Gamma^{n;d,\varepsilon,K}(\omega) \big\|_{2;\mu^{K;d}_n}
        \le 
         \big\| \widehat{Z}^{*,K,n;d} - \bar{z}^{n;d,\varepsilon,K} \big\|_{2;\mu^{K;d}_n}
        +  I^{d,\varepsilon,K}_n(\omega) \label{eq:total_estimation_neural_construct}. 
\end{equation}
Note that $\widehat{V}^{d,\varepsilon}_{n+1}\big(\widehat{f}^{i,n+1;d}_{k}(x,\cdot) \big)\Delta W^{i;d}_{t^n_k},i=1,\ldots,J$ are i.i.d. together with
\begin{align*}
     \big(\mathbb{E}\|X^d_{t^n_k}\|^2 \big)^{\frac{1}{2}} & = \Big(\mathbb{E}\Big[ \big(\mathbb{E}\big[ \|f_{t_n}^{t^n_k;d}(x,\cdot) \|^{\tilde{p}} \big] \big)^{\frac{2}{\tilde{p}}} \big|_{x=X^d_{t_n}} \Big] \Big)^{\frac{1}{2}} \le \bar{c}d^{\bar{q}} \big(1+ \| X^d_{t_n} \|_{L^2} \big) \\
     & \le \bar{c}d^{\bar{q}}(1+cd^{q}) \big(1+ \| X^d_{t_{n-1}} \|_{L^2} \big) \le  \cdots \le \bar{c}d^{\bar{q}}(1+cd^{q})^{n} (1+\|x^d_0\| ) \le \bar{c}d^{\bar{q}}(1+cd^{q})^{N} \\
\end{align*}
for all $ k\in \K $
by (2).(a) in Assumption~\ref{ass:dynamic_ass} and the similar argument in Theorem~\ref{theorem:neural_approx_V_new} and (2).(a) in Assumption~\ref{ass:stronger_ass}. 
Furthermore, as 
\[
\sum_{k=0}^{K - 1} a_k \le \sum_{k=0}^{K - 1} a_k + 2\sum_{0\le k<m\le K-1} \sqrt{a_k a_m} = \Big( \sum_{k=0}^{K - 1} a_k^{\frac{1}{2}} \Big)^2
\]
for every $a_k\ge 0,\; k\in \K $, which means that \( \big(\sum_{k=0}^{K - 1} a_k \big)^{\frac{1}{2}} \le \sum_{k=0}^{K - 1} a_k^{\frac{1}{2}} \).
Thus, for $I^{d,\varepsilon,K}_n(\omega) $, by the concavity of $x^{\frac{1}{2}} $, the H\"{o}lder inequality, \citep[Lemma~2.1]{Jentzen23} and \eqref{eq:measure_mu_def},
\begin{align*}
   \mathbb{E}[ I^{d,\varepsilon,K}_n] 
     & \le \frac{1}{(\Delta t)^{\frac{1}{2}}} \sum_{k=0}^{K-1}  \Big\|  \frac{1}{J }\mathbb{E}\big\|\widehat{V}^{d,\varepsilon}_{n+1}\big(\widehat{f}^{n+1;d}_{k}(x,*) \big)\Delta W^d_{t^n_k}  - \mathbb{E} \big[\widehat{V}^{d,\varepsilon}_{n+1} \big(\widehat{f}^{n+1;d}_{k}(x,\cdot) \big) \Delta W^d_{t^n_k} \big] \big\|^2 \Big\|_{2;h^{n;d}_k}  \\
     & \le \frac{1}{\sqrt{J\Delta t} } \sum_{k=0}^{K-1} \Big\| \mathbb{E} \big\|\widehat{V}^{d,\varepsilon}_{n+1} \big(\widehat{f}^{n+1;d}_{k}(x,\cdot) \big)\Delta W^d_{t^n_k} \big\|^2  \Big\|_{2;h^{n;d}_k}  \\
     & \le \frac{c_{n+1}d^{q_{n+1}}}{\sqrt{J\Delta t} } \varepsilon^{-\tau_{n+1}} \sum_{k=0}^{K-1} \Big\|  \mathbb{E}\big[\|\Delta W^d_{t^n_k}\|^2+\|\widehat{f}^{n+1;d}_{k}(x,\cdot)\|^2\|\Delta W^d_{t^n_k}\|^2 \big]  \Big\|_{2;h^{n,d}_k}  \\
     & \le \frac{c_{n+1}d^{q_{n+1}}}{\sqrt{J\Delta t} }\varepsilon^{-\tau_{n+1}} \sum_{k=0}^{K-1}\Big[(d\Delta t)^{\frac{1}{2}} +  \big\|  \big(\mathbb{E}\|\widehat{f}^{n+1;d}_{k}(x,\cdot)\|^{\Tilde{p}} \big)^{\frac{2}{\Tilde{p}}} \big(\mathbb{E}\|\Delta W^d_{t^n_k}\|^{\frac{2\Tilde{p}}{\Tilde{p}-2}} \big)^{\frac{\Tilde{p}-2}{\Tilde{p}}} \big\|_{2;h^{n;d}_k}   \Big] \\
     & \le \frac{2C_{\tilde{p}}   \bar{c}c_{n+1}d^{q_{n+1}+ \bar{q}+\frac{1}{2}}}{\sqrt{J} }\varepsilon^{-\tau_{n+1}  } \sum_{k=0}^{K-1}  \big(1 + \big\| X^d_{t^n_k} \big\|_{L^2} \big)
      \le J^{-\frac{1}{2}} \Tilde{c}_{n+1} d^{\Tilde{q}_{n+1}} \varepsilon^{ - \Tilde{\tau}_{n+1}}  K
\end{align*}
with $\Tilde{c}_{n+1}=2C_{\tilde{p}} \bar{c}c_{n+1}(1+ \bar{c}(c+1)^{N}) ,\; \Tilde{q}_{n+1} = q_{n+1}+2\bar{q} + qN +\frac{1}{2},\; \Tilde{\tau}_{n+1}=\tau_{n+1}  $ and $ \; C_{\tilde{p}}=\big(\frac{2^{\frac{\tilde{p}}{\tilde{p}-2}}}{\sqrt{\pi}}\tilde{\Gamma}( \frac{3\tilde{p}-2}{2(\tilde{p}-2)} ) \big)^{\frac{\tilde{p}-2}{2\Tilde{p}}} $.
Thus, let $J^{K,d} := \lceil 9(K+1)^2 (K)^2 (\Tilde{c}_{n+1})^2 d^{2\Tilde{q}_{n+1}} \varepsilon^{-\Tilde{\tau}_{n+1}-2 } \rceil  $, we obtain
\begin{equation*}
    \mathbb{E}[I^{d,\varepsilon,K}_n] \le \frac{1}{3K + 2}\varepsilon,\;  \mathbb{E}\Growth(\widehat{f}^{i,n+1;d}_k(*,\cdot) ) , \mathbb{E}\size(\widehat{f}^{i,n+1;d}_k(*,\cdot) ) \le \bar{c} d^{\bar{q}}, \; i=1,\ldots,J^{K,d} ,\; k\in \K  .
\end{equation*} 
Let us denote
\begin{align*}
    B^1_{K,d} &:= \prod_{k=0}^{K-1} \big\{ \max\limits_{i=1,\ldots,J^{K,d}}\big(\Growth (\widehat{f}^{i,n+1;d}_k(*,\cdot) ) \big)\le (3K+2)J^{K,d} \bar{c} d^{\bar{q}} \big\} , \\
    B^1_{K,d} &:= \prod_{k=0}^{K-1} \big\{ \max\limits_{i=1,\ldots,J^{K,d}}\big(\size (\widehat{f}^{i,n+1;d}_k(*,\cdot) ) \big) \le (3K+2)J^{K,d} \bar{c} d^{\bar{q}} \big\}  \\
    B^3_{K,d} &:= \prod_{k=0}^{K-1} \big\{ \max\limits_{i=1,\ldots,J^{K,d}}\| \Delta W^{i;d}_{t^K}\| \le (3K+2)J^{K,d} (d\Delta t)^{\frac{1}{2}} .
\end{align*}
By Lemma~\ref{lemma:prob_new}, we have
\(
    \P( I^{d,\varepsilon,K}_n \le  \varepsilon , B^1_{K,d} , B^2_{K,d} , B^3_{K,d} ) > 0  
\).
Then, there exists an $\omega_0 \in \Omega $, such that \( I^{d,\varepsilon,K}_n(\omega_0) \le \varepsilon \),
and for all $ x,y\in \mathbb{R}^d, \; i=1,\ldots,J^{K,d} ,\; k\in \K $, 
\[
    \Growth(\widehat{f}^{i,n+1;d}_{k}(*,\omega_0)) ,   \size(\widehat{f}^{i,n+1;d}_{k}(*,\omega_0)) \le (3K+2) J^{K,d} \bar{c} d^{\bar{q}} 
\]
and \( \big\|\Delta W^{i;d}_{t^n_k}(\omega_0) \big\|  \le (3K+2)J^{K,d}(d\Delta t)^{\frac{1}{2}}  \).
From \citep[Propositions 2.2, 2.3]{opschoor20}, we can realize deep ReLU networks for all $k\in \K $ as follows:
\[\gamma^{n;d,\varepsilon,K}_k(x) = \Gamma^{n;d,\varepsilon}_k(x)(\omega_0) = \frac{1}{J^{K,d} \Delta t}\sum_{i=1}^{J^{K,d}} \widehat{V}^{d,\varepsilon}_{n+1}\big(\widehat{f}_k^{i,n+1;d}(x,\omega_0) \big)\Delta W^{i;d}_{t^n_k}(\omega_0)  \]
with the size, growth rate, and Lipschitz bound determined as
\begin{align*}
    \size(\gamma^{n;d,\varepsilon,K}_k) & \le \sum_{i=1}^{J^{K,d}} \big(\size(\widehat{V}^{d,\varepsilon}_{n+1})  + \size(\widehat{f}_k^{i,n+1;d}(\cdot,\omega_0)) + d \big) \\
    & \le \sum\limits_{i=1}^{J^{K,d}} \big( c_{n+1} d^{q_{n+1}} \varepsilon^{-\tau_{n+1}} + (3K+2)J^{K,d} \bar{c} d^{\bar{q}} + d \big) 
     \le c_{n,1} d^{q_{n,1}} \varepsilon^{-\tau_{n,1}} (K)^9,
\end{align*}
where $c_{n,1}=9^2\cdot4 ((\Tilde{c}_{n+1})^2+1)^2 (c_{n+1} + \bar{c}+1) ,\; q_{n,1}=4\Tilde{q}_{n+1}+q_{n+1}+ \bar{q} +1 ,\; \tau_{n,1}=4(\Tilde{\tau}_{n+1}+1) + \tau_{n+1}   $,
\begin{align*}
    \|\gamma^{n;d,\varepsilon,K}_k(x) \| & \le \frac{1}{ J^{K,d} \Delta t} \sum_{i=1}^{J^{K,d}} \big|\widehat{V}^{d,\varepsilon}_{n+1}(\widehat{f}_k^{i,n+1;d}(x,\omega_0) ) \big| \big\|\Delta W^{i;d}_{t^n_k}(\omega_0) \big\| \\
    & \le \frac{c_{n+1}d^{q_{n+1}+\frac{1}{2}}}{(\Delta t)^{\frac{1}{2}}} \varepsilon^{-\tau_{n+1}}  (3K+2) \sum\limits_{i=1}^{J^{K,d} } \left(1+\|\widehat{f}_k^{i,n+1;d}(x,\omega_0)\| \right)  \\
    & \le \frac{ \bar{c}c_{n+1}d^{q_{n+1}+\frac{1}{2}+ \bar{q}}}{(\Delta t)^{\frac{1}{2}}} \varepsilon^{-\tau_{n+1}} (3K+2)^2 (J^{K,d})^2 (1+\|x\|) 
     \le c_{n,2}d^{q_{n,2}} \varepsilon^{-\tau_{n,2}} (K)^{\frac{21}{2}} (1+\|x\| ) 
\end{align*}
for all $ x\in \R^d $, 
where $c_{n,2}=18^2\cdot 25 ( ((\Tilde{c}_{n+1})^2 + 1  )^2 (\frac{N}{T})^{\frac{1}{2}} \bar{c}c_{n+1} ,\; q_{n,2}=4\Tilde{q}_{n+1}+ q_{n+1}+ \bar{q}+\frac{1}{2},\; $ and $ \tau_{n,2}=4(\Tilde{\tau}_{n+1} + 1)+\tau_{n+1}  $. 
For the first part of the target estimation, by Jensen inequality,
\begin{align*}
        \big\| \widehat{Z}^{*,K,n;d} - \bar{z}^{n;d,\varepsilon,K} \big\|_{2;\mu^{K;d}_n} 
        &\le \frac{1}{(\Delta t)^{\frac{1}{2}}}  
            \Big(\sum_{k=0}^{K-1} \mathbb{E} \big[\mathbb{E}\big[ \big\| \big(V^d_{n+1}(X^d_{t_{n+1}}) - \widehat{V}^{d,\varepsilon}_{n+1}(X^d_{t_{n+1}}) \big) \Delta W^d_{t^n_k} \big\|^2  \mid X^d_{t^n_k} \big] \big]   \Big)^{\frac{1}{2}} \\
        & \le \frac{1}{(\Delta t)^{\frac{1}{2}}}  \big\| V^d_{n+1}(x) - \widehat{V}^{d,\varepsilon}_{n+1}(x) \big\|_{2;\tilde{\rho}^{K;d}_{n+1}}  \le (\frac{N}{T})^{\frac{1}{2}} (K)^{\frac{1}{2}} \varepsilon. 
\end{align*}
Let $\widehat{z}^{d,\varepsilon,K}_n(t,x):=\Gamma^{n;d,\varepsilon,K}(t,x)(\omega_0) = \sum_{k=0}^{K-1} \gamma^{n;d,\varepsilon,K}_k  1_{t^n_k}(t) $. By plugging these results into the target estimation~\eqref{eq:total_estimation_neural_construct}, we obtain
\( \big\|  \widehat{z}^{d,\varepsilon,K}_n -  \widehat{Z}^{*,K,n;d}  \big\|_{2;\mu^{K;d}_n} \le [ (\frac{N}{T})^{\frac{1}{2}} (K)^{\frac{1}{2}}  + 1 ] \varepsilon \).
Then, for any $\bar{\varepsilon}\in (0,1] $, we take $\varepsilon =\bar{\varepsilon} [ (\frac{N}{T})^{\frac{1}{2}} (K)^{\frac{1}{2}} + 1 ]^{-1} $ and again choose the constants $\widehat{c}_n,\widehat{q}_n,\widehat{\tau}_n, $ and $ \widehat{m}_n \ge 1 $ independent of $k,d,\bar{\varepsilon}, $ and $ K $, such that (denote $ \widehat{z}^{d,\bar{\varepsilon},K}_n = \widehat{z}^{d,\varepsilon,K}_n $, $ \gamma^{n;d,\bar{\varepsilon},K}_k = \gamma^{n;d,\varepsilon,K}_k $ )
\[
      \big\| \widehat{z}^{d,\bar{\varepsilon},K}_n -   \widehat{Z}^{*,K,n;d} \big\|_{2;\mu^{K;d}_n } \le \bar{\varepsilon}  ,\; \; \text{and} \; \; \Growth( \gamma^{n;d,\bar{\varepsilon},K}_k) ,\; \size(\gamma^{n;d,\bar{\varepsilon},K}_k) \le  \widehat{c}_n d^{\widehat{q}_n}  \bar{\varepsilon}^{-\widehat{\tau}_n} (K)^{\widehat{m}_n} ,\; \forall \; k\in \K   ,   
\]
which completes the proof.
\end{proof}

\begin{proof}[Proof of Lemma~\ref{lem:time_indicator_realization}]

We first justify the deep ReLU network realization of $h^n_k $. The case of $ h^{n;K}_0 $ is immediate, since $h^{n;K}_0(t)=\frac{1}{t^n_1 - t^n_0 }(t^n_1 - t )^{+} $. For $ k\in \K \setminus \{0\} =: \K^{-1} $, 
define
\( g^{n;K}_k(t)  = \max \big( \min (t, t^n_{k+1}), t^n_{k-1} \big) \). Then
    \[
         h^{n;K}_k(t) = \frac{1}{t^n_{k+1}- t^n_k } [t^n_{k+1} - g^{n;K}_k(t) ] - \frac{t^n_{k+1} - t^n_{k-1} }{(t^n_{k+1}-t^n_k )(t^n_k - t^n_{k-1} )} \big(t^n_k -g^{n;K}_k(t) \big)^{+} .
    \]
Using \( \max(x,y) = (x-y)^{+} + y^{+} - (-y)^{+} \),
we obtain
\[
    \min(t,t^n_{k+1}) = -\max( -t, -t^n_{k+1} ) = - [ (t^n_{k+1} -t)^{+} + (-t^n_{k+1})^{+} - (t^n_{k+1})^{+} ] = - (t^n_{k+1}-t )^{+} + t^n_{k+1}   .
\]
Hence,
\( g^{n;K}_k(t) = \max[ - (t^n_{k+1}-t )^{+} + t^n_{k+1}, t^n_{k-1} ] = \big(- (t^n_{k+1}-t )^{+} + t^n_{k+1} - t^n_{k-1} \big)^{+} + t^n_{k-1} \)
and is therefore a deep ReLU network with $\size(g^{n;K}_k) = 6 $. Since $t^n_{k+1} - g^{n;K}_k(t) \ge 0 $, we have $t^n_{k+1} - g^{n;K}_k(t) = (t^n_{k+1} - g^{n;K}_k(t) )^{+} $. Combining this with \citep[Propositions 2.2, 2.3]{opschoor20}, we obtain the deep ReLU realization
\[
    h^{n;K}_k(t) = A_2 \; \sigma(A_1 g^{n;K}_k(t) + b_1 )
\]
where \( A_1 = - (1,1)^{\top} ,\; b_1 = (t^n_{k+1} , t^n_k)^{\top} \) and $ A_2 = ( \frac{1}{t^n_{k+1}-t^n_k} , -\frac{t^n_{k+1}-t^n_{k-1}}{(t^n_{k+1}-t^n_k)(t^n_{k} - t^n_{k-1})} ) $,
and where $\sigma$ denotes the component-wise ReLU activation function. Consequently, $h^n_k $ is a deep ReLU network with $\size(h^{n;K}_k)=\size(g^{n;K}_k) + 6 = 12 $, for all $ k \in \K^{-1} $.
\end{proof}


\begin{proof}[Proof of Theorem~\ref{theorem:joint nerual network realization_new}]

By Theorem~\ref{theorem:z_express_approx_new}, there exist constants $\widehat{c}_n,\widehat{q}_n,\widehat{\tau}_n,\widehat{m}_n \ge 1 $ independent of $d,K$, such that for any given $\bar{\varepsilon}\in (0,\frac{1}{2}) $ (to be specified later) and $K \in \mathbb{N}_{+} $ , there exists a family of deep ReLU networks $(\gamma^{n;d,\bar{\varepsilon},K}_k)_{k=0}^{K-1} $ and a joint function $\widehat{z}^{d,\bar{\varepsilon},K}_n $ such that $ \widehat{z}^{d,\bar{\varepsilon},K}_n(t^n_k,\cdot) \equiv \gamma^{n;d,\bar{\varepsilon},K}_k $ and
\[
          \big\| \widehat{z}^{d,\bar{\varepsilon},K}_n -   \widehat{Z}^{*,K,n;d} \big\|_{2;\mu^{K;d}_n } \le \bar{\varepsilon}  ,\; \;  \text{and} \; \; \Growth( \gamma^{n;d,\bar{\varepsilon},K}_k) ,\; \size(\gamma^{n;d,\bar{\varepsilon},K}_k) \le  \widehat{c}_n d^{\widehat{q}_n}  \bar{\varepsilon}^{-\widehat{\tau}_n} (K)^{\widehat{m}_n} ,\; \forall \; k\in \K .
\]

Let $ (h^{n;K}_k)_{k=0}^{K-1} $ be the deep ReLU networks constructed in Lemma~\ref{lem:time_indicator_realization}. 
Since $h^{n;K}_k(t^n_p)= \delta_{kp} $ (Kronecker symbol), the joint function $\widehat{z}^{d,\bar{\varepsilon},K}_n $ admits the representation
\begin{equation}
    \widehat{z}^{d,\bar{\varepsilon},K}_n(t,x) = \sum_{k=0}^{K-1} h^{n;K}_k(t) \gamma^{n;d,\bar{\varepsilon},K}_k(x) =: \Gamma^{d,\bar{\varepsilon},K}_n(t,x) , \; 
    \; \mu^{K;d}_n  \text{-a.e.}  \label{eq:z_realization_mu_new}
\end{equation}
To verify \eqref{eq:z_realization_mu_new}, note that
\begin{align*}
     \big\| \widehat{z}^{d,\bar{\varepsilon},K}_n -  \Gamma^{d,\bar{\varepsilon},K}_n\big\|^2_{2;\mu^{K;d}_n} & = \E \Big[\sum_{p=0}^{K-1} \Big( \gamma^{n;d,\bar{\varepsilon},K}_p(x) - \sum_{k=0}^{K-1} \delta_{kp} \gamma^{n;d,\bar{\varepsilon},K}_k(x) \Big)^2 \Delta t \Big] \\
     & = \mathbb{E}\Big[\sum_{p=0}^{K-1}\big(\gamma^{n;d,\bar{\varepsilon},K}_p(x) -  \gamma^{n;d,\bar{\varepsilon},K}_p(x) \big)^2 \Delta t\Big] = 0 .
\end{align*}


Then, by \citep[Proposition 4.1]{opschoor20} and \citep[Lemma 4.1]{gonon23}, there exists a constant $c \ge 1 $, and for the above given $\bar{\varepsilon} $ and any $M \ge 1 $ (will be chosen later), there exists a deep ReLU network $ \calN^{\bar{\varepsilon},M} : \mathbb{R}^{2} \rightarrow \mathbb{R} $, such that
\[
      \sup_{t,y \in [-M,M ]}| \calN^{\bar{\varepsilon},M}(t,y) - ty|   \le \bar{\varepsilon} ,\quad \size(\calN^{\bar{\varepsilon},M}) \le c (\log(\bar{\varepsilon}^{-1}) + \log(M ) + 1 ) \le c (\bar{\varepsilon}^{-1} + M -1 )   ,
\]
and satisfies for all $t,t^{'},y,y^{'} \in \mathbb{R} $
\begin{align}
    | \calN^{\bar{\varepsilon},M}(t,y) - \calN^{\bar{\varepsilon},M}(t^{'},y^{'}) | & \le M c ( |t-t^{'}| + | y-y^{'} | ) ,\; \text{and} \label{eq:Lip_prod_new} \\
    \calN^{\bar{\varepsilon},M}(t,0) &= \calN^{\bar{\varepsilon},M}(0,y)=0. \label{eq:zero_prod_new}
\end{align}

Let 
\( \tilde{z}^{d,\bar{\varepsilon},K,M}_n(t,x) = \sum\limits_{k=0}^{K-1} \mathbf{n}^{d,\bar{\varepsilon},M} \big( h^{n;K}_k(t) , \gamma^{n;d,\bar{\varepsilon},K}_k(x) \big) \),
where for all $ (t_1,y_1) \in \R^{1+d} $, 
\[
    \mathbf{n}^{d,\bar{\varepsilon},M}( t_1 , y_1 ) = \Big( \calN^{\bar{\varepsilon},M}\big(t_1, (y_1)_1\big) , \ldots ,  \calN^{\bar{\varepsilon},M}\big(t_1, (y_1)_d\big)\Big)^{\top}  .
\]
As $0\le h^{n;K}_k(t) \le 1 \le M,\; \forall \; t\in [t_n,t_{n+1}] $, then for any $x\in \mathbb{R}^d $ that satisfies $\|x\| \le M_0 := M(\widehat{c}_n)^{-1}d^{-\widehat{q}_n} \bar{\varepsilon}^{\widehat{\tau}_n} (K)^{-\widehat{m}_n} - 1 $ (note $M $ here should be $> \widehat{c}_n d^{\widehat{q}_n} \bar{\varepsilon}^{-\widehat{\tau}_n} (K)^{\widehat{m}_n}$) , we have 
\[ \| \gamma^{n;d,\bar{\varepsilon},K}_k(x)\|\le \widehat{c}_nD^{\widehat{q}_n}\bar{\varepsilon}^{-\widehat{\tau}_n} (K)^{\widehat{m}_n} (1+\|x\|)\le M, \; \forall \; k\in \K ,\]
which immediately implies for all $ i=1,\ldots,d $, $ k \in \K $, \(  \big|\big( \gamma^{n;d,\bar{\varepsilon},K}_k(x) \big)_i | \le M \).
Let $B(0,M_0):= \left\{ x\in \mathbb{R}^d : \|x\| \le M_0 \right\} $. Then, for all $ i=1,\ldots,d $, $ k \in \K $, $t\in [t_n,t_{n+1}] $, $x\in B(0,M_0) $, the following 
\[ \Big| \calN^{\bar{\varepsilon},M} \Big(h^{n;K}_k(t), \big(\gamma^{n;d,\bar{\varepsilon},K}_k(x) \big)_i \Big) - h^{n;K}_k(t) \big(\gamma^{n;d,\bar{\varepsilon},K}_k(x) \big)_i \Big| \le \bar{\varepsilon} \]
holds, and thus \( \big\| \mathbf{n}^{d,\bar{\varepsilon},M} \big( h^{n;K}_k(t) , \gamma^{n;d,\bar{\varepsilon},K}_k(x) \big) - h^{n;K}_k(t) \gamma^{n;d,\bar{\varepsilon},K}_k(x) \big\| \le d^{\frac{1}{2}} \bar{\varepsilon} \) holds for all \( t\in [t_n,t_{n+1}] ,\; x\in B(0,M_0) \).
Immediately, for all \( t\in [t_n,t_{n+1}] ,\; x\in B(0,M_0) \),
\[ \big\| \widehat{z}^{d,\bar{\varepsilon},K}_n(t,x) - \tilde{z}^{d,\bar{\varepsilon},K,M}_n(t,x) \big\| \le \sum_{k=0}^{K-1} \big\| \mathbf{n}^{d,\bar{\varepsilon},M} \big( h^{n;K}_k(t) , \gamma^{n;d,\bar{\varepsilon},K}_k(x) \big) - h^{n;K}_k(t) \gamma^{n;d,\bar{\varepsilon},K}_k(x) \big\| \le K d^{\frac{1}{2}} \bar{\varepsilon}  . \]
By \eqref{eq:Lip_prod_new} and \eqref{eq:zero_prod_new}, for all $t\in \mathbb{R} $ and $ x\in \mathbb{R}^d $,
\begin{align*}
    \Big| \calN^{\bar{\varepsilon},M} \Big(h^{n;K}_k(t), \big(\gamma^{n;d,\bar{\varepsilon},K}_k(x) \big)_i \Big) \Big| & = 
    \Big| \calN^{\bar{\varepsilon},M} \Big(h^{n;K}_k(t), \big(\gamma^{n;d,\bar{\varepsilon},K}_k(x) \big)_i \Big) - \calN^{\bar{\varepsilon},M} (0,0) \Big|  \\
     & \le M c \big(|h^{n;K}_k(t)| + \big|(\gamma^{n;d,\bar{\varepsilon},K}_k(x))_i \big| \big) \\
    & \le 2M c\widehat{c}_nD^{\widehat{q}_n}\bar{\varepsilon}^{-\widehat{\tau}_n} (K)^{\widehat{m}_n}(1+\|x\|) ,
\end{align*}
from which we immediately deduce the growth rate for $\tilde{z}^{d,\bar{\varepsilon},K,M}_n $ as: for all $ t\in [t_n,t_{n+1}],\; x\in \R^d $,
\begin{align*}
    \big\| \tilde{z}^{d,\bar{\varepsilon},K,M}_n(t,x)  \big\| & \le \sum\limits_{k=0}^{K-1} \big\|\mathbf{n}^{d,\bar{\varepsilon},M} ( h^{n;K}_k(t) , \gamma^{n;d,\bar{\varepsilon},K}_k(x) )  \big\|  \\
    & \le \sum_{k=0}^{K-1} \Big(\sum\limits_{i=1}^d \big| \calN^{\bar{\varepsilon},M} \big(h^{n;K}_k(t), (\gamma^{n;d,\bar{\varepsilon},K}_k(x) )_i \big) \big|^2 \Big)^{\frac{1}{2}}  \\
    & = 2 M c\widehat{c}_n  d^{\widehat{q}_n + \frac{1}{2} }\bar{\varepsilon}^{-\widehat{\tau}_n} (K)^{\widehat{m}_n+1}(1+\|x\|)  ,
\end{align*}
and the growth bound for $ \widehat{z}^{d,\bar{\varepsilon},K}_n  $ as: for all $ t\in [t_n,t_{n+1}],\; x\in \R^d $,
\begin{align*}
    \| \widehat{z}^{d,\bar{\varepsilon},K}_n(t,x)\| & \le \sum_{k=0}^{K-1} \big|h^{n;K}_k(t) \big|  \| \gamma^{n;d,\bar{\varepsilon},K}_k(x) \| 
     \le  \widehat{c}_n  d^{\widehat{q}_n}\bar{\varepsilon}^{-\widehat{\tau}_n} (K)^{\widehat{m}_n+1}(1+\|x\|)  .
\end{align*}
Then, by the H\"older inequality, the following estimation holds:
\begin{align*}
     & \big\| \widehat{z}^{d,\bar{\varepsilon},K}_n - \tilde{z}^{d,\bar{\varepsilon},K,M}_n \big\|^2_{2;\mu^{K;d}_n}  
       =  \int_{\mathbb{R}\times B(0,M_0) } \big\| \widehat{z}^{d,\bar{\varepsilon},K}_n(t,x) - \tilde{z}^{d,\bar{\varepsilon},K,M}_n(t,x) \big\|^2 \mu^{K;d}_n(dtdx) \\
     & \quad \quad \quad \quad + \int_{\mathbb{R}\times B^c(0,M_0) } \big\| \widehat{z}^{d,\bar{\varepsilon},K}_n(t,x) - \tilde{z}^{d,\bar{\varepsilon},K,M}_n(t,x) \big\|^2 \mu^{K;d}_n(dtdx) \\
    & \le (K)^2 d (\bar{\varepsilon})^2 \mu^{K;d}_n(\mathbb{R}^{1+d} ) + \int_{\mathbb{R}\times B^c(0,M_0) } \big(\| \widehat{z}^{d,\bar{\varepsilon},K}_n(t,x)\| + \|\tilde{z}^{d,\bar{\varepsilon},K,M}_n(t,x)\| \big)^2 \mu^{K;d}_n(dt dx) \\
    & \le 8M^2(c\widehat{c}_n)^2 d^{2\widehat{q}_n + 1 }\bar{\varepsilon}^{-2\widehat{\tau}_n} (K)^{2(\widehat{m}_n+1)}     \Big(\int_{\mathbb{R}^{1+d} } (1 + \|x\| )^4 \mu^{K;d}_n(dt dx) \Big)^{\frac{1}{2}} \big(\mu^{K;d}_n (\mathbb{R}\times B^c(0,M_0) ) \big)^{\frac{1}{2}} \\
    & \quad \quad \quad \quad + (K)^2 d (\bar{\varepsilon})^2 \frac{T}{N}.
\end{align*}
By Assumption~\ref{ass:dynamic_ass}, Assumption~\ref{ass:stronger_ass}, and the similar argument in the proof of Theorem~\ref{theorem:z_express_approx_new},
\begin{align*}
   \big\| X^d_{t^n_k} \big\|_{L^{\tilde{p}}}  &  = \Big(\mathbb{E}\Big[ \mathbb{E}\big[\big\|f_{t_n}^{t^n_k;d}(x,\cdot) \big\|^{\tilde{p}} \big] \mid_{x=X_{t_n}} \Big] \Big)^{\frac{1}{\tilde{p}}} 
   \le  \bar{c}d^{\bar{q}} \big[ 1 +  \big\| X^d_{t_n} \big\|_{L^{\tilde{p}}} \big] 
     \le \bar{c}d^{\bar{q}}( 1 + cd^q)^{N} ,\; k\in \K .
\end{align*}
Using the monotonicity of $L^p $-norm, we obtain for all $ k\in \K $, \( \big\| X^d_{t^n_k} \big\|_{L^4} \le \big\| X^d_{t^n_k} \big\|_{L^{\tilde{p}}} \le \bar{c}d^{\bar{q}}( 1 + cd^q)^{N} \).
Then,
\begin{align*}
     \Big(\int_{\mathbb{R}^{1+d} } (1 + \|x\| )^4 \mu^{K;d}_n(dt dx) \Big)^{\frac{1}{2}}  
     & \le \big(\mu^{K;d}_n(\mathbb{R}^{1+d}) \big)^{\frac{1}{4}} +  \Big(\sum_{k=0}^{K-1} \big\| X^d_{t^n_k} \big\|^4_{L^4} \Delta t \Big)^{\frac{1}{4}} \\
    & \le (\frac{T}{N})^{\frac{1}{4}}\left[1+ \bar{c}d^{\bar{q}}( 1 + cd^q)^{N} \right] 
    \le \widehat{c}_{n,1} d^{\widehat{q}_{n,1}},
\end{align*}
where $\widehat{c}_{n,1}=(\frac{T}{N} )^{\frac{1}{4}}[1+\bar{c}(1+c)^{N} ] $ and $ \; \widehat{q}_{n,1}=Nq + \bar{q} $. For $ \mu^{K;d}_n (\mathbb{R}\times B^c(0,M_0) ) $, we require $M\ge 2\widehat{c}_n d^{\widehat{q}_n}\bar{\varepsilon}^{-\widehat{\tau}_n} (K)^{\widehat{m}_n} $, which leads to $M_0=M(\widehat{c}_n d^{\widehat{q}_n} \bar{\varepsilon}^{-\widehat{\tau}_n} (K)^{\widehat{m}_n})^{-1} -1 \ge \frac{1}{2}M \left(\widehat{c}_n d^{\widehat{q}_n} \bar{\varepsilon}^{-\widehat{\tau}_n}(K)^{\widehat{m}_n} \right)^{-1} $.
Then, by the Markov inequality, we have
\begin{align*}
    \big(\mu^{K;d}_n (\mathbb{R}\times B^c(0,M_0) )  \big)^{\frac{1}{2}} 
    & =  \Big( \sum\limits_{k=0}^{K-1} \mathbb{P} ( \|X^d_{t^n_k}\| > M_0
    ) \Delta t  \Big)^{\frac{1}{2}} 
    \le \Big( \sum\limits_{k=0}^{K-1} \frac{1}{(M_0)^{\Tilde{p}}} \big\| X^d_{t^n_k} \big\|^{\tilde{p}}_{L^{\tilde{p}}} \Delta t \Big)^{\frac{1}{2}} \\
    & \le \frac{1}{(M_0)^{\frac{\Tilde{p}}{2}} }(\frac{T}{N})^{\frac{1}{2}}( 1 + c d^q)^{\frac{\Tilde{p}}{2}(N+1)} 
     \le \frac{1}{M^{\frac{\Tilde{p}}{2}}} \widehat{c}_{n,2}d^{\widehat{q}_{n,2}}\bar{\varepsilon}^{-\widehat{\tau}_{n,2}}  (K)^{\widehat{m}_{n,2}} ,
\end{align*}
where $\widehat{c}_{n,2}=  (\frac{T}{N})^{\frac{1}{2}}(2\widehat{c}_n)^{\frac{\Tilde{p}}{2}}(1+c)^{\frac{\Tilde{p}}{2}(N+1)} ,\; \widehat{q}_{n,2}= \frac{\Tilde{p}}{2}\widehat{q}_n +\frac{\Tilde{p}}{2}(N+1)q,\; \widehat{\tau}_{n,2}=\frac{\Tilde{p}}{2}\widehat{\tau}_n,\; $ and $ \widehat{m}_{n,2}=\frac{\Tilde{p}}{2}\widehat{m}_n  $. 
Then, by combining the results, we obtain
\begin{align*}
    \big\| \widehat{z}^{d,\bar{\varepsilon},K}_n - \tilde{z}^{d,\bar{\varepsilon},K,M}_n \big\|^2_{2;\mu^{K;d}_n}  
    & \le (K)^2 d (\bar{\varepsilon})^2 \frac{T}{N} + \frac{1}{M^{\frac{\tilde{p}-4}{2}}} \widehat{c}_{n,3}d^{\widehat{q}_{n,3}} \bar{\varepsilon}^{-\widehat{\tau}_{n,3}} (K)^{\widehat{m}_{n,3}},
\end{align*}
where $\widehat{c}_{n,3}=8(c\widehat{c}_n\widehat{c}_{n,1})^2 \widehat{c}_{n,2},\; \widehat{q}_{n,3}=2(\widehat{q}_{n} + \widehat{q}_{n,1} )+ \widehat{q}_{n,2} +1,\; \widehat{\tau}_{n,3}=\widehat{\tau}_{n,2}+2\widehat{\tau}_{n}$ and $ \widehat{m}_{n,3} = \widehat{m}_{n,2}+2(\widehat{m}_n + 1) $. 
We let $\widehat{c}_{n,4}=\max(2\widehat{c}_n , (\widehat{c}_{n,3})^{\frac{2}{\tilde{p}-4}} ),\; \widehat{q}_{n,4}=\max(\widehat{q}_n , \frac{2}{\tilde{p}-4}\widehat{q}_{n,3} ),\; \widehat{\tau}_{n,4}= \frac{4}{\tilde{p}-4} +\frac{2}{\tilde{p}-4}\widehat{\tau}_{n,3}  ,\; $ and $ \widehat{m}_{n,4}= \max(\widehat{m_n} , \frac{2}{\tilde{p}-4}\widehat{m}_{n,3} ) $ and choose $ M_{d,\bar{\varepsilon},K} \ge \widehat{c}_{n,4} d^{\widehat{q}_{n,4}} \bar{\varepsilon}^{-\widehat{\tau}_{n,4}} (K)^{\widehat{m}_{n,4} } $, and then have ($ \tilde{z}^{d,\bar{\varepsilon},K }_n := \tilde{z}^{d,\bar{\varepsilon},K, M_{d,\bar{\varepsilon},K} }_n $), \(  \big\| \widehat{z}^{d,\bar{\varepsilon},K}_n - \tilde{z}^{d,\bar{\varepsilon},K }_n \big\|^2_{2;\mu^{K;d}_n}  \le \bar{\varepsilon}^2 (\frac{T}{N} + 1 ) d (K)^{2} \).
Then, 
\[
     \big\| \widehat{Z}^{*,K,n;d}_n - \tilde{z}^{d,\bar{\varepsilon},K }_n \big\|^2_{2;\mu^{K;d}_n} \le  \big\|  \widehat{Z}^{*,K,n;d}_n - \widehat{z}^{d,\bar{\varepsilon},K}_n   \big\|^2_{2;\mu^{K;d}_n} +  \big\| \widehat{z}^{d,\bar{\varepsilon},K}_n - \tilde{z}^{d,\bar{\varepsilon},K }_n \big\|^2_{2;\mu^{K;d}_n} \le  \bar{\varepsilon} [1 + (\frac{T}{N} + 1 )  d^{\frac{1}{2}} K ]  .
\]
Note that $(\gamma^{n;d,\bar{\varepsilon},K}_k )_{i},i=1,\ldots,d $ are also deep ReLU networks with $\size((\gamma^{n;d,\bar{\varepsilon},K}_k)_i) \le \size(\gamma^{n;d,\bar{\varepsilon},K}_k)$, $i=1,\ldots,d $.
Thus, for any given $\varepsilon \in (0,1) $, we choose $\bar{\varepsilon}= \varepsilon [1 + (\frac{T}{N} + 1 ) d^{\frac{1}{2}} K ]^{-1}  $ (easy to verify $\bar{\varepsilon}<\frac{1}{2} $). By denoting $ \tilde{z}^{d, {\varepsilon},K }_n := \tilde{z}^{d,\bar{\varepsilon},K }_n $, then
\[
    \big\| \widehat{Z}^{*,K,n;d}_n - \tilde{z}^{d, {\varepsilon},K }_n \big\|_{2;\mu^{K;d}_n} \le \varepsilon
\]
together with (after applying \citep[Propositions 2.2, 2.3]{opschoor20})
\begin{align*}
    \size( \tilde{z}^{d, {\varepsilon},K }_n ) & \le 2\sum_{k=0}^{K-1}\sum_{i=1}^d \big[ \size( \calN^{\bar{\varepsilon},M_{d,\bar{\varepsilon},K}} ) + \big( \size(h^{n;K}_k) + \size((\gamma^{n;d,\bar{\varepsilon},K}_k)_i) \big)  \big] \\
    & \le 2\sum_{k=0}^{K-1}\sum_{i=1}^d \big[  c(\bar{\varepsilon}^{-1} +  M_{d,\bar{\varepsilon},K}  -1) + ( 12 + \widehat{c}_n d^{\widehat{q}_n} \bar{\varepsilon}^{-\widehat{\tau}_n} (K)^{\widehat{m}_n} )  \big] 
    \le \widehat{c}_{n,5}d^{\widehat{q}_{n,5}}\varepsilon^{-\widehat{\tau}_{n,5}} (K)^{\widehat{m}_{n,5}} ,
\end{align*}
where $\widehat{c}_{n,5}= 2 \left[c \left(2+\frac{T}{N} + \widehat{c}_{n,4}(2+\frac{T}{N})^{\widehat{\tau}_{n,4}} -1 \right) +  12 + \widehat{c}_n(2+\frac{T}{N})^{\widehat{\tau}_n}  \right] ,\; \widehat{q}_{n,5}= \frac{1}{2}(1+\widehat{\tau}_{n,4} + \widehat{\tau}_n) + \widehat{q}_{n,4} + \widehat{q}_n +1   ,\; \widehat{\tau}_{n,5}= 1+\widehat{\tau}_{n,4}+\widehat{\tau}_n ,\; $ and $ \widehat{m}_{n,5}=1+ \widehat{m}_{n,4} + \widehat{\tau}_{n,4} + \widehat{m}_n + \widehat{\tau}_n  $, and for any $x \in \mathbb{R}^d $,
\begin{align*}
    \| \tilde{z}^{d, {\varepsilon},K }_n(t,x)  \| & \le 2 M_{d,\bar{\varepsilon},K} c\widehat{c}_n  d^{\widehat{q}_n + \frac{1}{2} } \bar{\varepsilon}^{-\widehat{\tau}_n}(K)^{\widehat{m}_n+1} (1+\|x\|)  \le \widehat{c}_{n,6} d^{\widehat{q}_{n,6}} \varepsilon^{-\widehat{\tau}_{n,6}}(K)^{\widehat{m}_{n,6}}  (1+\|x\|), 
\end{align*}
where $\widehat{c}_{n,6}=2c\widehat{c}_n(2+\frac{T}{N})^{ \widehat{\tau}_n + \widehat{\tau}_{n,4}}\widehat{c}_{n,4} ,\; \widehat{q}_{n,6}=\widehat{q}_n + \frac{1}{2}(1 + \widehat{\tau}_{n,4}+\widehat{\tau}_n) + \widehat{q}_{n,4}  ,\; \widehat{\tau}_{n,6}= \widehat{\tau}_{n,4}+\widehat{\tau}_n,\; $ and $ \widehat{m}_{n,6}=\widehat{m}_n + 1  + \widehat{m}_{n,4} + \widehat{\tau}_{n,4}+\widehat{\tau}_n $. By choosing the same constants, we complete the proof. 
\end{proof}

\subsubsection{Detailed proof of the expressivity of DeepMartingale}

\begin{proof}[Proof of Theorem~\ref{thm:express_deep_mtg}]
We directly combine Theorems~\ref{theorem:express_N_0_new} and \ref{theorem:joint nerual network realization_new}.

    \paragraph{ \it 1. Applying Theorem~\ref{theorem:joint nerual network realization_new}}
    There exist constants $\bar{c}_n,\bar{q}_n,\bar{\tau}_n,\bar{m}_n > 0 ,\; n \in \N^{-1} $, such that for any $\varepsilon >0$ and $K \in \bbN_{+}$, there exist deep ReLU networks $ \tilde{z}^{d,\varepsilon,K}_n :\mathbb{R}^{1+d}\rightarrow \mathbb{R}^d,\; n \in  \N^{-1} $ satisfying
    \[
        \big\|  \tilde{z}^{d,\varepsilon,K}_n - \widehat{Z}^{*,K,n;d} \big\|_{2;\mu^{K;d}_n } \le \frac{1}{2} \varepsilon , \; \; \text{and} \; \; \Growth(\tilde{z}^{d,\varepsilon,K}_n(t,\cdot)) ,\; \size(\tilde{z}^{d,\varepsilon,K}_n ) \le \bar{c}_n d^{\bar{q}_n} \varepsilon^{-\bar{\tau}_n}(K)^{\bar{m}_n} 
    \]
    for all $ t\in [t_n,t_{n+1}] $.
    
    \paragraph{ \it 2. Applying Theorem~\ref{theorem:express_N_0_new}}
    By the structure of the dynamic process $X^d$ and the terminal function $ g^d $, we already have positive constants $B^{*},Q^{*} $. Hence, for the above $\varepsilon$, there exists $ K_{d,\varepsilon} \le B^{*}d^{Q^{*}} \varepsilon^{-1} $ such that for all $ n\in \N^{-1} $,
    \[
        \big\| \widehat{Z}^{*,K_{d,\varepsilon},n ; d} - Z^{*;d} \big\|^2_{[t_n,t_{n+1}],2} \le \frac{1}{2} \varepsilon .
    \]

    \paragraph{ \it 3. Combining these results}
    Let $ \tilde{z}^{d,\varepsilon}_n := \tilde{z}^{d,\varepsilon,K_{d,\varepsilon} }_n $. Then, for any $ n\in \N^{-1} $, $t\in [t_n,t_{n+1}] $, and $x\in \mathbb{R}^d $,
    \begin{align*}
        \| \tilde{z}^{d,\varepsilon}_n(t,x)\| & \le \bar{c}_n d^{\bar{q}_n}\varepsilon^{-\bar{\tau}_n}(B^{*}d^{Q^{*}} \varepsilon^{-1} )^{\bar{m}_n} (1+\|x\|) 
        \le \Tilde{c}_n d^{\Tilde{q}_n} \varepsilon^{\Tilde{r}_n} (1+\|x\|), \; \text{and} \\
        \size( \tilde{z}^{d,\varepsilon}_n )  & \le \bar{c}_n d^{\bar{q}_n} \varepsilon^{-\bar{\tau}_n}(B^{*}d^{Q^{*}} \varepsilon^{-1})^{\bar{m}_n}  \le  \Tilde{c}_n d^{\Tilde{q}_n} \varepsilon^{\Tilde{r}_n} ,
    \end{align*}
    where $\Tilde{c}_n = \bar{c}_n (B^{*})^{\bar{m}_n},\; \Tilde{q}_n = \bar{q}_n + Q^{*}\bar{m}_n $ and $ \Tilde{r}_n = \bar{\tau}_n + \bar{m}_n $. Taking maxima over $n\in \N^{-1}$, we may use the same constants $\Tilde{c},\Tilde{q},\Tilde{r} $ for all $ n \in \N^{-1} $.

    Next, since $Y^{*;d}_n = \widetilde{U}^d_n(M^{*;d}) $ by Lemma~\ref{lemma:sure_optimal}, Lemma~\ref{lem-error-propagate} together with It\^o isometry yields
    \begin{align*}
         \big\| \widetilde{U}^d_n(\widetilde{M}^{d,\varepsilon}  ) - Y^{*;d}_n \big\|_{L^2} & \le   \big\| \widetilde{U}^d_n(\widetilde{M}^{d,\varepsilon}  ) - \widetilde{U}^d_n(\widehat{M}^{*,K_{d,\varepsilon} }) \big\|_{L^2} +  
         \big\| \widetilde{U}^d_n(\widehat{M}^{*,K_{d,\varepsilon}}) - \widetilde{U}^d_n(M^{*;d}) \big\|_{L^2} \\
        & \le \sum_{m=n}^{N-1} \Big\|\sum_{k=0}^{K_{d,\varepsilon}-1} \big( \widehat{Z}^{*,K_{d,\varepsilon},m ; d}(t^m_k,X^d_{t^m_k}) - \tilde{z}^{d,\varepsilon}_m(t^m_k,X^d_{t^m_k}) \big) \cdot \Delta W^d_{t^m_k} \Big\|_{L^2}  \\
        & \quad \quad  + 
        \sum_{m=n}^{N-1}  \Big\|\sum_{k=0}^{K_{d,\varepsilon}-1} \int_{t^m_k}^{t^m_{k+1}} \big( \widehat{Z}^{*,K_{d,\varepsilon},m ; d}(t^m_k,X^d_{t^m_k}) - Z^{*;d}_s \big) \cdot dW^d_s  \Big\|_{L^2}  \\
        & = \sum_{m=n}^{N-1}\Big[ \big\|  \tilde{z}^{d,\varepsilon,K_{d,\varepsilon}}_n - \widehat{Z}^{*,K_{d,\varepsilon},n;d} \big\|_{2;\mu^{K;d}_n } + \big\| \widehat{Z}^{*,K_{d,\varepsilon},n ; d} - Z^{*;d} \big\|^2_{[t_n,t_{n+1}],2}  \Big] \\
        & \le \frac{1}{2}(N-n)\varepsilon + \frac{1}{2}(N-n)\varepsilon  \le (N-n)\varepsilon 
    \end{align*}
    as claimed.
\end{proof}

\subsubsection{Detailed proof of DeepMartingale's expressivity for AID }

We first recall some important propositions of the affine function and AID discussed in \citep{Jentzen23}.
\begin{lemma}
\label{lemma:equiv_affine_new}
$a^d:\mathbb{R}^{d} \to \mathbb{R}^{d} $, $b^d:\mathbb{R} \to \mathbb{R}^{d\times d} $ are affine vector(matrix)-valued functions if and only if there exists $A^{1}_d \in \mathbb{R}^{d\times d},\; b^1_d \in \mathbb{R}^d,\; A^2_d \in \mathbb{R}^{d\times d \times d},\; $ and $ b^2_d \in \mathbb{R}^{d\times d} $, such that for all $z \in \mathbb{R}^d $,
\[ 
     a^d(z) = A^1_d z + b^1_d , \; \text{and} \; \; b^d(z) = ( A^{2;1}_d z , \ldots , A^{2;d}_d z ) + b^2_d
\]
In particular, $A^1_d, b^1_d, A^2_d, $ and $ b^2_d $ have the following forms:
\begin{align*}
    A^1_d & = \big(a^d(e_1)-a^d(0),\ldots,a^d(e_d)-a^d(0) \big), \; b^1_d=a^d(0),\;  \text{and}\\
    A^2_d & =  ( A^{2;1}_d  , \ldots , A^{2;d}_d ) =\big(b^d(e_1)-b(0),\ldots,b^d(e_d)-b^d(0) \big), \; b^2_d =b^d(0),
\end{align*}
where $e_1=(1,0,\ldots,0 )^{\top} ,\; e_2=(0,1,\ldots,0)^{\top}, \ldots ,\; e_d=(0,0,\ldots,1)^{\top} $.

\end{lemma}

According to \citep{Jentzen23}, the following propositions hold.
\begin{proposition}[Existence of a dynamic process with a continuous sample path]
    \label{theorem:continous_sample_path_new}
    For AID (Definition~\ref{def:AID}), there exist up to indistinguishable unique $\mathbb{F}$-adapted 
    stochastic processes with continuous sample paths $X^d :[0,T]\times \Omega \to \mathbb{R}^d $ that satisfies the following:
    for all $x^d_0 \in \mathbb{R}^d$ and $ t\in [0,T]$, SDE~\eqref{eq:SDE} holds $\mathbb{P} $-a.s.. 
\end{proposition}

\begin{proposition}[Linear RanNN representation of AID]
\label{proposition:affine_process_new}
   For any $0\le s\le T $, if $X^{s,x;d}_t, \; $ where $ s\le t\le T $ denotes an AID with continuous sample path starting at $s $ for any initial value $x$, then for any $s\le t\le T$, there exists a random matrix and random vector $A^{s;d}_t:\Omega \to \mathbb{R}^{d\times d},\; b^{s;d}_t:\Omega \to \mathbb{R}^d $, such that
   \begin{equation}\label{eq:affine_linear_expand}
       X^{s,x;d}_{t}(\omega) = A^{s;d}_t(\omega) x + b^{s;d}_t(\omega),\; \forall \; x\in \mathbb{R}^d,\; \omega \in \Omega.
   \end{equation}
\end{proposition}
By direct verification, AID-log satisfies Assumption~\ref{ass:N_0_structural_ass_model_deter}. Thus, according to Theorem~\ref{theorem:dynamic_bound_new}, we have the linear growth rate bound for AID-log with the expression rate. 
\begin{proposition}[Dynamic bound]
    \label{proposition:dynamic_bound_new}
    If $X^d $ is an AID-log, given any $\bar{p} \in [2,\infty) $, there exists positive constants $B_{\bar{p}},Q_{\bar{p}},R_{\bar{p}} $ that are only dependent on $T,\bar{p}$, such that 
    \begin{equation*}
        \mathbb{E}|X^{*,0;d}_T|^{\bar{p}} \le B_{\bar{p}}d^{Q_{\bar{p}}}(\log d)^{R_{\bar{p}}}(1+\|x^d_0\|^{\bar{p}}).
    \end{equation*}
    If $X^{s,l,x;d}_t,\; 0\le s\le t\le l\le  T$ follows the same coefficient functions assumption as AID with $\frac{1}{2}$-$\log $ growth rate, which is that it starts at $s$ with the value $x$, then a similar argument holds for the same constants:
\begin{equation*}
    \mathbb{E}|X^{*,s;d}_l|^{\bar{p}} \le B_{\bar{p}}d^{Q_{\bar{p}}}(\log d)^{R_{\bar{p}}}(1+\|x\|^{\bar{p}}).
\end{equation*}
\end{proposition}

According to Lemma~\ref{lemma:equiv_affine_new}, we can utilize the fundamental matrix of the linear SDE (e.g., \citep{MAO201191}) to further derive the result of AID-log, especially for the Lipschitz bound.  We denote $\|\cdot\|_2 $ as the square norm of the matrix induced by vector ($||A||_2 := \sup_{\| x\|=1} \|Ax\| $) and the fundamental matrix for the homogeneous linear SDE ($b^1_d,b^2_d=0 $ in Lemma~\ref{lemma:equiv_affine_new}) $ $ and $ \Phi^d(t)=\Phi^d(\omega,t) $ (omit $\omega$ for simplicity), which satisfies the following matrix-valued linear SDE ( \citep[Chapter 3]{MAO201191}): $\Phi^d(s)\equiv \mathbf{I}_d $ and
\begin{equation*}
    d\Phi^d(t) = A^1_d\Phi^d(t) dt + \sum_{i=1}^d A^{2;i}_d\Phi^d(t) dW^{i;d}_t ,\; s \le t\le l.
\end{equation*}
The following proposition provides the expressivity result for the AID-log fundamental matrix and subsequently the Lipschitz bound for AID-log.
\begin{proposition}\label{prop:fundamental_solution_bound}
    Under Definition~\ref{def:ADI_express}, for the fundamental matrix $\Phi^{s;d}_t $ on $[s,l],\; s<l \le T $, the following expressivity result holds: given any ${\bar{p}}\ge 2 $, there exists positive constants $B_{\bar{p}},Q_{\bar{p}},R_{\bar{p}} $ that are only dependent on $T,{\bar{p}} $ (chosen to be the same as in Proposition~\ref{proposition:dynamic_bound_new}), such that
    \begin{equation}\label{eq:fundamental_matrix_est}
        \mathbb{E}| \Phi^{*,s;d}_l |^{\bar{p}} \le B_{\bar{p}}d^{Q_{\bar{p}}}(\log d)^{R_{\bar{p}}}.
    \end{equation}
    Then, for the Lipschitz bound of AID-log (which is only determined by $A^s_t$ in Equation~\eqref{eq:affine_linear_expand}), we have
    \begin{equation}\label{eq:weights_lip_bound}
        \mathbb{E}\|A^{s;d}_t\|_2^{\bar{p}} \le B_{\bar{p}}d^{Q_{\bar{p}}}(\log d)^{R_{\bar{p}}}.
    \end{equation}
\end{proposition}

\begin{proof}[Proof of Proposition~\ref{prop:fundamental_solution_bound}]
It is easy to verify that the fundamental matrix $\Phi^{s;d}_t $ satisfies Assumption~\ref{ass:N_0_structural_ass_model}; thus, we can directly apply Theorem~\ref{thm:dynamic_X_multidim_tensor}, which immediately derives \eqref{eq:fundamental_matrix_est} by choosing the constants. For \eqref{eq:weights_lip_bound}, obviously, $ A^{s;d}_t x = X^{s,x;d}_t - X^{s,0;d}_t  $ is the solution of the following linear SDE: $Y^d_s = x $ and
    \begin{equation*}
        dY^d_t = A^1_d Y^d_t dt + \sum_{i=1}^d A^{2;i}_d Y^d_t dW^{i;d}_t ,\; s \le t\le l.
    \end{equation*}
    By \citep[Theorem 2.1, uniqueness, Chapter 3]{MAO201191}, we have \(  A^{s;d}_t x = \Phi^{s;d}_t x ,\; \forall \; x\in \R^d \);
    then, 
    \begin{equation*}
        \mathbb{E}\|A^{s;d}_t\|_2^{\bar{p}} = \mathbb{E}\|\Phi^{s;d}_t\|_{2}^{\bar{p}} \le \mathbb{E}\|\Phi^{s;d}_t\|_{\H}^{\bar{p}} \le  B_{\bar{p}}d^{Q_{\bar{p}}}(\log d)^{R_{\bar{p}}}
    \end{equation*}
    as claims.
\end{proof}

Building on the above preparation, we here prove Lemma~\ref{lemma:AID-log_ass1} and Theorem~\ref{thm:AID-log_express}.

\begin{proof}[Proof of Lemma~\ref{lemma:AID-log_ass1}]
The argument for Assumption~\ref{ass:N_0_structural_ass_model_deter} is trivial. For Assumption~\ref{ass:dynamic_ass} and Assumption~\ref{ass:stronger_ass}, 
we know for any fixed $\omega$ , $(x,t,s)\mapsto X^{s,x ;d}_t(\omega)=:f^{t;d}_s(x,\omega) $ is $\mathcal{B}(\mathbb{R}^{d+2}) $-measurable. By Proposition~\ref{proposition:affine_process_new}, $f^{t;d}_s(x,\omega)=X^{s,x;d}_{t}(\omega) = A^{s;d}_t(\omega) x + b^{s;d}_t(\omega), \; \forall \; x\in \mathbb{R}^d,\; \omega \in \Omega $, which can obviously be represented by a RanNN $\widehat{f}^{t;d}_s = f^{t;d}_s $ with depth $I_s^{t;d} \le 1$ and $\size(\widehat{f}^{t;d}_s(\cdot,\omega)) \le d(d+1), \; \forall \; \omega $. In addition, for any $\tilde{p} \ge 2 $, by Proposition~\ref{proposition:dynamic_bound_new},
\[
     \| f^{t;d}_s(0,\omega) \| = \|X^{s,0;d}_t(\omega)\| = \|b^{s;d}_t(\omega)\| ,\; \;  \| f^{t;d}_s(x,\omega) - f^{t;d}_s(y,\omega) \| = \|A^{s;d}_t(\omega)(x-y)\|  .
\]
    \begin{equation}
        \big\| f^{t;d}_s(0,\cdot) \big\|_{L^{\tilde{p}}}
          = \big\| X^{s,0;d}_t \big\|_{L^{\tilde{p}}} \le (B_{\tilde{p}})^{\frac{1}{\tilde{p}}}d^{\frac{1}{\tilde{p}}(Q_{\tilde{p}}+R_{\tilde{p}})}. \label{eq:initial_bound_AID-type}
    \end{equation}
    For the Lipschitz bound, by Proposition~\ref{prop:fundamental_solution_bound}, for any $\tilde{p}\ge 2 $,
    \[
         \big\| \Lip(f^{t;d}_s(*,\cdot)) \big\|_{L^{\tilde{p}}} =  \big\| \sup_{\|x\|=1}\|A^{s;d}_tx\| \big\|_{L^{\tilde{p}}}  = \big(\mathbb{E} \|A^{s;d}_t \|_2^{\tilde{p}}   \big)^{\frac{1}{\tilde{p}}} \le (B_{\tilde{p}})^{\frac{1}{\tilde{p}}}d^{\frac{1}{\tilde{p}}(Q_{\tilde{p}}+R_{\tilde{p}})}  .
    \]
    Thus, by
    \[
    \frac{\|f^{t;d}_s(x,\omega)  \|}{1+\|x\|} \le \frac{\|f^{t;d}_s(x,\omega) - f^{t;d}_s(0,\omega)  \|}{\|x\|} \frac{\|x\|}{1+\|x\|} + \frac{\|f^{t;d}_s(0,\omega)\|}{1+\|x\|} \le \Lip(f^{t;d}_s(*,\omega)) + \|f^{t;d}_s(0,\omega)\|  ,
    \]
    we know $ \big\| \Growth(f^{t;d}_s(*,\cdot)) \big\|_{L^{\tilde{p}}} \le 2(B_{\tilde{p}})^{\frac{1}{\tilde{p}}}d^{\frac{1}{\tilde{p}}(Q_{\tilde{p}}+R_{\tilde{p}})} $.
    Then, AID-log satisfies Assumption~\ref{ass:dynamic_ass} for any $p>2 $ and Assumption~\ref{ass:stronger_ass} for any $\tilde{p}>4 $. 
\end{proof}

\begin{proof}[Proof of Theorem~\ref{thm:AID-log_express}]
By Lemma~\ref{lemma:AID-log_ass1}, we know that AID-log satisfies Assumption~\ref{ass:N_0_structural_ass_model_deter},  Assumption~\ref{ass:dynamic_ass} for any $p>2 $ and Assumption~\ref{ass:stronger_ass} for any $\tilde{p}>4 $. By Corollary~\ref{coro:g_growth_bound_new}, Corollary~\ref{coro:V_growth_bound_new}, and Proposition~\ref{prop:V_Lip_growth}, we can directly apply Theorem~\ref{thm:express_deep_mtg}, which completes the proof. 
\end{proof}

\end{appendix}

\end{document}